\documentclass[a4paper,12pt]{article}
\usepackage[T2A]{fontenc}
\usepackage[cp1251]{inputenc}
\begin{document}
\author{S.V. Ludkovsky.}
\title{Groups of diffeomorphisms and geometric loops of manifolds
over ultra-normed fields.}
\date{19 December 2007}
\maketitle

\begin{abstract} The article is devoted to the
investigation of groups of diffeomorphisms and loops of manifolds
over ultra-metric fields of zero and positive characteristics.
Different types of topologies are considered on groups of loops and
diffeomorphisms relative to which they are generalized Lie groups or
topological groups. Among such topologies pairwise incomparable are
found as well. Topological perfectness of the diffeomorphism group
relative to certain topologies is studied. There are proved theorems
about projective limit decompositions of these groups and their
compactifications for compact manifolds. Moreover, an existence of
one-parameter local subgroups of diffeomorphism groups is
investigated.

\end{abstract}

\section{Introduction.}
\par Non-archimedean analysis has rather long history, but
it is much less developed in comparison with the classical analysis
over the fields $\bf R$ and $\bf C$. Therefore, the theory of groups
on manifolds over non-archimedean fields is not so well investigated
as for Riemann or complex manifolds \cite{pontr,hew,fell,banyaga}.
\par As it is known fields with multiplicative ultra-norms such as
the field of $p$-adic numbers were first introduced by K. Hensel
\cite{hensel}. Several years later on it was proved by A. Ostrowski
\cite{ostrow} that on the field of rational numbers each
multiplicative norm is either the usual norm as in $\bf R$ or is
equivalent to a non archimedean norm $|x|=p^{-k}$, where
$x=np^k/m\in \bf Q$, $n, m, k\in \bf Z$, $p\ge 2$ is a prime number,
$n$ and $m$ and $p$ are mutually pairwise prime numbers. Each
locally compact infinite field with a non trivial non archimedean
valuation is either a finite algebraic extension of the field of
$p$-adic numbers or is isomorphic to the field ${\bf F}_{p^k}(\theta
)$ of power series of the variable $\theta $ with expansion
coefficients in the finite field ${\bf F}_{p^k}$ of $p^k$ elements,
where $p\ge 2$ is a prime number, $k\in \bf N$ is a natural number
\cite{roo,weil}. The valuation group $\Gamma _{\bf K}:= \{ |x|: x\in
{\bf K}, x\ne 0 \} \subset (0,\infty )$ of a locally compact field
$\bf K$ is discrete. Non locally compact fields are wide spread as
well and among them there are fields with the valuation group
$\Gamma _{\bf K}=(0,\infty )$ \cite{diar,roo,sch1}. The
non-archimedean analog of the field of complex numbers is $\bf C_p$
which is complete algebraically and as the uniform space relative to
its multiplicative norm and $\Gamma _{\bf C_p} = \{ x\in {\bf Q}:
x>0 \} $.
\par The importance of transformation groups of manifolds
in the non-Archimedean functional analysis, representation theory
and mathematical physics is clear and also can be found in the
references given below
\cite{aref,luseamb,lutmf99,luseamb03,luijmms0342,roo,roosch}. This
article is devoted to several aspect of such groups. One of them is
on their structure from the point of view of the non-archimedean
compactficiation (see also about Banaschewski compactification in
\cite{roo}). Though a new topology used for compactification may be
different from the initial topology of a group or even may be
non-comparable, because on the same group it is possible an
existence of several different topologies making it a topological
group. This is useful also for studying their representations as
restrictions of representations of non-archimedean
compactifications, which are constructed below such that they also
are groups.  Apart from previous works
\cite{luseamb,lutmf99,luanmbp, luseamb03,luijmms0342}, where the
characteristic $char ({\bf K})=0$ was zero, in this paper groups on
manifolds over fields with non-zero characteristics also are defined
and investigated. Different types of topologies are considered on
groups of loops and diffeomorphisms relative to which they are
generalized Lie groups or topological groups. Among such topologies
pairwise incomparable are found as well. It is caused by the fact
that repeated application of projective and inductive limits of
topological spaces generate topologies and spaces in general
dependent from an order of taking limits and their types, so that
such topologies may appear incomparable on a subset contained in
these topological spaces. It is proved, that relative to the
$C^{\infty }$ bounded-open topology groups of geometric loops and
groups of diffeomorphisms of manifolds over ultra-normed fields are
the generalized Lie groups.
\par Previously one-parameter subgroups over fields of non-zero
characteristics were not studied. This article contains as well
results on one-parameter subgroups over $\bf K$ with $char ({\bf
K})=p>1$ using its multiplicative subgroup ${\bf K}^* := {\bf
K}\setminus \{ 0 \} $. It is proved below that the diffeomorphism
group of a compact manifold is topologically simple relative to the
$C^{\infty }$ compact-open topology, that develops previous results
\cite{luseamb03}, where topological simplicity and perfectness was
proved over fields of zero characteristic.
\par At first in Section 2 we remind basic facts and notations,
which are given in detail in references
\cite{roo,sch1,luseamb,lutmf99,luanmbp}. A loop group $L_t(M,N)$ is
defined as a quotient space of a family of mappings $f: M\to N$ of
class $C^t$ of one Banach manifold  $M$ into another $N$ over the
same local field $\bf K$ such that $\lim_{z\to s}(\bar \Phi
^mf)(z;h_1,...,h_m;\zeta _1,..., \zeta _m)=0$ or $\lim_{z\to
s}(\Upsilon ^mf)(z^{[m]})=0$ for each $0\le m\le t$, where $M$ and
$N$ are embedded into the corresponding Banach spaces $X$ and $Y$,
$cl(M)=M\cup \{ s \} $, $cl(M)$ and $N$ are clopen in $X$ and $Y$
respectively, $0\in N$, $(\bar \Phi ^mf)(z;h_1,...,h_m;\zeta
_1,..,\zeta _m)$ and $(\Upsilon ^mf)(z^{[m]})$ are continuous
extensions of difference quotients by variables corresponding to $z$
or by all appearing inductively variables over a non-archimedean
field $\bf K$ of zero characteristic or of $char ({\bf K})=p>0$
respectively, $z\in M$, $h_1,...,h_m$ are nonzero vectors in $X$,
$\zeta _1,...,\zeta _m\in \bf K$ such that $z+\zeta _1h_1+...+\zeta
_mh_m\in M$, $z^{[m+1]}:=(z^{[m]},v^{[m]},\zeta _{m+1})$,
$z^1=(z,v^{[0]},\zeta _1)$, $z^{[m]}, v^{[m]}\in {\bf K}^{[m]}$,
${\bf K}^{[m+1]}={\bf K}^{[m]}\oplus {\bf K}^{[m]}\oplus {\bf K}$,
$U^{[1]}=U$, $z^{[m]}+\zeta _{m+1}v^{[m]}\in U^{[m]}$ (see also \S
2.1).
\par In Section 2 preliminary investigations on structures of
$Diff^t(M)$ as the topological groups and Lie groups are studied.
Non-archimedean completions of clopen subgroups $W$ of loop groups
$G$ and diffeomorphism groups $G$ are considered in Sections 3 and
7. Completions are considered relative to uniformities associated
with projective decompositions. They produce topologies incomparable
with the initial one. Relative to them they remain topological
groups. In the case of the loop group the non-archimedean completion
produces a new topological group $V$ in which the initial group $W$
is embedded as a dense subgroup such that $V\ne W$. Such topologies
have purely non-archimedean origin related with non-archimedean
uniformities or families of non-archimedean semi-norms on spaces of
continuous or more narrow classes of functions between
non-archimedean manifolds. In the classical case over $\bf R$ one
might expect instead of this some repeated combination of an
inductive and a projective limits, which is quite different thing.

\par For the compact manifold $M$ in the case of the diffeomorphism
group the non-archimedean completion of $W$ produces profinite
group. For the locally compact manifolds $M$ and $N$ in the case of
the loop group $L_t(M,N)$ one of the non-archimedean completion of
$W$ produces its embedding into ${\bf Z_p}^{\bf N}$ and also there
exists the completion isomorphic with $(\nu {\bf Z})^{\aleph _0}$,
where $\nu {\bf Z}$ is the one-point Alexandroff compactification of
$\bf Z$. When $W$ is bounded relative to the corresponding metric in
$L_t(M,N)$, then $W$ is embedded into ${\bf Z_p}^{\bf N}$. Moreover,
topologies of $Diff^w(M)$ and $Diff^t(M)$ or $L_t(M,N)$ and
$L^w(M,N)$ are incomparable for compact manifolds $M$ and $N$, where
the groups $Diff^w(M)$ and $L^w(M,N)$ are supplied with the weak
projective limit topologies $\tau _w$.
 The group $Diff^t(M)$ is topologically simple, on the other hand,
the group $L_t(M,N)$ is commutative. \par An existence of one
parameter subgroups of $Diff^t(M)$ is investigated in Section 5. It
is proved in Section 6, that $Diff^t(M)$ is topologically simple
relative to its $C^t$ compact-open topology, as well as the theorem
about continuous automorphisms of $Diff^t(M)$ is proved.
\par The notation given below and the corresponding definitions are
given in detail in \cite{luseamb,luanmbp}. All results of this paper
over the fields of positive characteristics are obtained for the
first time.
\section{Groups of diffeomorphisms.}
\par {\bf 1. Definitions.} Let $\bf K$ be an infinite field with a
non trivial non archimedean valuation, let also $X$ and $Y$ be
topological vector spaces over $\bf K$ and $U$ be an open subset in
$X$. For a function $f: U\to Y$ consider the associated function
\par $f^{[1]}(x,v,t) := [f(x+tv) - f(x)]/t$ \\
on a set $U^{[1]}$ at first for $t\ne 0$ such that $U^{[1]} := \{
(x,v,t)\in X^2\times {\bf K}, x\in U, x+tv\in U \} $. If $f$ is
continuous on $U$ and $f^{[1]}$ has a continuous extension on
$U^{[1]}$, then we say, that $f$ is continuously differentiable or
belongs to the class $C^1$. The $\bf K$-linear space of all such
continuously differentiable functions $f$ on $U$ is denoted
$C^{[1]}(U,Y)$. By induction we define functions $f^{[n+1]}:=
(f^{[n]})^{[1]}$ and spaces $C^{[n+1]}(U,Y)$ for $n=1,2,3,...$,
where $f^{[0]}:=f$, $f^{[n+1]}\in C^{[n+1]}(U,Y)$ has as the domain
$U^{[n+1]} := (U^{[n]})^{[1]}$.
\par The differential $df(x): X\to Y$ is defined as
$df(x)v := f^{[1]}(x,v,0)$.
\par Define also partial difference quotient operators $\Phi ^n$
by variables corresponding to $x$ only such that
\par $\Phi ^1f(x;v;t) = f^{[1]}(x,v,t)$ \\
at first for $t\ne 0$ and if $\Phi ^1f$ is continuous for $t\ne 0$
and has a continuous extension on $U^{[1]}=:U^{(1)}$, then we denote
it by ${\bar {\Phi }}^1f(x;v;t)$. Define by induction \par $\Phi
^{n+1} f(x;v_1,...,v_{n+1};t_1,...,t_{n+1}):= \Phi ^1(\Phi
^nf(x;v_1,...,v_n;t_1,...,t_n))(x;v_{n+1};t_{n+1})$ \\
at first for $t_1\ne 0,...,t_{n+1}\ne 0$ on $U^{(n+1)}:= \{
(x;v_1,...,v_{n+1};t_1,...,t_{n+1}): x\in U; v_1,...,v_{n+1}\in X;
t_1,...,t_{n+1}\in {\bf K}; x+v_1t_1\in
U,...,x+v_1t_1+...+v_{n+1}t_{n+1}\in U \} $. If $f$ is continuous on
$U$ and partial difference quotients $\Phi ^1f$,...,$\Phi ^{n+1}f$
has continuous extensions denoted by ${\bar {\Phi }}^1f$,..., ${\bar
{\Phi }}^{n+1}f$ on $U^{(1)}$,...,$U^{(n+1)}$ respectively, then we
say that $f$ is of class of smoothness $C^{n+1}$. The $\bf K$ linear
space of all $C^{n+1}$ functions on $U$ is denoted by
$C^{n+1}(U,Y)$, where $\Phi ^0f := f$, $C^0(U,Y)$ is the space of
all continuous functions $f: U\to Y$. \par Then the differential is
given by the equation $d^nf(x).(v_1,...,v_n) := n! {\bar {\Phi
}}^nf(x;v_1,...,v_n;0,...,0)$, where $n\ge 1$, also denote
$D^nf=d^nf$. Shortly we shall write the argument of $f^{[n]}$ as
$x^{[n]}\in U^{[n]}$ and of ${\bar {\Phi }}^nf$ as $x^{(n)}\in
U^{(n)}$, where $x^{[0]}=x^{(0)}=x$, $x^{[1]}=x^{(1)}=(x,v,t)$,
$v^{[0]}=v^{(0)}=v$, $t_1=t$, $x^{[k]}=(x^{[k-1]},v^{[k-1]},t_k)$
for each $k\ge 1$, $x^{(k)} := (x;v_1,...,v_k;t_1,...,t_k)$.
\par Subspaces of uniformly $C^n$ or $C^{[n]}$ bounded continuous
functions together with ${\bar {\Phi }}^kf$ or $\Upsilon ^kf$ on
bounded open subsets of $U$ and $U^{(k)}$ or $U^{[k]}$ for
$k=1,...,n$ we denote by $C^n_b(U,Y)$ or $C^{[n]}_b(U,Y)$
respectively.
\par Consider partial difference quotients of products and compositions
of functions and relations between partial difference quotients and
differentiability of both types. Denote by $L(X,Y)$ the space of all
continuous $\bf K$-linear mappings $A: X\to Y$. By $L_n(X^{\otimes
n},Y)$ denote the space of all continuous $\bf K$ $n$-linear
mappings $A: X^{\otimes n}\to Y$, particularly,
$L(X,Y)=L_1(X^{\otimes 1},Y)$. If $X$ and $Y$ are normed spaces,
then $L_n(X^{\otimes n},Y)$ is supplied with the operator norm: $ \|
A \| := \sup_{h_1\ne 0,...,h_n\ne 0; h_1,...,h_n\in X} \|
A.(h_1,...,h_n) \|_Y/ (\| h_1 \| _X... \| h_n \| _X)$.
\par {\bf 2. Lemma.} {\it The spaces $C^{[1]}(U,Y)$ and $C^1(U,Y)$
are linearly topologically isomorphic. If $f\in C^n(U,Y)$, then
${\bar {\Phi }}^nf(x;*;0,...,0): X^{\otimes n}\to Y$ is a $\bf K$
$n$-linear $C^0(U,L_n(X^{\otimes n},Y))$ symmetric map.}
\par {\bf Proof.} From Definition 1 it follows, that
$f^{[1]}(x,v,t)= {\bar {\Phi }}^1f(x;v;t)$ on $U^{[1]}=U^{(1)}$, so
both $\bf K$-linear spaces are linearly topologically isomorphic. On
the other hand, due to its definition ${\bar {\Phi
}}^nf(x;*,0,...,0)$ is the $\bf K$ $n$-linear symmetric mapping for
each $x\in U$ and it belongs to $C^0(U,L_n(X^{\otimes n},Y))$, since
${\bar {\Phi }}^nf(x;v_1,...,v_n;t_1,...,t_n)$ is continuous on
$U^{(n)}$ and for each $x\in U$ and $v_1,...,v_n\in X$ there exist
neighborhoods $V_i$ of $v_i$ in $X$ and $W$ of zero in $\bf K$ such
that $x+WV_1+...+WV_n\subset U$.
\par {\bf 3. Lemma.} {\it Operators $\Upsilon ^n(f) := f^{[n]}$
from $C^{[n]}(U,Y)$ into $C^0(U^{[n]},Y)$ and ${\bar {\Phi }}^n:
C^n(U,Y)\to C^0(U^{(n)},Y)$ are $\bf K$-linear and continuous.}
\par {\bf Proof.} Since $[(af+bg)(x+vt)-(af+bg)(x)]/t=
a(f(x+vt)-f(x))/t + b(g(x+vt)-g(x))/t$ for each $f, g\in C^1(U,Y)$
and each $a, b\in \bf K$, then applying this formula by induction
and using definitions of operators $\Upsilon ^n$ and ${\bar {\Phi
}}^n$ we get their $\bf K$-linearity. Indeed, \par $\Upsilon
^n(af+bg)(x^{[n]})=\Upsilon ^1(\Upsilon
^{n-1}(af+bg)(x^{[n-1]}))(x^{[n]})=\Upsilon
^1(af^{[n-1]}+bg^{[n-1]})(x^{[n]})=af^{[n]}(x^{[n]})+bg^{[n]}(x^{[n]})$
and
\par ${\bar {\Phi }}^n(af+bg)(x^{(n)})={\bar {\Phi }}^1(
{\bar {\Phi }}^{n-1}(af+bg)(x^{(n-1)}))(x^{(n)})={\bar {\Phi
}}^1(af^{(n-1)}+bg^{(n-1)})(x^{(n)})=af^{(n)}(x^{(n)})+bg^{(n)}(x^{(n)})$.
\\ The continuity of $\Upsilon ^n$ and ${\bar {\Phi }}^n$ follows
from definitions of spaces $C^{[n]}(U,Y)$ and $C^n(U,Y)$
respectively.
\par {\bf 4. Definitions.} Let $M$ be a manifold modelled on a
topological vector space $X$ over $\bf K$ such that its atlas $At
(M) := \{ (U_j,\mbox{ }_M\phi _j): j\in \Lambda _M \} $ is of class
$C^{\alpha ' }_{\beta '}$, that is the following four conditions are
satisfied:
\par $(M1)$ $\{ U_j: j\in \Lambda _M \} $ is an open covering of $M$,
$U_j=\mbox{ }_MU_j$, \par $(M2)$ $\bigcup_{j\in \Lambda _M} U_j=M$,
\par $(M3)$ $\mbox{ }_M\phi _j :=\phi _j:
U_j\to \phi _j(U_j)$ is a homeomorphism for each $j\in \Lambda _M$,
$\phi _j(U_j)\subset X$,
\par $(M4)$ $\phi _j\circ \phi _i^{-1}\in C^{\alpha '}_{\beta '}$ on
its domain for each $U_i\cap
U_j\ne \emptyset $, \\
where $\Lambda _M$ is a set, $C^{\infty }:= \bigcap_{l=1}^{\infty
}C^l_{\beta }$, $C^{[\infty ]}_{\beta }:=\bigcap_{l=1}^{\infty
}C^{[l]}$, $\alpha '\in \{ n, [n]: 1\le n \le \infty \} $, $\beta
\in \{ \emptyset , b \} $, $C^{\alpha '}_{ \emptyset } := C^{\alpha
'}$.
\par Supply $C^{\alpha }_{\beta }(U,Y)$ with the bounded-open
$C^{\alpha }_{\beta }$ topology (denoted by $\tau _{\alpha, \beta }$
generally or $\tau _{\alpha }$ for $\beta =\emptyset $ or for
compact $U$) with the base $W(P,V)= \{ f\in C^{\alpha }_{\beta
}(X,Y): S^kf|_P\in V, k=0,...,n \} $ of neighborhoods of zero, where
$P$ is bounded and open in $U\subset X$, $P\subset U$, $V$ is open
in $Y$, $0\in V$, $S^k={\bar {\Phi }}^k$ or $S^k=\Upsilon ^k$ for
$\alpha =n$ or $\alpha =[n]$ respectively, $v_1,...,v_n\in (P-y_0)$,
$v^{[k]}_l\in (P-y_0)$ for each $k, l$ for some marked $y_0\in P$
and $|t_j|\le 1$ for every $j$.
\par If $M$ and $N$ are $C^{\alpha '}_{\beta }$ manifolds
on topological vector spaces $X$ and $Y$ over $\bf K$, then consider
the uniform space $C^{\alpha }_{\beta }(M,N)$ of all mappings $f:
M\to N$ such that $f_{j,i}\in C^{\alpha }_{\beta }$ on its domain
for each $j\in \Lambda _N$, $i\in \Lambda _M$, where $f_{j,i}:=
\mbox{ }_N\phi _j\circ f\circ \mbox{ }_M\phi _i^{-1}$ is with values
in $Y$, $\alpha \le \alpha '$. The uniformity in $C^{\alpha }_{\beta
}(M,N)$ is inherited from the uniformity in $C^{\alpha }_{\beta
}(X,Y)$ with the help of charts of atlases of $M$ and $N$. If $M$ is
compact, then $C^{\alpha }_b(M,N)$ and $C^{\alpha }(M,N)$ coincide.
\par The family of all homeomorphisms $f: M\to M$ of class
$C^{\alpha }_{\beta }$ denote by $Diff^{\alpha }_{\beta }(M)$.
\par Let $\gamma $ be a set, then denote by $c_0(\gamma ,{\bf K})$
the normed space consisting of all vectors $ x= \{ x_j\in {\bf K}:
j\in \gamma , \mbox{ for each } \epsilon >0 \mbox{ the set } \{ j:
|x_j|>\epsilon \} \mbox{ is finite } \} $, where $\| x \| :=
\sup_{j\in \gamma } |x_j|$. In view of the Kuratowski-Zorn lemma it
is convenient to consider $\gamma $ as an ordinal. Henceforth,
suppose that $X=c_0(\gamma _X,{\bf K})$ and $Y=c_0(\gamma _Y,{\bf
K})$.
\par {\bf 5. Theorem.} {\it The uniform space $Diff^{\alpha
}_b(M)$ (see \S 4) is the topological group relative to compositions
of mappings.}
\par {\bf Proof.} The group operation in $Diff^{\alpha }_{\beta
}(M)$ is $(f,g)\mapsto f\circ g$, where $f\circ g(x):=f(g(x))$ for
each $x\in M$. Then $f=id$ is the unit element in $Diff^{\alpha
}_{\beta }(M)$, where $id(x)=x$ for each $x\in M$. Since the
composition of mappings is associative, then $f\circ (g\circ h)=
(f\circ g)\circ h$ is associative as the group operation. For each
$f\in Diff^{\alpha }_{\beta }$ there exists its inverse mapping
$f^{-1}$ such that $f^{-1}(y)=x$ for each $y=f(x)$, $x\in M$, since
$f: M\to M$ is the homeomorphism. It remains to verify that
$f^{-1}\in Diff^{\alpha }_{\beta }(M)$ for each $f\in Diff^{\alpha
}_{\beta }(M)$ and the composition $(Diff^{\alpha }_{\beta })^2\ni
(f,g)\mapsto f\circ g\in Diff^{\alpha }_{\beta }(M)$ and inversion
$f\mapsto f^{-1}$ are continuous operations.
\par In the normed space $Y=c_0(\gamma _Y,{\bf K})$ a subspace
$span_{\bf K} \{ e_j: j\in \gamma _X \} $ consisting of all finite
$\bf K$-linear combinations of vectors $e_j=(0,...,0,1,0,...)$ with
$1$ on the $j$-th place is everywhere dense. Therefore, each $f\in
C^{\alpha }_b(U,Y)$ is the uniform limit of mappings
$(f_1,...,f_j,0,...)\in C^{\alpha }_b(U,Y)$ together with ${\bar
{\Phi }}^kf$ or $\Upsilon ^kf$ on bounded subsets of $U$ and
$U^{(k)}$ or $U^{[k]}$ for $1\le k\le n$, $n\in \bf N$, $n\le \alpha
$. In particular, consider $\mbox{ }_N\phi _l\circ f\circ \mbox{
}_M\phi _i^{-1}$ for $f\in C^{\alpha }_b(M,N)$ taking $U$ as a
finite union of $\mbox{ }_M\phi _i(\mbox{ }_MU_j)$. Consider all
possible embeddings of ${\bf K}^v$ into $X$, particularly,
containing $x^{(k)}$ or $x^{[k]}$ for each $0\le k\le n$, where
$n\in \bf N$, $n\le \alpha $. In view of Formulas 6$(1)$ or 7$(1)$
of the Appendix using restrictions on different embedded subspaces
${\bf K}^{[k]}$ or ${\bf K}^{(k)}$ into $X^{[k]}$ or $X^{(k)}$ and
uniform continuity of $\Upsilon ^k$ and ${\bar {\Phi }}^k$ on
bounded open subsets, $k=0,1,2,...$, we get, that $f\circ g\in
C^{\alpha }_b(M,M)$ for each $f, g\in Diff^{\alpha }_b(M)$, since
$f_{l,s}\circ g_{s,i}\in C^{\alpha }_b(U_{l,s,i},X)$ on
corresponding domains $U_{l,s,i}$ in $X$. From $f, g\in Hom(M,M)$ it
follows, that $f\circ g\in Hom (M,M)$, hence $f\circ g\in
Diff^{\alpha }_b(M)$.  Applying to $id_{l,i}=f^{-1}_{l,s}\circ
f_{s,i}$ on corresponding domains Formulas A.6$(1)$ or A.7$(1)$ and
restricting on different embedded subspaces ${\bf K}^{[k]}$ or ${\bf
K}^{(k)}$ in $X^{[k]}$ or $X^{(k)}$ and using uniform continuity on
bounded open subsets for $\Upsilon ^k$ or ${\bar {\Phi }}^k$,
$k=0,1,2,...$, to both sides of this equality gives that $f^{-1}\in
Diff^{\alpha }_b(M)$ for each $f\in Diff^{\alpha }_b(M)$.
\par The space $C^{\phi }_b(U,Y)$ is normed for $U$ bounded in
$X$ for $\phi \in \{ n, [n] \} $ with $n\in \bf N$ such that
\par $(1)$ $\| f \|_{C^n_b(U,Y)} := \sup_{0\le k\le n;
z\in V^{(k)}} \| {\bar {\Phi }}^kf(z) \| _Y$ or
\par $(2)$ $\| f \|_{C^{[n]}_b(U,Y)} :=
\sup_{0\le k\le n; z\in V^{[k]}} \| \Upsilon ^kf(z) \| _Y$, \\
where $V^{(k)}:= \{ z\in U^{(k)}: z=(x;v_1,v_2,...;t_1,t_2,...), \|
v_j \|_X=1 \forall j \} $, $V^{[k]} := \{ z=x^{[k]}\in U^{[k]}: \|
v^{[k-1]}_1 \|_X=1, |\mbox{ }_lv^{[q]}_2t_{q+1}|\le 1,
|v^{[q]}_3|\le 1 \quad \forall l, q \} $. The uniformity of
$C^{\infty }_b(U,Y)$ or $C^{[\infty ]}_b(U,Y)$ is defined by the
family of such norms. Then the uniformity in $C^{\alpha }_b(M,N)$ is
induced by the uniformity in $C^{\alpha }_b(U,Y)$ by all bounded
subsets $U$ in finite unions of $\mbox{ }_M\phi _i(\mbox{ }_MU_i)$,
since to each $f\in C^{\alpha }_b(M,N)$ there corresponds $f_{j,i}=
\mbox{ }_N\phi _j\circ f\circ \mbox{ }_M\phi _i^{-1}$ in $C^{\alpha
}_b(U_{j,i},Y)$ with a corresponding domain $U_{j,i}\subset X$.
 Then application of Formulas A.6$(1)$ or A.7$(1)$ by induction on
$k$ and restricting on different embedded subspaces ${\bf K}^{[k]}$
or ${\bf K}^{(k)}$ in $X^{[k]}$ or $X^{(k)}$ and using uniform
continuity on bounded open subsets gives that $(f,g)\mapsto
f^{-1}\circ g$ is $C^{\alpha }_b$ uniformly continuous on bounded
subsets of $U$ and $U^{(k)}$ or $U^{[k]}$, where $U$ is a finite
union of charts $U_j$ of $M$.
\par {\bf 6. Definition.} A topological group $G$ is called a
$C^{\alpha }_{\beta }$ Lie group if and only if $G$ has a structure
of a $C^{\alpha '}_{\beta }$ manifold and the mapping $G^2\ni (f,g)
\mapsto f^{-1}g\in G$ is of class of smoothness $C^{\alpha }_{\beta
}$, where $\alpha \le \alpha '$.
\par {\bf 7. Theorem.} {\it If $M$ is a $C^{\alpha }_b$
manifold on $X=c_0(\gamma _X,{\bf K})$, where either $\alpha =\infty
$ or $\alpha =[\infty ]$, then $Diff^{\alpha }_b(M)$ is the
$C^{\alpha }_b$ Lie group.}
\par {\bf Proof.} In view of Theorem 5 it remains to demonstrate that
$G$ can be supplied with a structure of $C^{\alpha }_b$ manifold and
that $G^2\ni (f,g) \mapsto f^{-1}g\in G$ is of class of smoothness
$C^{\alpha }_b$. It is possible to take an equivalent atlas of $M$
consisting of clopen (closed and open simultaneously) charts $U_j$
shrinking it a little in case of necessity. Take a base $\cal W$ of
neighborhoods of $id$ in $Diff^{\alpha }_b(M)$ from the proof of
Theorem 5. Then consider a subgroup $\Omega $ in $Diff^{\alpha
}_b(M)$ such that $\Omega W$ covers $Diff^{\alpha }_b(M)$ for each
$W\in \cal W$, where $\Omega W:= \{ gW: g\in \Omega \} $, $gW:=\{
gf: f\in W \} $. Therefore, the base $\Omega {\cal W}$ generates a
topology in $Diff^{\alpha }_b(M)$ equivalent with the initial one
(see Chapter 8 in \cite{eng}). Moreover, $\Omega {\cal W}$ generates
a left uniformity in the group of diffeomorphisms. \par  Consider a
subset $W_U := \{ f|_U\in W: f\in Diff^{\alpha }_b(M) \} $ for a
bounded open $U$ in $U_j$ for some $j\in \Lambda _M$. Then $\phi
_i\circ (gW_U)\circ \phi _j^{-1} \subset C^{\alpha }_b(V_j,X)$ and
$\phi _i\circ (gW_U)\subset C^{\alpha }_b(U,X)$ for each $i$, where
$V_j := \phi _j(U_j)\subset X$. Therefore, $gW_U\cap hW_V= \{ f\in
Diff^{\alpha }_b(M): g^{-1}f|_U\in W, h^{-1}f|_V\in W \} $, $W_U\cap
W_V = \{ f\in Diff^{\alpha }_b(M): f|_{U\cup V}\in W \} $. Put $\psi
_{i,g,j} := \phi _i\otimes g^{-1}\otimes \phi _j^{-1}$ such that
$\psi _{i,g,j}(gW_U):= \phi _i\circ g^{-1}((gW_U)\circ \phi
_j^{-1})\subset C^{\alpha }_b(A,X)$, where $A=\phi _j(U)$, $\psi
^{-1}_{l,h,b}:=\phi _l^{-1}\otimes h\otimes \phi _b$, hence $\psi
_{l,h,b}^{-1}\circ \psi _{i,g,j}(gW_U\cap hW_V) =(\phi _l^{-1}\circ
\phi _i)\circ h (W_U\cap g^{-1}hW_V)\circ (\phi _b^{-1}\circ \phi
_j)^{-1}$ for $U\subset U_j$, $V\subset U_b$ with $U_j\cap U_b\ne
\emptyset $ and $U_l\cap U_i\ne \emptyset $. Thus $\psi
_{i,g,j}\circ \psi ^{-1}_{l,h,b}$ gives the transition mapping for
$Diff^{\alpha }_b(M)$. On the other hand, each $C^{\alpha }_b(A,X)$
has the natural embedding into $C^{\alpha }_b(X,X)$, since $X$ is
totally disconnected and $A$ can be taken clopen in $X$ such that
each $f\in C^{\alpha }_b(A,X)$ has a $C^{\alpha }_b(X,X)$ extension.
For $U_j\cap U_b\ne \emptyset $ the intersection $W_U\cap W_V$ is
non void.
\par  Take $At (Diff^{\alpha }_b(M)):= \{ (W_{g,U}, \psi _{i,g,j}):
g\in \Omega , i, j\in \Lambda _M \} $ with charts $W_{g,U}:=gW_U$,
$W\in \cal W$ and $g\in \Omega $, $U$ bounded in some $U_j$ and with
transition mappings $\psi _{i,g,j}\circ \psi _{l,h,b}^{-1}$ for
$W_{h,U}$ and $W_{g,V}$ when $U\cap V\ne \emptyset $. Since $\psi
_{i,g,j}\circ \psi _{l,h,b}^{-1}\in C^{\alpha }_b$, then this is the
$C^{\alpha }_b$ atlas. The mapping $(f,g)\mapsto f^{-1}g$ is of
class $C^{\alpha }_b(Diff^{\alpha }_b(M),Diff^{\alpha }_b(M))$ due
to Formulas A.6$(1)$ and A.7$(1)$, since the mappings
$(f_{i,j},g_{i,l})\mapsto f_{i,j}^{-1}\circ g_{i,l}$ from $C^{\alpha
}_b(U,X)$ into $C^{\alpha }_b(V,X)$ with $U\subset \phi _l(U_l)\ne
\emptyset $ and $V\subset \phi _j(U_j)$ are of class $C^{\alpha
}_b$.
\par {\bf 8. Theorem.} {\it If an ultrametric field $\bf K$ is complete
relative to its multiplicative norm and $\alpha \ge 1$, then
$Diff^{\alpha }_b(M)$ is complete as a left uniform space.}
\par {\bf Proof.} Recall some facts about uniform spaces.
A subset $A$ of the product $S\times S$ of a set $S$ is called a
relation in $S$. The relation inverse to $A$ is denoted $-A$ such
that $-A= \{ (x,y): (y,x)\in A \} $ and the composition of relations
is denoted $A+B$ such that $A+B = \{ (x,z): \mbox{ there exists a }
y\in S\mbox{ such that }(x,y)\in A\mbox{ and }(y,z)\in B \} $.
Denote by $\Delta := \{ (x,x): x\in S \} $ the diagonal of the
product $S\times S$. Every subset in $S\times S$ containing $\Delta
$ is called an entourage of the diagonal $\Delta $. The family of
all entourages of the diagonal is denoted by ${\cal D}_S$. One
writes $|x-y|<V$ if $(x,y)\in V$ and one says that $x$ and $y$ are
at a distance less than $V$. If the condition $|x-y|<V$ is not
satisfied, then one writes $|x-y|\ge  V$. If $A\subset S$ and
$|x-y|<V$ for each $x, y\in A$, then one says that the diameter
$\delta (A)$ of $A$ is less than $V$. Denote $1A:=A$,
$nA:=(n-1)A+A$.

A uniformity $\cal U$ in a set $S$ is a non-empty subfamily in
${\cal D}_S$ satisfying the following four conditions: \par $(U1)$
If $V\in \cal U$ and $V\subset W\in {\cal D}_S$, then $W\in \cal U$;
\par $(U2)$ If $V_1, V_2\in \cal U$, then $V_1\cap V_2\in \cal U$;
\par $(U3)$ For each $V\in \cal U$ there exists $W\in \cal U$ such
that $2W\subset V$; \par $(U4)$ $\bigcap_{V\in \cal U} V=\Delta $.
\par A topological space $S$ is called a $T_1$ space if for each
$x\ne y\in S$ there exists an open set $U$ in $S$ such that $x\in U$
and $y\notin U$. A topological space $S$ is called a Tychonoff space
and it is denoted by $T_{3\frac{1}{2}}$ if it is a $T_1$ space and
for each point $x\in S$ and each closed set $J$ in $S$ with $x\notin
J$ there exists a continuous function $f: S\to [0,1]\subset \bf R$
such that $f(x)=0$ and $f(y)=1$ for each $y\in J$. If $S$ is a
Tychonoff space then for every finite family of functions
$f_1,...,f_n\in C^0(S,{\bf R})$ or $C^0_b(S,{\bf R})$ the formula
$\rho _{f_1,...,f_n}(x,y) := \max_{i=1}^n |f_i(x)-f_i(y)|$ defines a
pseudo-metric in $S$. The families of such pseudo-metrics denote by
$P$ and $P_b$ respectively. They generate uniformities denoted by
$\cal C$ and ${\cal C}_b$ correspondingly (see Chapter 8 in
\cite{eng}).
\par If $\cal U$ is a uniformity in $S$, then the family
${\cal O} := \{ U\subset S: \forall x\in U \exists V\in {\cal U}
\mbox{ such that } B(x,V)\subset U \} $ is the topology in $S$ and
it is called the topology induced by the uniformity $\cal U$.
\par Every covering of $S$ for which there exists $V\in \cal U$
such that ${\cal C}(V)$ is a refinement of it is called a uniform
covering relative to $\cal U$, where ${\cal C}(V):= \{ B(x,V): x\in
S \} $, $B(x,V):= \{ y\in S: |x-y|<V \} $.
\par If $(S,{\cal U})$ is a uniform space and $\cal F$ is a family
of subsets of $S$, then one says that $\cal F$ contains arbitrary
small sets if for every $V\in \cal U$ there exists $F\in \cal F$
such that $\delta (F)<V$. A uniformity $\cal U$ in $S$ is called
complete or $(S,{\cal U})$ is called complete, if for each family of
closed subsets $\{ F_u: u\in \Psi \} $ of the topological space $S$
with the topology induced by $\cal U$ which has the finite
intersection property and contains arbitrary small sets the
intersection $\bigcap_{u\in \Psi } F_u$ is non-empty.
\par If $C$ is an ordered set and $\{ x_u: u\in C \} $ is a
net in $S$ such that for each $V\in {\cal U}$ there exists $u_0\in
C$ for which $|x_u-x_v|<V$ for each $C\ni u, v\ge u_0$, then it is
called the Cauchy net. In accordance with theorem of Chapter 8
\cite{eng} a uniform space $(S,{\cal U})$ is complete if and only if
each Cauchy net in this space is convergent. On the other hand, each
closed subset $A$ of a complete uniform space $(S,{\cal U})$ is
complete relative to the uniformity ${\cal U}_A$ in $A$ inherited
from $\cal U$ in $S$.
\par If $G$ is a topological group and ${\cal B}(e)=\cal B$ is a
base of neighborhoods at the unit element $e$, then each $F\in \cal
B$ determines three coverings of the topological space $G$: ${\cal
C}_l(F)= \{ xF: x\in G \} $, ${\cal C}_r(F):= \{ Fx: x\in G \} $,
${\cal C}(F):= \{ xFy: x, y \in G \} $. By ${\sf C}_l$, ${\sf C}_r$
and $\sf C$ are denoted the families of those coverings of $G$ which
have refinements of the form ${\cal C}_l(F)$, ${\cal C}_r(F)$ or
${\cal C}(F)$ respectively. Each of these families generates a
uniformity in $G$. The topology of each of these uniformities is the
same as the initial topology in $G$.

\par The uniform space $C^{\alpha }_b(M,M)$ is complete for complete $M$.
If $U$ is a bounded canonical closed subset of $M$ contained in a
finite number of charts of $M$, then $Diff^{\alpha }_b(U)$ for each
$\alpha \ge 1$ is the neighborhood of $id|_U$ in $C^{\alpha
}_b(U,X)$, since from $f\in C^{\alpha }_b(U,X)$ with $ \|
(id_{i,j}-f_{i,j})|_{\phi _j(U)} \|_{n,U,X}\le |\pi | $ for each
$i,j $ and $1\le n\le \alpha $, $n\in \bf N$, it follows that $f\in
Diff^{\alpha }_b(U)$, where $\pi \in \bf K$, $0<|\pi |<1$, $\|
*\|_{n,U,X}$ is either $\| * \|_{C^n_b(U,X)}$ or $ \| *
\|_{C^{[n]}_b(U,X)}$ for $\alpha \in {\bf N}\cup \{ \infty \} $ or
$\alpha \in \{ [1],[2],... \} \cup \{ [\infty ] \} $ respectively.

\par Let $\{ f_w: w\in C \} $ be a Cauchy net in $Diff^{\alpha
}_b(M)$, where $C$ is a directed set. This means that for each
neighborhood $W$ of $id$ in $Diff^{\alpha }_b(M)$ there exists
$w_0\in C$ such that $f_w^{-1}f_v\in W$ for each $w_0\le w, v\in C$,
hence $f_v\in f_{w_0}W$ for each $v\ge w_0$. Since $Diff^{\alpha
}_b(M)\subset C^{\alpha }_b(M,M)$ and $C^{\alpha }_b(M,M)$ is
complete, then $f_v$ converges in $C^{\alpha }_b(M,M)$ to a function
$f$. Uniformities in the latter space and in the group of
diffeomorphisms are related by formulas outlined above. We can take
as $W$ a canonical closed subset in $Diff^{\alpha }_b(M)$, $id \in
W$. For each bounded $U$ in $M$ as above the restriction $f|_U$ is
the diffeomorphism of $U$ onto $f(U)$, $f|_U\in f_{w_0}W_U$.
\par Relative to the left uniformity in $Diff^{\alpha }_b(M)$
the left shift mapping $L_h$ is uniformly continuous, where $L_hg :=
hg$ for each $h, g$ in the group. There exists a canonical closed
neighborhood $U$ of the unit element $id$ in $C^{\alpha }_b(M,X)$
such that it is contained in the group, since $\alpha \ge 1$ and due
to definitions of uniformities in them, while $M$ is modeled on the
Banach space. Take a canonical closed symmetric neighborhood $V$ in
$Diff^{\alpha }_b(M)$, $V= cl (Int (V))$, $V^{-1}=V$, such that
$V^2\subset U$. If $ \{ g_n: n \} $ is a Cauchy net in the group,
then for each $\epsilon
>0$ there exists $n_0$ such that $\rho (g_m,g_n)<\epsilon $ for all
$n, m> n_0$, where $\rho $ is the left-invariant metric in the
group. This is equivalent with the fact that for each $\epsilon
>0$ there exists $n_0$ so that $\| g_n^{-1}g_m - id \| _{C^{\alpha}_b}
<\epsilon $ in the corresponding space for each $n, m> n_0$. On the
other hand, there exists $h$ in the group and $N>0$ so that $L_hg_m
\in V$ for every $m>N$. If $X$ is a Banach space, then from the
completeness of $C^{\alpha }_b(M,X)$, and hence of $U$, since $U$ is
closed in $C^{\alpha }_b(M,X)$, it follows, that $Diff^{\alpha
}_b(M)$ is also complete as the uniform space (see also Theorems
8.3.6 and 8.3.20 \cite{eng}). 
\par Since $\bf K$ is complete, then $X$ is complete, hence $X$ is the
Banach space. Thus $C^{\alpha }_b(U,X)$ is complete for each $U$
bounded clopen subset in $X$. Consider a neighborhood $W$ from the
base of neighborhoods of $id$, $W=W_{U,n,\epsilon ,i,j} = \{ f\in
Diff^{\alpha }_b(M): \| f_{i,j}\| _{n,U,X}<\epsilon  \} $, where $i,
j\in \Lambda _M$, $n\in \bf N$, $n\le \alpha $, $U$ is a clopen
bounded subset in $X$ such that $U\subset \phi _j(U_j)$, $0<\epsilon
<\infty $, $\| * \| _{n,U,X}:=\| *\|_{C^n_b(U,X)}$ or $\|
* \| _{n,U,X}:=\| *\|_{C^{[n]}_b(U,X)}$ respectively.
Take without loss of generality $0<\epsilon \le |\pi |$. The
manifold $M$ is on the normed space $X$, hence $M$ is paracompact
and it has a locally finite refinement of its atlas (see
\cite{eng}). Each $\psi _i\circ f_w\circ \phi _j^{-1}$ is converging
to a function denoted $f_{i,j}$ and by transfinite induction this
consistent family $\{ f_{i,j} \} $ induces $f\in C^{\alpha }_b(M,M)$
such that $f_w$ converges to $f$. Since $\Upsilon ^k\phi _i\circ
f_w\circ \phi _j^{-1}(V^{[k]})$ is bounded for each $V^{[k]}$
corresponding to bounded $U$, then $\Upsilon ^k\phi _i\circ f\circ
\phi _j^{-1}(V^{[k]})$ is bounded and $\| f_{i,j}\| _{n,U,X}<\infty
$, analogously for ${\bar {\Phi }}^kf$. Thus $f\in C^{\alpha
}_b(M,M)$, consequently, $f\in Diff^{\alpha }_b(M)$.
\par May be simply it can be proved using the tangent bundle
of the group of diffeomorphisms for $\alpha \ge 1$.
\par {\bf 9. Theorem.} {\it If a manifold $M$ has a finite atlas
having clopen bounded $\phi _j(U_j)\subset X$, $0\le \alpha <\infty
$, then $Diff^{\alpha }_b(M)$ is metrizable by a left invariant
metric.}
\par {\bf Proof.} The metrization Theorem 8.3 \cite{hew}
states that if $G$ is a $T_0$ topological group, then $G$ is
metrizable if and only if there is a countable open basis at $e$. In
this case, the metric can be taken left-invariant. If $At(M)$ is
finite such that each $\phi _j(U_j)$ is a clopen bounded subset in
$X$, then the base of neighborhoods of $id$ in $Diff^{\alpha }_b(M)$
is countable and $Diff^{\alpha }_b(M)$ is metrizable by a left
invariant metric in accordance with the general metrization theorem.
Practically take as the metric $\rho (f,g):=\rho (id,f^{-1}g) :=
\max_{i,j\in \Lambda _M} \| id_{i,j}- (f^{-1}g)_{i,j} \| _{C^{\alpha
}_b(V_j,X)}$, where $V_j:=\phi _j(U_j)$.
\par {\bf 10. Remark.} As it is known left and right uniformities
in a topological group may be different. Here a left uniformity was
considered above. On the same group there may exist several
topologies supplying it with structures of a topological group and
these topologies need not be comparable.

\section{Projective decomposition of diffeomorphism groups.}

\par {\bf 1. Notations and Notes.} Let $M$ and $N$ be compact manifolds
over a locally compact field $\bf K$. Suppose that $M$ and $N$ are
embedded into $B({\bf K^{\sf m}},0,1)$ and $B({\bf K^{\sf n}},0,1)$
as clopen (closed and open at the same time) subsets
\cite{boum,luum985}, where ${\sf m}, {\sf n}\in \bf N$, $B(X,y,r):=
\{ z: z\in X; d_X(y,z)\le r \} $ denotes a clopen ball in a space
$X$ with an ultra-metric $d_X$.
\par The unit ball $B({\bf K^{\sf n}},0,1)$ has the ring structure with
coordinate wise addition and multiplication, where $char ({\bf
K})=0$ or $char ({\bf K})=p>1$ is a prime number. This ring is
isomorphic to a subring of diagonal matrices in the ring $M_n({\bf
K})$ of $n\times n$ square matrices over $\bf K$. Then $B({\bf
K}^n,0,|\pi |^k)$ for $k\ge 1$, $\pi \in \bf K$, $0<|\pi |<1$ is its
two-sided ideal, since $\bf K$ is commutative and $|*|=|*|_{\bf K}$
is the multiplicative norm in $\bf K$. Thus there exists the
quotient ring $B({\bf K}^n,0,1)/B({\bf K}^n,0,|\pi |^k)$
\cite{bacht}. The ring $B({\bf K^{\sf n}},0,1)$ is algebraically
isomorphic with the projective limit $B({\bf K^{\sf
n}},0,1)=pr-\lim_k{\bf S_{p^k}}^{\sf n}$, $\bf S_{p^k}$ is a finite
ring consisting of $p^{kc}$ elements such that ${\bf S_{p^k}}={\bf
S_{p^k}}({\bf K})$ is equal to the quotient ring $B({\bf
K},0,1)/B({\bf K},0,p^{-k})$, ${\bf S_{p^k}}^{\sf n}={\bf
S_{p^k}}^{\otimes \sf n}$ is an external product of $\sf n$ copies
of ${\bf S_{p^k}}$, $c$ is a natural number. Though their structure
depends on $char ({\bf K})$ we denote these rings by the same symbol
depending on $\bf K$ and omitting it, when a field is specified. In
particular $B({\bf F}_{p^n}(\theta ),0,1)/B({\bf F}_{p^n}(\theta
),0,p^{-k})=({\bf F}_{p^n})^{\otimes p^k}={\bf S_{p^{nk}}}={\bf
S_{p^{nk}}}({\bf F}_{p^n}(\theta ))$ and $B({\bf Q_p},0,1)/B({\bf
Q_p},0,p^{-k}) ={\bf Z_p}/(p^k{\bf Z_p})={\bf S_{p^k}}={\bf
S_{p^k}}({\bf Q_p})$ are finite rings consisting of $p^{nk}$ and
$p^k$ elements respectively, $\bf Z_p$ is the ring of $p$-adic
integer numbers, $aB:= \{ x: x=ab, b\in B \} $ for a multiplicative
group $B$ and its element $a\in B$, $k\in \bf N$ \cite{roo,weil}.
The quotient mapping $\pi _k: {\bf K}\to {\bf K}/B({\bf
K},0,p^{-k})$ is defined as well, where $k\in \bf N$.
\par Decompositions of continuous and differentiable functions on
compact subsets of locally compact fields $\bf K$ of zero and
non-zero characteristics with values in $\bf K$ into series of
polynomials were studied in \cite{ami,sch1,cacha} and references
therein. Each function $f\in C^t(M,N)$ has a $C^t(B({\bf K^{\sf
m}},0,1), {\bf K^{\sf n}})$-extension by zero on $B({\bf K^{\sf
m}},0,1)$. Therefore, it has the decomposition
\par $(1)$ $f=\sum_{l,m}f^l_m{\bar Q}_me_l$ \\ in the Amice's basis
${\bar Q}_m$, where $e_l$ is the standard orthonormal in the
non-archimedean sense \cite{roo} basis in $\bf K^{\sf n}$ such that
$e_l=(0,...,0,1,0,...)$ with $1$ in the $l$-th place, $m\in \bf
Z^{\sf n}$, $m_l\ge 0$, $m=(m_1,...,m_{\sf n})$, $f^l_m\in \bf K$
are expansion coefficients such that $\lim_{l+|m|\to \infty }
|f^l_m|_{\bf K}J(t,m)=0$, ${\bar Q}_m$ are polynomials on $B({\bf
K^{\sf m}},0,1)$ with values in $\bf K$, $J(t,m):=\| {\bar Q}_m
\|_{C^t_b(B({\bf K^{\sf m}},0,1),{\bf K})}$. The space $C^t(M,N)$ is
supplied with the uniformity inherited from the space $C^t_b({\bf
K^{\sf m}},{\bf K^{\sf n}})$.
\par {\bf 2. Lemma.} {\it Each $f\in
C^t(M,N)$ is a projective limit $f=pr-\lim_kf_k$ of polynomials
$f_k=\sum_{l,m}f^l_{m,k}{\bar Q}_{m,k}e_l$ on rings ${\bf
S_{p^k}^{\sf m}}={\bf S_{p^k}^{\sf m}}({\bf K})$ with values in
${\bf S_{p^k}^{\sf n}}={\bf S_{p^k}^{\sf n}}({\bf K})$, where
$f^l_{m,k}\in \bf S_{p^k}$ and ${\bar Q}_{m,k}$ are polynomials on
$\bf S_{p^k}^{\sf m}$ with values in $\bf S_{p^k}$.}
\par {\bf Proof.} For each $m\ge k$ consider the quotient mappings (ring
homomorphisms): $\pi _m: B({\bf K},0,1)\to {\bf S_{p^m}}$ and $\pi
^m_k: {\bf S_{p^m}}\to \bf S_{p^k}$ (see \S 1). This induces the
quotient mappings $\pi _m: N\to N_m$ and $\pi ^m_k: N_m\to N_k$,
where $N_m\subset {\bf S_{p^m}}$, $\pi ^m_k\circ \pi _m=\pi _k$,
$\pi ^k_k=id_k: {\bf S_{p^k}}\to {\bf S_{p^k}}$.
\par Let now $M$ and $N$ be two analytic compact manifolds embedded
into $B({\bf K^{\sf m}},0,1)$ and $B({\bf K^{\sf n}},0,1)$
respectively as clopen subsets and $f\in C^t(M,N)$, where $C^t(M,N)$
denotes the space of functions $f: M\to N$ of class $C^t$, $t\ge 0$.
There exists $s\in \bf N$ such that if $x\in M$ and $y\in N$, then
$B({\bf K^{\sf m}},x,p^{-s})\subset M$ and $B({\bf K^{\sf
n}},y,p^{-s})\subset N$. Therefore, consider the cofinal set
$\Lambda _s:= \{ k: k\ge s, k\in {\bf N} \} $ in $\bf N$. For an
integer $t$ it is a space of $t$-times continuously differetiable
functions in the sense of partial difference quotients (see Section
2 and \cite{luseamb,luanmbp,sch1}). Thus $f=pr-\lim_kf_k$, where
$f_k:= \pi _k^*(f)$, $\pi _k^*$ is naturally induced by $\pi _k$
using the polynomial expansion of $f$ (see in details below), where
such decomposition exists for each continuous $f: M\to N$ due to
Formula 1$(1)$ (for the limit of an inverse sequence, see
\cite{bacht}, \S 2.5 \cite{eng} and \S \S 3.3, 12.202 \cite{nari}).
Put $C^t(M_k,N_k) := \pi _k^* \circ C^t(M,N) = \{ f_k: f\in C^t(M,N)
\} $, hence \par $(1)$ $C^t(M,N)\subset pr-\lim_kC^t(M_k,N_k)$ \\
algebraically without taking into account topologies. Thus write it
in the form:
\par $(2)$ $C^t(M,N)= T-pr-\lim \{ C^t(M_k,N_k), \pi ^k_l, \Lambda _s
| C^t(M,N) \} $ algebraically, where
\par $(3)$ $T-pr-\lim \{ P_k, \pi ^k_l, \Lambda | G \} :=
\{ f: pr-\lim \{ f_k,\pi ^k_l, \Lambda \} =f, f\in G, f_k\in P_k ~
\forall k\in \Lambda \} $ \\
denotes the conditional projective limit with a condition $G$, since
$T-pr-\lim \{ P_k, \pi ^k_l, \Lambda | G \} =G\cap pr-\lim \{ P_k,
\pi ^k_l,\Lambda \} $.
\par Indeed, in accordance with \S 1  $f_k=\pi _k^*(f)$ and
\par $(4)$ $\pi _k^*(f(x))=\sum_{l,m}(\pi _k (f^l_m)(\pi _k^* {\bar
Q}_m(x))e_l$,\\ since $\pi _k$ is the ring homomorphism and $\pi
_k(e_l)=e_l$. Then $\pi _k(ax^m)=a_kx^m(k)$ for each $a\in \bf K$
and $x\in B({\bf K^{\sf m}},0,1)$, where $x^m=x_1^{m_1}...x_{\sf
m}^{m_{\sf m}},$ $a_k=\pi _k(a)$ with $a_k\in \bf S_{p^k}$ and
$x^m(k):= \pi _k(x^m)$ with $x(k)\in \bf S_{p^k}^{\sf m},$ hence
\par $(5)$ $\pi _k^*({\bar Q}_m(x))={\bar Q}_{m,k}(x(k)).$ \\
If ${\bar Q}_m(x)=\sum_{q, 0\le q_j\le m_j\forall j}b_qx^q$, where
$q=(q_1,...,q_{\sf m})$, $m=(m_1,...,m_{\sf m})$, $m_j\in \{ 0, 1,
2,... \} =:{\bf N_0}$, $j=1,...,{\sf m}$, $b_q\in \bf K$,
$x^q:=x_1^{q_1}...x_{\sf m}^{q_{\sf m}}$, then
\par ${\bar Q}_{m,k}(x(k))=\sum_{q, 0\le q_j\le m_j\forall j}\pi
_k(b_q)x(k)^q$,
\\ $x(k)=(\pi _k(x_1),...,\pi _k(x_{\sf m}))$, since $\pi _k: {\bf
K}\to {\bf K}/B({\bf K},0,p^{-k})$ is defined as well. The series
for $f_k$ is finite, since $\pi _k(a)=0$ for each $a\in \bf K$ with
$|a|<p^{-k}$ and \par $(6)$ $\lim_{l+|m|\to \infty } |f^l_m|_{\bf
K}J(t,m)=0$. Therefore,
\par $(7)$ $f(x)= pr-\lim \{ f_k(x(k)), \pi ^k_l, \Lambda _s \} $
for each $x\in M$.
\par As shows this proof the index $t$ can be omitted from
$C^t(M_k,N_k)$, where $k\in \Lambda _s$. More precisely we have the
following corollary.
\par {\bf 3. Corolary.} {\it The space $C^t(M_k,N_k)$ is independent
from $t$ and it is algebraically isomorphic with the space
$N_k^{M_k}$ of all mappings from $M_k$ into $N_k$ for each $k\in
\Lambda _s$. Moreover, $({\bf S_{p^k}^{\sf n}})^{({\bf S_{p^k}^{\sf
m}})}$ is a finite-dimensional space over the ring $\bf S_{p^k}$.}
\par {\bf Proof.} In view of Lemma 2
in the module $C^t({\bf S_{p^k}^{\sf m}},{\bf S_{p^k}})$ of the ring
$\bf S_{p^k}$ there is only a finite number of $\bf
S_{p^k}$-linearly independent polynomilas ${\bar Q}_{m,k}(x(k))$,
since the rings ${\bf S_{p^k}^{\sf m}}$ and $\bf S_{p^k}$ are
finite, also $z^a=z^b$ for each natural numbers $a$ and $b$ such
that $a=b\mbox{ } (mod\mbox{ }(p^k))$ and each $z\in \bf S_{p^k}$.
The space $C^t(M_k,N_k)$ is isomorphic with $N_k^{M_k}$, since $M_k$
and $N_k$ are discrete.
\par Therefore, denote $C^t(M_k,N_k)$ by $C(M_k,N_k)$.
\par {\bf 4. Corollary.} {\it There exists the group
$\pi _k^* \circ Diff^{\alpha }(M)$ isomorphic with the symmetric
group $\Sigma _{\sf n_k}$ for every $k\in \Lambda _s$, where $\sf
n_k$ is the cardinality of $M_k$, $\alpha \in \{ n, [n]: 0\le n\le
\infty \} $.}
\par {\bf Proof.} For each $k\in \Lambda _s$ there exists the
mapping $\pi _k: M\to M_k$ (see the proof of lemma 2). Since $\pi
_k$ is the quotient continuous mapping, then $\pi _k(M)=M_k$. If $f:
M\to M$ is a continuous epimorphism, $f(M)=M$, then $\pi
_k(f(M))=M_k$. Let $z\in M_k$, $k\in \Lambda _s$, then there exists
$x\in M$ such that $\pi _k(x)=z$, hence $f_k(z)=\pi_k^*(f(\pi
_k^{-1}(z)))$ due to Formulas 2$(4,5)$, since $\pi _k(\pi
_k^{-1}(z))=z$ and $\pi _k(x)=x(k)=z$. Then $\pi _k^{-1}(M_k)=M$ and
${\pi _k^*}^{-1}(f_k(M_k)) = \bigcup_{z\in M_k} \{ f(x): x\in \pi
_k^{-1}(z) \} $, consequently, ${\pi _k^*}^{-1}(f_k(M_k))=M$, since
for each $x\in M$ there exists $z\in M_k$ such that $\pi _k(x)=z$.
Therefore, if $h\in Diff^{\alpha }(M)$, then $h_k(M_k)=M_k$, since
$h(M)=M$ and $M_k$ is finite. In accordance with Theorem 3.2.14
\cite{eng} if $ \{ \phi ,f_{\sigma '} \} $ is a mapping of an
inverse spectrum $ \{ X_{\sigma }, \pi ^{\sigma }_{\rho }, \Psi \} $
of compacts into an inverse spectrum $ \{ Y_{\sigma '}, \pi ^{\sigma
'}_{\rho '}, \Psi ' \} $ of $T_1$-spaces and all $f_{\sigma '}$ are
epimorphisms on $Y_{\sigma '}$, then the limit mapping $f=pr-\lim \{
\phi , f_{\sigma '} \} $ is also the epimorphism. In view of
Corollary 3 $\pi _k^*\circ Diff^{\alpha }(M)$ is algebraically
isomorphic with the following discrete group $Hom(M_k)$ of all
homeomorphisms $h_k$ of $M_k$, that is, bijective surjective
mappings $h_k: M_k\to M_k$. Using an enumeration of elements of
$M_k$ we get an isomorphism of $Hom(M_k)$ with $\Sigma _{\sf n_k}$.
\par {\bf 5.} Suppose that $C_w(M,N):=pr-\lim_k \{ N_k^{M_k}, \pi ^k_l,
\Lambda _s \} $ is an uniform space of continuous mappings $f: M\to
N$ supplied with an uniformity and projective weak topology as well
inherited from products of uniform spaces $\prod_{k=1}^{\infty
}N_k^{M_k}$ (see also \S 8.2 \cite{eng} and \S 1 above). Denote the
corresponding projective weak topology in $C_w(M,N)$ by $\tau _w$.
The spaces $C^{\alpha }(M,N)$ and $C_w(M,N)$ are subsets of $\bf
K$-linear spaces $C^{\alpha }(M,{\bf K^{\sf n}})$ and $C^0(M,{\bf
K^{\sf n}})$ respectively. Supply with algebraic structures subsets
of the latter $\bf K$-linear spaces as inherited from them.
\par {\bf Corollary.} {\it The uniform space $C^{\alpha }(M,N)$ is not
algebraically isomorphic with $C_w(M,N)$, when $\alpha >0$. The
topological space $C_w(M,N)$ is compact.}
\par {\bf Proof.} In accordance with appendix in \cite{bacht} and
\S 2.5 \cite{eng} and Formulas 2$(1-6)$ above the spaces $C^0(M,N)$
and $C_w(M,N)$ coincide algebraically, since the connecting mappings
$\pi ^m_n$ are uniformly continuous for each $m\ge n$. The space
$C^0(M,{\bf K^{\sf n}})$ is $\bf K$-linear and its uniformity is
completely determined by a neighborhood base of zero. The space
$C_w(M,N)$ is uniformly homeomorphic with $pr-\lim_k({\bf
S_{p^k}})^{M_k}$, which is compact by the Tychonoff Theorem 3.2.4
\cite{eng}. Since $C^0(M,N)\ne C^{\alpha }(M,N)$ for $\alpha >0$,
then $C_w(M,N)$ and $C^{\alpha }(M,N)$ are different algebraically
see Formulas 2$(1-3)$).
\par {\bf Remark.} In general in two consequtive projective limits
of topological spaces two limits may be non commuting.
\par {\bf 6.} Suppose that $Diff_w(M):=pr-\lim_k Hom(M_k)$
is supplied with the uniformity inherited from $C_w(M,M)$. The group
$Diff_w(M)$ is called the non-archimedean compactification of
$Diff^{\alpha }(M)$.
\par {\bf Theorem.} {\it The group $Diff_w(M)$ is the compact
topological group and it is the compactification of $Diff^{\alpha
}(M)$ in the projective weak topology $\tau _w$. If $\alpha >0$,
then $Diff^{\alpha }(M)$ does not coincide with $Diff_w(M)$.}
\par {\bf Proof.} From $Diff^{\alpha }(M)\subset C^{\alpha }(M,M)$
it follows that $Diff^{\alpha }(M)$ has the corresponding algebraic
embedding into $C_w(M,M)$ as the set. Since $C_w(M,M)$ is compact
and $Hom(M)$ is a closed subset in $C_w(M,M)$, then due to Corollary
5 $Hom(M)\cap C_w(M,M)=Diff_w(M)$ is compact. The topological space
$C^{\alpha }(M,M)$ is dense in $C^0(M,M)$, consequently,
$Diff^{\alpha }(M)$ is dense in $Diff_w(M)$. If $\alpha >0$, then
$Diff^{\alpha }(M)\ne Hom(M)$, hence $Diff^{\alpha }(M)$ and
$Diff_w(M)$ do not coincide algebraically. It remains to verify,
that $Diff_w(M)$ is the topological group in its projective weak
topology $\tau _w$. If $f, g \in C^{\alpha }(M,N)$, then $\pi
_k^*({\bar Q}_m(g(x)))={\bar Q}_{m,k}(g_k(x(k))$ due to Formula
2$(4)$, consequently, $\pi _k^*(f\circ g)=\sum_{l,m}\pi
_k(f^l_m){\bar Q}_{m,k}(g_k(x(k))e_l$, hence \par $(1)$ $(f\circ
g)_k=f_k\circ g_k.$ \\ Since $\pi _k(x)=x(k)$, then $\pi
_k^*(id(x))=id_k(x(k)),$ where $id(x)=x$ for each $x\in M$.
Therefore, for $f=g^{-1}$ we have $(f\circ g)_k=f_k\circ g_k=id_k$,
hence \par $(2)$ $\pi _k^*(g^{-1})=g_k^{-1}$.\\ The associativity of
the composition $(f_k\circ g_k)\circ h_k=f_k\circ (g_k\circ h_k)$ of
all functions $f_k, g_k, h_k\in Hom(M_k)$ together with others
properties given above means, that $Diff_w(M)$ is the algebraic
group. \par Indeed, inverse limits of mappings $f=pr-\lim_kf_k$,
$g=pr-\lim_kg_k$ and $h=pr-\lim_kh_k$ satisfy the associativity
axiom as well, each $f$ has the inverse element
$f^{-1}=pr-\lim_kf_k^{-1}$ such that $f^{-1}(f(x))=id$ and $e=id$ is
the unit element. By the definition of the weak topology in
$Diff_w(M)$ for each neighborhood of $e=id$ in  $Diff_w(M)$ there
exists $k\in \bf N$ and a subset $W_k\subset Hom(M_k)$ such that
$e_k\in \pi _k^{-1}(W_k)\subset W$. On the other hand, $Hom(M_k)$ is
discrete, hence there are $e_k\in V_k\subset Hom(M_k)$ and $e_k\in
U_k \subset Hom(M_k)$ such that $V_kU_k\subset W_k$, hence there are
neighborhoods $e\in V\subset Diff_w(M)$ and $e\in U\subset
Diff_w(M)$ such that $VU\subset W$, where $V=\pi _k^{-1}(V_k)$,
$U=\pi _k^{-1}(U_k)$ and $VU= \{ h: h=f\circ g, f\in V, g\in U \} $.
Consider a neighborhood $W'$ of $f^{-1}$, then $V:=W'f^{-1}$ is the
neighborhood of $e$ and there exists $k\in \bf N$ such that $\pi
_k^{-1}(e_k)=:U\subset V^{-1}$, since $e_k^{-1}=e_k$ and $\pi _k$ is
the homomorphism. Thus, $fU:=W$ is the neighborhood of $f$ such that
$W^{-1}\subset W'$, hence the inversion operation $f\mapsto f^{-1}$
is continuous.
\par {\bf 7. Theorem.} {\it The initial $C^{\alpha }$ topology
$\tau _{\alpha }$ and the weak projective limit topology $\tau _w$
in $Diff^{\alpha }(M)$ are incomparable, where $M$ is compact.}
\par {\bf Proof.} Remind that a topology ${\cal O}_1$ in a topological
space $F$ is called weaker, than a topology ${\cal O}_2$, if ${\cal
O}_1\subset {\cal O}_2$ or one says that ${\cal O}_2$ is stronger
than ${\cal O}_1$. Up to a diffeomorphism as above consider $M$
clopen in $B({\bf K}^{\sf n},0,1)$. Since $Diff^{\beta }(M)$ is
contained in $Diff^{\alpha }(M)$ for each $\beta >\alpha \ge 0$ and
$\tau _{\beta }$ is stronger than $\tau _{\alpha }$ in $Diff^{\beta
}(M)$, then it is sufficient to prove this theorem for $\alpha \ge
1$. Each projection $\pi _k^*: C^{\alpha }(M,{\bf K^{\sf n}}) \to
({\bf K_k^{\sf n}})^{M_k}$ induces the quotient metric $\rho _k$ in
the $\bf K_k$-module $({\bf K_k^{\sf n}})^{M_k}$ such that $\rho
_k(f_k,g_k):=\inf_{z, \pi _k(z)=0} \| f-g+z \| _{C^{\alpha } (M,{\bf
K^{\sf n}})},$ where ${\bf K_k}:={\bf K}/B({\bf K},0,p^{-k})$ is the
quotient ring and $\pi _k$ is induced by such quotient mapping from
$\bf K$ onto $\bf K_k$.
\par In view of the Kaplansky theorem which is the non-archimedean analog
of the Stone-Weierstrass theorem and true over fields of zero and
non-zero characteristics in \\ $C^{\alpha }(B({\bf K}^{\sf
n},0,1),{\bf K})$ the set of polynomials is dense \cite{sch1,cacha}.
Thus after application of the quotient mapping we get that in the
module $C({\bf S_{p^k}}^{\sf n}, {\bf S_{p^k}})$ over the finite
ring ${\bf S_{p^k}}$ the set of polynomials is dense. Moreover, each
continuous function has a decomposition into converging series of
polynomials. The locally compact field $\bf K$ is commutative
\cite{weil}, hence the ring ${\bf S_{p^k}}$ is commutative. This
means that the multiplicative group ${\bf S_{p^k}}^*:={\bf
S_{p^k}}\setminus \{ 0 \} $ is commutative and consists of $p^k-1$
elements. Thus, if $0\ne x\in {\bf S_{p^k}}$, then $x^{p^k}=1$.
Therefore, over ${\bf S_{p^k}}$ the set of polynomials is finite
dimensional, hence each $f_k\in C({\bf S_{p^k}}^{\sf n}, {\bf
S_{p^k}})$ is polynomial over ${\bf S_{p^k}}$.
\par Again consider $C^{\alpha }(B({\bf K}^{\sf n},0,1),{\bf K})$,
which is the algebra over $\bf K$. If $f_i(x)=x_i$, where
$x=(x_1,...,x_{\sf n})\in B({\bf K}^{\sf n},0,1)$, then
$\prod_{i=1}^{\sf n}f_i(x)^{s_i}=\prod_{i=1}^{\sf n}x_i^{s_i}\in
C^{\alpha }(B({\bf K}^{\sf n},0,1),{\bf K})$, where
$x_i^s:=\prod_{j=1}^s(x_i)_j$ with $(x_i)_j=x_i$ for each $j$, where
$s, s_i\in \bf N$. In particular, it contains a subalgebra ${\cal
A}_k$ containing all constants from $\bf K$ and all polynomials of
the form \par $(1)$ $\sum_{k_1=p^k l(1),...,k_{\sf n}=p^k l({\sf n})
} a_{k_1,...,k_{\sf n}} x_1^{k_1}...x_{\sf n}^{k_{\sf n}}$,\\ where
$a_{k_1,...,k_{\sf n}}\in \bf K$, $l(i)\in \{ 0, 1, 2,... \} $ with
\par $(2)$ $\sum_{k_1=p^k l(1),...,k_{\sf n}=p^k l({\sf n}) }
a_{k_1,...,k_{\sf n}}=0$, $1\le k\in \bf Z$. \par This algebra
${\cal A}_k$ over $\bf K$ separates points in $C^{\alpha }(B({\bf
K}^{\sf n},0,1),{\bf K})$ and by the Kaplansky theorem ${\cal A}_k$
is everywhere dense in $C^{\alpha }(B({\bf K}^{\sf n},0,1),{\bf
K}))$. But $\pi _k(x_i^{k_i})=(\pi _k(x_i^{p^k}))^{l(i)}=1$ for each
$i$ and $k_i=p^k l(i)$ with $l(i)\in \{ 0, 1, 2,... \} $. Consider a
polynomial with values in ${\bf K}^{\sf n}$ of the form
\par $(3)$ $f=id + \sum_{i=1}^{\sf n} \sum_{k_1=p^k l(1),...,k_{\sf
n}=p^k l({\sf n}) }
a_{i, k_1,...,k_{\sf n}} x_1^{k_1}...x_{\sf n}^{k_{\sf n}}e_i$,\\
where $(f-id)\in {\cal A}_k^{\sf n}$, $|a_{i, k_1,...,k_{\sf n}}|\le
|\pi |$ for each $i, l(1),..., l({\sf n})$, where $\pi \in B({\bf
K},0,1)$, $|\pi |<1$, $|\pi |=\max \{ y\in \Gamma _{\bf K}: 0<y<1 \}
$. Then $\| f -id \|_{C^{\alpha }(B({\bf K}^{\sf n},0,1),{\bf
K})}\le |\pi |$, consequently, $f$ is the isometry and $f\in
Diff^{\alpha }(M)$ and inevitably \par $f_k := \pi _k^*(f)=\pi
_k^*(id) + \sum_{i=1}^{\sf n} \sum_{k_1=p^k l(1),...,k_{\sf n}=p^k
l({\sf n})} \pi _k(a_{i, k_1,...,k_{\sf n}})e_i=\pi _k^*(id)$\\ due
to Condition (2). Therefore, $(\pi _k^*)^{-1}(e_k)$ is everywhere
dense in a neighborhood of $e=id$ in $Diff^{\alpha }(M)$, where
$e_k\in \pi _k^*(Diff^{\alpha }(M))=Hom(M_k)$ is the unit element,
$k>1$, $card (M_k)>1$, $Hom (M_k)$ is the symmetric group of $M_k$
elements.
\par On the other hand, there is ${\sf c}=card ({\bf R})$ elements
$f\in Diff^{\alpha }(M)$ with \par $\| f -id \|_{C^{\alpha }(B({\bf
K}^{\sf n},0,1),{\bf K})}\le |\pi |$ such that \par $f=id
+\sum_{i=1}^{\sf n} \sum_{k_1=p^k l(1),...,k_{\sf n}=p^k l({\sf n})
} a_{i, k_1,...,k_{\sf n}} x_1^{k_1}...x_{\sf n}^{k_{\sf n}}e_i$,
but with \par $\sum_{k_1=p^k l(1),...,k_{\sf n}=p^k l({\sf n}) }
a_{i, k_1,...,k_{\sf n}}\ne 0$ for which $\pi _k^*(f)\ne e_k$.\\
Thus the set ${\pi _k^*}^{-1}(e_k)$ is open in $(Diff^{\alpha
}(M),\tau _w)$, but ${\pi _k^*}^{-1}(e_k)\notin \tau _{\alpha }$. At
the same time, $\{ f\in Diff^{\alpha }(M): \| f-id \| _{C^{\alpha
}(B({\bf K}^{\sf n},0,1),{\bf K}^{\sf n})} \le |\pi | \} $ is open
in $(Diff^{\alpha }(M),\tau _{\alpha })$, but it is not open in
$\tau _w$ topology. This proves the assertion of this theorem, since
neither $\tau _{\alpha }$ nor $\tau _w$ is weaker among two of them
and they are different.
\par {\bf 8. Theorem.} {\it Let $M$ be a $C^{\alpha }$ manifold finite
dimensional over a non-archimedean infinite field $\bf K$ with a non
trivial multiplicative norm complete relative to its uniformity,
where $\bf K$ may be non locally compact. Suppose that $Diff^{\alpha
}_b(M)$ is the group of all uniformly $C^{\alpha }$ continuous
diffeomorphisms of $M$, where $M$ is embedded into ${\bf K}^n$ as
the bounded clopen subset, $\alpha \in \{ t, [t] \} $, $1\le t\le
\infty $. Then $G_s:= \{ g\in Diff^{\alpha }_b(M): \| g-id \|
_{C^{\alpha }(M,{\bf K}^n)}<|\pi |^s \} $ has a non-archimedean
completion, which is a group, where $\pi \in \bf K$, $|\pi |<1$,
$1\le s\in \bf N$.}
\par {\bf Proof.} If $g\in G_s$, then $g$ is the isometry:
$|g(x)-g(y)|=|x-y|$ for each $x, y \in M\subset {\bf K}^n$, since
$1\le t$ and $g(x)-g(y)={\bar {\Phi }}^1g(x;x-y;1)$, where $1\le s$.
Without loss of generality up to an affine diffeomorphism of ${\bf
K}^n$ we can consider, that $M\subset B({\bf K}^n,0,1)$. Consider
the ring homomorphism $\pi _k: B({\bf K}^n,0,1)\to B({\bf
K}^n,0,1)/B({\bf K}^n,0,|\pi |^k)$, where the quotient ring is
discrete and can be supplied with the quotient norm. Then $\pi _k(M)
=: M_k$ is the discrete topological space. \par Since $M$ is clopen
in ${\bf K}^n$ then each $g\in Diff^{\alpha }_b(M)$ has the
extension as $id$ on ${\bf K}^n\setminus M$. Therefore, without loss
of generality consider $M$ such that there exists $q\in \bf N$ for
which from $x\in M$ it follows that $B({\bf K}^n,x,|\pi |^q) \subset
M$. Consider the cofinal set $\Lambda _q := \{ k\in {\bf N}: k\ge q
\} $. If $g\in G_s$, then \par $(1)$ $g(B({\bf K}^n,x,|\pi |^k)) =
B({\bf K}^n,g(x),|\pi |^k)$ for each $B({\bf K}^n,x,|\pi |^k)\subset
M$, since $g$ is the isometry. Consequently, \par $(2)$ $f\circ
g(B({\bf K}^n,x,|\pi |^k)) =f(B({\bf K}^n,g(x),|\pi |^k))=B({\bf
K}^n,f(g(x)),|\pi |^k)$ and \par $g^{-1}(B({\bf K}^n,x,|\pi
|^k)=B({\bf K}^n,g^{-1}(x),|\pi |^k)$ for each $f, g \in G_s$.
\par Therefore, $g$ and $\pi _k$ generate the natural mapping $\mbox{
}_kg: M_k\to M_k$ such that it is bijective and epimorphic, since
$\pi _k^{-1}(z)=B({\bf K}^n,x,|\pi |^k)$ for each $z\in \pi
_k(B({\bf K}^n,0,1))$, and \par $(3)$ $B({\bf K}^n,x,|\pi
|^k)=x+B({\bf K}^n,0,|\pi |^k)$, where $x\in B({\bf K}^n,0,1)$ is
such that $\pi _k(x)=z$, $\pi _k^*(g(x))=\mbox{ }_kg(z):= \pi
_k\circ g\circ \pi _k^{-1}(z)$ for each $z\in M_k$. Thus due to
Equations $(1-3)$ there exists the discrete group $\pi
_k^*(G_s)=:\mbox{ }_kG_s$ and $\pi ^l_k (\mbox{ }_lG_s)=\mbox{
}_kG_s$ for each $l\ge k$, where $\pi ^l_k$ are mappings of the
inverse system such that $\pi ^l_k\circ \pi _l=\pi _k$ such that
$\pi ^l_k$ and $\pi _l, \pi _k$ are algebraic homomorphisms, $\pi
^l_k$ are written without star for simplicity of notation. As in
Theorem 6 it gives the inverse sequence of discrete groups $ S= \{
\mbox{ }_lG_s, \pi ^l_k, \Lambda _q \} $. In view of Lemma 2.5.9
\cite{eng} if $\{ \phi , f_{\sigma '} \} $ is a mapping of an
inverse system $S= \{ X_{\sigma }, \pi ^{\sigma }_{\rho }, \Psi \} $
into an inverse system $S'= \{ Y_{\sigma '}, \pi ^{\sigma '}_{\rho
'}, \Psi ' \} $ and all mappings $f_{\sigma '}$ are injective, then
the limit mapping $f=\lim \{ \phi ,f_{\sigma '} \} $ is also
injective. If moreover, all $f_{\sigma '}$ are surjective, then $f$
is also surjective.
\par  Each discrete topological group is complete
relative to its left uniformity generated by its topology. Thus, the
limit $\lim S$ of the inverse system of discrete groups is the
Tychonoff topological group relative to the projective weak topology
inherited from the product Tychonoff topology $\tau _w$, $G\subset
\lim S \subset \prod_{l\in \Lambda _q}\mbox{ }_lG_s$. Moreover,
$\lim S=:G^w_s$ is the complete uniform space with the left
uniformity ${\cal T}_w$ generated by the left shifts and the
neighborhood base of $e$ in $G^w_s$, since each $\mbox{ }_kG_s$ is
complete (see Theorems 2.5.13, 8.3.6 and 8.3.9 \cite{eng}). We have
that algebraically $G_s\subset G^w_s$, $G_s$ is the topological
group relative to the topology inherited from $G^w_s$. In view of
Theorem 7 the $\tau _w|_{G_s}$ topology is incomparable with the
$C^{\alpha }_b$ uniformly bounded continuous topology.
\par {\bf 9. Corollary.} {\it Let $Diff^{\alpha }_b(M)$
be the group as in Theorem 8. Then $Diff^{\alpha }_b(M)$ has the
non-archimedean completion which is the topological group.}
\par {\bf Proof.} Let $G_s$ be the subgroup as in Theorem 8, then
$G_s$ is clopen in $Diff^{\alpha }_b(M)$. There exists a family
$g_j\in Diff^{\alpha }_b(M)$ such that $\bigcup_{j\in \Psi
}g_jG_s=Diff^{\alpha }_b(M)$, where $\Psi $ is a set. We have that
$gG_s$ is clopen in $Diff^{\alpha }_b(M)$ for each $g\in
Diff^{\alpha }_b(M)$, since $L_g: Diff^{\alpha }_b(M)\to
Diff^{\alpha }_b(M)$ is the homeomorphism, where $L_g(h)=gh$ for
each $g, h\in Diff^{\alpha }_b(M)$. Then $g_jG_s\cap g_iG_s=
g_j(G_s\cap g_j^{-1}g_iG_s)$. Let $g_{j_0}=id$ and $g_j\ne id$ for
each $j\ne j_0$. \par The group $Diff^{\alpha }_b(M)$ is metrizable,
since $C^{\alpha }_b(M,{\bf K}^n)\subset C^{\alpha }_b(B({\bf
K}^n,0,1),{\bf K}^n)$ is metrizable, hence it is paracompact. Each
two balls in $C^{\alpha }_b(M,{\bf K}^n)$ either are disjoint or
coincide, since it is the normed space. We have that $Diff^{\alpha
}_b(M)$ is contained in $C^{\alpha }_b(M,{\bf K}^n)$ and at the same
time $Diff^{\alpha }_b(M)$ is the neighborhood of $id$ in $C^{\alpha
}_b(M,{\bf K}^n)$, $B(C^{\alpha }_b(M,{\bf K}^n),id,|\pi |^s)=G_s$
for $s\ge 1$ and $1\le t$, $\alpha \in \{ t, [t] \} $. Therefore
choose $\Psi $ such that $g_jG_s\cap g_iG_s=\emptyset $ for each
$i\ne j\in \Psi $, hence $g_i^{-1}g_j\notin G_s$ for each $i\ne j$.
\par The minimal group $gr (\bigcup_{j\in \Psi
}g_jG_s)=Diff^{\alpha }_b(M)$ and it is contained in the minimal
algebraic group $gr(\bigcup_{j\in \Psi }g_jG^w_s)=:G^w$. Supply the
latter group with the uniformity ${\cal T}_w$ induced from the base
of neighborhoods of $e$ in $G^w_s$ with the help of left shifts.
Then the restriction of ${\cal T}_w$ on $G^w_s$ coincides with that
of ${\cal T}_w$ in Theorem 8. Since $G^w_s$ is clopen in $G^w$ and
$G^w_s$ is complete, then for each Cauchy net $ \{ h_q: q\in \nu \}
$ in $G^w$, where $\nu $ is an ordered set, there exists $q_0$ such
that for each $q, l\ge q_0$ there is the inclusion $h_q^{-1}h_l\in
G^w_s$, hence $ \{ h_{q_0}^{-1}h_l: l >q_0 \} $ converges in $G^w_s$
to some $g_s\in G^w_s$, since $G^w_s$ is complete, consequently, $
\{ h_q: q\in \nu \} $ converges in $G^w$ to $h_{q_0}g_s\in G^w$ and
inevitably $G^w$ is complete (see also Theorem 8.3.20 \cite{eng}).
This $G^w$ is the desired non-archimedean completion.

\section{Example of the group of diffeomorphisms.}
\par This section contains the example of the group of diffeomorphisms.
It illustrates the general theory. For the group of diffeomorphisms
of $\bf Z_p$ of class $C^t$, where $0\le t \le \infty $, formulas
for expansion coefficients in the Mahler base of compositions
$g\circ f$ and inverse elements $f^{-1}$ of diffeomorphisms $f$ and
$g$ are found.
\par Let $C^t$ be a class of smoothness
of functions $f: M\to \bf K$ as in \S 1 (see also
\cite{ami,sch1,luseamb,lutmf99}), where $M$ is a Banach manifold
over a complete (as an uniform space) non-Archimedean infinite field
$\bf K$ with non-trivial valuation and of zero characteristic
$char({\bf K})=0$. Suppose that $L_bf(x):=f(x+b)$ is a shift
operator, $\Delta _bf(x):=(L_b-I)f(x)$ is a difference operator,
$L^0=I$, $\Delta ^0=I$, $L:=L_1$, $\Delta :=\Delta _1$, where $x,
b\in \bf K$. For a product of two functions $f, g: {\bf K}\to \bf K$
there are formulas:
\par $(1)$ $\Delta [f(x)g(x)]=(\Delta f)(x)(Lg)(x)+f(x)(\Delta g)(x)$ and
\par $(2)$ $L\Delta =\Delta L$. Therefore,
\par $\Delta ^k[f(x)g(x)]=\sum_{j=0}^k{k\choose j}[\Delta ^jf(x)]
L^j[\Delta ^{k-j}g(x)]$ for each $k\in \bf N$, consequently,
\par $(3)$ $\Delta ^kf_1(x)...f_n(x)=\sum_{k_1+...+k_n=k}
(k!/(k_1!...k_n!))[\Delta ^{k_1}f_1(x)]L^{k_1}[\Delta ^{k_2}f_2(x)]
L^{k_2}...$ \\
$[\Delta ^{k_{n-1}}f_{n-1}(x)]L^{k_{n-1}} [\Delta ^{k_n}f_n(x)]$,
where $L^{k_j}$ acts on all functions situated on the right from it.
If $f\in C^t({\bf Z_p},{\bf Q_p})$, then there is its expansion in
the Mahler base ${x \choose n}$ as a series:
\par $(4)$ $f(x)=\sum_{j=0}^{\infty }f_j {x \choose j}$,
where $f_j\in \bf Q_p$,
\par $(5)$ $\lim_{j\to \infty }|f_j|j^t=0$ for $0\le t<\infty $ and
\par $(6)$ $f_j=[\Delta ^jf(x)]|_{x=0}$, since
\par $(7)$ $\Delta {x \choose j}={x \choose {j-1}}$
for each $j\in \bf N$, ${x \choose 0}=1$ and ${x \choose m} :=0$ for
each $0>m\in \bf Z$ (see \S 52 \cite{sch1}, \cite{riordan} and
\cite{lutmf99}). Therefore,
\par $(8)$ $(g\circ f)_k=\sum_ng_n\Delta ^k { {\sum_{m=0}^{\infty }
f_m {x \choose m} }\choose n}|_{x=0}$. Since
\par $(9)$ ${x\choose n}=x(x-1)...(x-n+1)/n!$, then
\par $(10)$ $\Delta ^k{ {{f(x)}\choose m}\choose n}=\sum_{k_1+...+k_n=k}
(k!/(k_1!...k_n!))[\Delta ^{k_1}{ {f(x)} \choose m}]L^{k_1} [\Delta
^{k_2}({ {f(x)} \choose m} -1)] L^{k_2}...[\Delta ^{k_{n-1}}({
{f(x)} \choose m} -n+2)]L^{k_{n-1}} [\Delta ^{k_n}({ {f(x)} \choose
m} -n+1)]/n!$,  where
\par $(11)$ $\Delta ^k {{f(x)}\choose m} =m^{-1}\sum_{l=0}^k
{k\choose l} [\Delta ^lL^{k-l}{ {f(x)}\choose {m-1}} -\delta
_{l,0}(m-1)]$. Since
\par $(12)$ $\Delta ^lL^{k-l}f(x)|_{x=0}=\sum_{m=0}^{\infty }
f_m\Delta ^l{{x+k-l} \choose m}|_{x=0}=\sum_mf_m{{k-l}\choose
{m-l}}$, then
$$(13)\mbox{ }\Delta ^k{ {{f(x)}\choose m}\choose n}|_{x=0}
=(n!)^{-1}\sum_{l_1+...+l_n=k} {k\choose {l_1}} {{k-l_1}\choose
{l_2}} ... { {k-l_1-...-l_{n-2}} \choose {l_{n-1}} }$$
$$ [\sum_{m=0}^{\infty } f_m{{k-l_1}\choose
{m-l_1}} - \delta _{l_1,0} (n-1)] [\sum_{m=0}^{\infty }
f_m{{k-l_1-l_2}\choose {m-l_2}} - \delta _{l_2,0} (n-2)] ...$$
$$[\sum_{m=0}^{\infty } f_m{{k-l_1-...-l_{n-1}}\choose
{m-l_{n-1} } } - \delta _{l_{n-1},0}]f_{l_n},$$ so that it is
necessary to evaluate coefficients
$$(14)\mbox{ }\Omega ^{k,n}_{m_1,...,m_n}:=\sum_{l_1+...+l_{n-1}=k-m_n}
{k\choose {l_1}} {{k-l_1}\choose {l_2}} ... { {k-l_1-...-l_{n-2}}
\choose {l_{n-1}} } {{k-l_1}\choose {m_1-l_1}}$$
$${{k-l_1-l_2}\choose {m_2-l_2}}
.. {{k-l_1-...-l_{n-1}}\choose {m_{n-1}-l_{n-1} } } .$$ There are
identities:
\par $(15)$ $ {m\choose l} = (-1)^l {{l-m-1}\choose l} $,
\par $(16)$ $(1+x)^e(1+(1+x)^a)^b= \sum {b\choose c} {{ac+e}\choose d}x^d$,
such that for $e=b$ and $a=-1$ this gives
$$(17)\mbox{ }\sum_{l_1,l_2} { {k-l_1} \choose {m_1-l_1}}
{{k-l_1-l_2}\choose {m_2-l_2}} {k\choose {l_1}} {{k-l_1} \choose
{l_2}}=$$
$$\sum_{l_1} {{k-l_1} \choose {m_1-l_1}} {k\choose {l_1}}
{{k-l_1}\choose {m_2}}2^{m_2}=$$
$${k\choose {m_1}} \sum_{l_1} {{m_1}\choose {l_1}} {{k-l_1}
\choose {m_2}}2^{m_2},\mbox{ since}$$
$$\sum_l{{b-l}\choose {m-l}} {b\choose l}=\sum_lb!/[(m-l)!(b-m+l)!l!].$$
General expressions are too complicated, but it is suffucient to
construct a generating function for these coefficients:
$$(18)\mbox{ }y_1^{m_1}...y_{n-2}^{m_{n-2}}y_{n-1}^{m_{n-1}+m_n}
(x_1y_2...y_n+(y_1+z_2)(x_2y_3...y_{n-1}+(y_2+z_3)(x_3y_4...y_{n-1}+$$
$$(y_3+z_4)(x_4y_5...y_{n-1}+...))...)^k=$$
$$y_1^{m_1}...y_{n-2}^{m_{n-2}}y_{n-1}^{m_{n-1}+m_n}
\sum_{l_1} {k\choose {l_1}} (x_1y_2...y_{n-1})^{l_1}
(y_1+z_2)^{k-l_1}(x_2y_3...y_{n-1}+(y_2+z_3) (x_3y_4...y_{n-1}+$$
$$(y_3+z_4)(x_4y_5...y_{n-1}...+(x_{n-1}+y_{n-1})...
)^{k-l_1}=$$
$$y_2^{m_2}...y_{n-2}^{m_{n-2}}y_{n-1}^{m_{n-1}+m_n}
\sum_{l_1,q_1,l_2}{k\choose {l_1}} {{k-l_1}\choose
{q_1-l_1}}{{k-l_1}\choose {l_2}} (x_1y_2...y_{n-1})^{l_1}
y_1^{k-q_1+m_1}z_2^{q_1-l_1}$$
$$(x_2y_3...y_{n-1})^{l_2}
(y_2+z_3)^{k-l_1-l_2}(x_3y_3...y_{n-1}+(y_3+...)...)^{k-l_1-l_2}=$$
$$y_3^{m_3}...y_{n-2}^{m_{n-2}}y_{n-1}^{m_{n-1}+m_n}
\sum_{l_1,q_1,l_2,q_2} {k\choose {l_1}} {{k-l_1}\choose
{q_1-l_1}}{{k-l_1}\choose {l_2}} {{k-l_1-l_2}\choose
{q_2-l_2}}(x_1y_3...y_{n-1})^{l_1}y_1^{k+m_1-q_1}$$
$$z_2^{q_1-l_1}(x_2y_3...y_{n-1})^{l_2}y_2^{k+m_2-q_2}z_3^{q_2-l_2}
(x_3y_4...y_{n-1}+(y_3+z_4(x_4y_5...y_{n-1}+...)...)^{k-l_1-l_2}=...
.$$ In this series coefficients in front of $x_1^{m_1}x_2^{m_2}...
x_{n-1}^{m_{n-1}}(y_1...y_{n-1})^k$ are equal to $\Omega
^{k,n}_{m_1,...,m_n}$, where $x_j=z_{j+1}$ for each $j=1,...,n-1$,
$x_j$ and $y_j$ are variables. Therefore, the generating function
has the form:
$$(19)\mbox{ }y_1^{m_1}...y_{n-2}^{m_{n-2}}y_{n-1}^{m_{n-1}+m_n}
(x_1y_2...y_n+(x_1+y_1)(x_2y_3...y_{n-1}+(x_2+y_2)(x_3y_4...y_{n-1}+$$
$$(x_3+y_3)(x_4y_5...y_{n-1}...+(x_{n-1}+y_{n-1})...)^k=$$
$$\sum_{m_1,...,m_{n-1}}\Omega ^{k,n}_{m_1,...,m_n}x_1^{m_1}x_2^{m_2}...
x_{n-1}^{m_{n-1}}(y_1...y_{n-1})^k+$$
$$\sum_{m_1,q_1,...,m_{n-1},q_{n-1}}\Upsilon ^{k,n}_{m_1,q_1,...,m_{n-1},
q_{n-1}}x_1^{q_1}x_2^{q_2}... x_{n-1}^{q_{n-1}} \times $$
$$y_1^{k+m_1-q_1}...y_{n-2}^{k+m_{n-2}-q_{n-2}}
y_{n-1}^{k+m_{n-1}-q_{n-1}+m_n},$$ where coefficients $\Upsilon
^{k,n}_{m_1,q_1,...,m_{n-1},q_{n-1}}$ are given by Equation $(18)$.
\par In particular $id=f^{-1}\circ f$ and $\Delta ^kid(x)|_{x=0}=
\delta _{k,1}$, hence
\par $(21)$ $\delta _{k,1}=\sum_{n=0}^{\infty }(n!)^{-1}(f^{-1})_n
Q^{k,n}(f)$, where coefficients $Q^{k,n}$ are given by Equations
$(8,13,14,19)$, that is
$$(22)\mbox{ }Q^{k,n}(f)=\sum_{m_1,...,m_n}\sum_{i_1<...<i_{n-p}; p}
\Omega
^{k,n}_{m_{i_1},...,m_{i_{n-p}}}f_{m_{i_1}}...f_{m_{i_{n-p}}}$$
$$\prod_{(j_s\in (1,...,n)\setminus \{ i_1,...,i_{n-p} \} ; j_1<j_2<...;
l_{j_1}=0,...,l_{j_{p-1}}=0)}(j_s-n){{k-l_1-...-l_{s-1}}\choose
{m_s}}.$$
\par For analytic functions there are equalities:
$$(23)\mbox{ }\Delta ^kx^n=\sum_{l_1=0}^{n-1}{n\choose {l_1}}\Delta ^{k-1}
x^{l_1}=$$
$$\sum_{l_1=0}^{n-1}\sum_{l_2=0}^{l_1-1}...\sum_{l_k=0}^{l_{k-1}-1}
{n\choose {l_1}}{{l_1}\choose {l_2}} ... {{l_{k-1}}\choose {l_k}}
x^{l_k},\mbox { consequently,}$$
$$(24)\mbox{ }(\Delta ^kx^n)_{x=0}=
\sum_{l_1=0}^{n-1}\sum_{l_2=0}^{l_1-1}...\sum_{l_{k-1}=0}^{l_{k-2}-1}
{n\choose {l_1}}{{l_1}\choose {l_2}} ... {{l_{k-2}}\choose
{l_{k-1}}}=: T_{n,k}.$$ On the other hand,
\par $(25)$ ${x\choose m}=(m!)^{-1}\sum_{l=0}^mx^{m-l}(-1)^l\alpha _l
(1,...,m-1)=$ \\
$\sum_lS_{m,l}x^l$, where
\par $(26)$ $\alpha _l(z_1,...,z_m):=\sum_{i_1<i_2<...<i_l}z_{i_1}...
z_{i_l}$, consequently,
\par $(27)$ $\sum_lS_{m,l}T_{l,j}=\delta _{m,j}$ and
\par $(28)$ $\sum_lT_{m,l}S_{l,j}=\delta _{m,j}$. Then
for analytic $g$ and $f$:
\par $(29)$ $g(x)=\sum_ma_m(g)x^m$ and
$$(30)\mbox{ }(g\circ f)(x)=\sum_{m,l_j,k_j,n}a_m(g)[m!/(l_1!...l_m!)]
a_{k_1}^{l_1}(f)...a_{k_m}^{l_m}(f)]x^n,$$ where in the last series
$k_1l_1+...+k_ml_m=n$, $l_1+...+l_m=m,$ $0\le l_j, k_j \in \bf Z$.
\par For estimations of
$\Omega ^{k,n}_{m_1,...,m_n}f_{m_1}...f_{m_n}$ it can be used
\par $(31)$ $|{k\choose l}|_p=p^{-\lambda (k)+\lambda (q)+\lambda (k-q)}$,
where
\par $(32)$ $\lambda (n)=(n-s_n)/(p-1)$, $n=a_0+a_1p+...+a_jp^j$,
$s_n:=a_0+a_1+...+a_j$, where $a_i\in \{ 0,1,...,p-1 \} $, $\lambda
(q)+\lambda (k-q)-\lambda (k)=(s_k-s_q-s_{k-q})/(p-1)$ (see also \S
25 in \cite{sch1} and \cite{riordan}).

\section{One-parameter subgroups of diffeomorphism groups.}
\par {\bf 1. Theorem.} {\it If $M$ is a compact manifold
over a locally compact field $\bf K$ of characteristic $char ({\bf
K}) = p>1$, $1\le t\in \bf N$, $\alpha \in \{ t, [t] \} $, then
there exists a clopen subgroup $W$ in $Diff^{\alpha }(M)$ such that
each element $g\in W$ lies on a local one-parameter subgroup $g^x$
continuous by $x$ relative to the multiplicative group $({\bf
K}^*,\times )$ with $x\in {\bf K}^*$.}
\par {\bf Proof.} Since $M$ is compact, then $M$ is finite dimensional
$dim_{\bf K}M=n$ over $\bf K$ and the manifold $M$ can be supplied
with a finite disjoint analytic atlas, that is, with disjoint clopen
charts a finite union of which covers $M$. Therefore, there exists a
natural $C^{\alpha }$ embedding of $M$ into ${\bf K}^n$ as the
clopen subset. As it was proved above $Diff^{\alpha }(M)$ is
mertizable with the left-invariant metric $\rho ^{\alpha }$. Choose
$k_0\in \bf N$ and the clopen subgroup $W := \{ g\in Diff^{\alpha
}(M): \rho ^{\alpha }(g,id)<|\pi |^{k_0} \} $ such that $g\in W$
implies $\| g-id \| _{C^{\alpha }_b(M,{\bf K}^n)}<|\pi |^{k_0}$.

\par Consider $x\in B({\bf K},1,|\pi |^s)$, where $1\le s\in \bf N$,
then $|x|=1$, since $|x-1|\le |\pi |^s$ and due to the ultrametric
inequality. If $g\in W$, then for the proving of this theorem it is
necessary to find a one-parameter local subgroup $\{ g^x: x\in
B({\bf K},1,|\pi |^s) \} $ satisfying Conditions $(1-3)$:
\par $(1)$ $g^1=id$;
\par $(2)$ $g^{x_1}g^{x_2} = g^{x_1x_2}$ for each $x_1,
x_2\in B({\bf K},1,|\pi |^s)$;
\par $(3)$ $g^{x_0}=g$ for some $x_0\in B({\bf K},1,|\pi |^s)$. \\
Each $g\in Diff^{\alpha }(M)$ has the form $g=(g_1,...,g_n)$, where
$g_j: M\to \bf K$ for each $j\in \{ 1,...,n \} $. Since $M$ is
embedded as the compact subset in ${\bf K}^n$, then its covering by
balls $B({\bf K}^n,x,r)\subset M$ has a finite subcovering $\{
B({\bf K}^n,x_j,r_j): j=1,...,m \} $, where $x\in M$, $0<r_j<\infty
$, $m\in \bf N$, consequently, $\min_{j=1,...,m}r_j = r>0$. \par If
$g^{x_1}$ and $f^{x_2}$ are two commuting local one-parameter
subgroups for each $x_1, x_2\in B({\bf K},1,|\pi |^s)$ in
$Diff^{\alpha }(M)$, then $g^xf^x=:(gf)^x$ is a local one-parameter
subgroup in $Diff^{\alpha }(M)$, since
$(gf)^{x_1}(gf)^{x_2}=(g^{x_1}f^{x_1})(g^{x_2}f^{x_2})=
g^{x_1}g^{x^2}f^{x_1}f^{x_2}=g^{x_1x_2}f^{x_1x_2}=(gf)^{x_1x_2}$ for
each $x_1, x_2\in B({\bf K},1,|\pi |^s)$, also $g^1=f^1=e=id$. For
example, if $supp (g^x)\subset A$ and $supp (f^x)\subset C$ for each
$x\in B({\bf K},1,|\pi |^s)$ such that $A\cap C=\emptyset $, where
$A$ and $C$ are closed subsets in $M$, $supp (g) := cl_M \{ y\in M:
g(y)\ne y \} $, $cl_M(V)$ denotes the closure of a subset $V$ in
$M$, then $g^{x_1}f^{x_2}=f^{x_2}g^{x_1}$ for each $x_1\ne x_2\in
B({\bf K},1,|\pi |^s)$. On the other hand, each $g\in Diff^{\alpha
}(M)$ can be decomposed into the product $g=h_1...h_m$, where $supp
(h_j)\subset B({\bf K}^n,x_j,r_j)$ for each $j=1,...,m$. Therefore,
the proof of this theorem reduces to the case $Diff^{\alpha }(B)$,
where $B$ is a clopen compact ball in ${\bf K}^n$, since $W$
decomposes into the internal direct product of its subgroups $W_j :=
\{ g\in W: supp (g)\subset B({\bf K}^n,x_j,r_j) \} $ and $W_j$ has
the natural embedding into $Diff^{\alpha }(B({\bf K}^n,x_j,r_j))$
for each $j=1,...,m$. For $g=id$ put $id^x=id$ for each $x\in B({\bf
K},1,|\pi |^s)$. Therefore, due to Equations 3.8$(1,2)$ it remains
the case of $g\ne id$ in $W$, consequently, $x_0\ne 1$ for such $g$.
\par In accordance with Lemma 7.6 \cite{herstein} if $\bf K$ is a
field and $G$ is a finite subgroup of the multiplicative group of
nonzero elements of $\bf K$, then $G$ is a cyclic group. Then by
Theorem 7.b the multiplicative group of nonzero elements of a finite
field is cyclic. In view of the Wedderburn Theorem 7.c a finite
division ring is necessarily a commutative field \cite{herstein}.
The field $\bf K$ is isomorphic with ${\bf F}_{p^u}(\theta )$, where
${\bf F}_{p^u}$ is the finite field consisting of $p^u$ elements,
$u\in \bf N$. The finite field ${\bf F}_{p^u}$ is the splitting
field of the polynomial $x^{p^u}-x$. Thus $\chi : {\bf K}^* \to {\bf
K}^*$ is a continuous multiplicative character, where ${\bf K}^* =
{\bf K}\setminus \{ 0 \} $, if and only if it can be written in the
form $\chi (x)=\phi (x)\psi (x)$, where $\phi $ and $\psi $ are
continuous multiplicative characters such that $\phi : {\bf
F}_{p^u}^*\to {\bf F}_{p^u}^*$ is some homomorphism of the
multiplicative group of ${\bf F}_{p^u}$ and $\phi (\theta )=\theta
$, $\psi (\theta )=\theta ^k$ for some nonnegative integer $k$ and
the restriction of $\psi $ on ${\bf F}_{p^u}$ is the identity
mapping, since there exists the natural embedding of ${\bf F}_{p^u}$
into $\bf K$. Therefore, $|\chi (\theta )| =|\theta |^k$. Denote by
$\Omega $ the family of all continuous multiplicative characters
$\chi : {\bf K}^*\to {\bf K}^*$. In particular, $\chi _1(x)=1$ and
$\chi _{id}(x)=x$ for each in $x\in {\bf K}^*$ are the trivial
character and the identity character respectively.
\par If $g^x$ is a local one-parameter subgroup in $Diff^t(M)$, then
$g^{\chi (x)}$ is also a local one-parameter subgroup in
$Diff^{\alpha }(M)$. If $g^{x_0}=g$, then $g^z=g$ for each $\chi
(z)=x_0$ with $\chi \in \Omega $, since $\chi (1)=1$. Therefore, if
$x_0\in B({\bf K},1,|\pi |^s)$ is a marked point, then it is
sufficient to satisfy Condition $(3)$ for some $x\in V_s := \{ \chi
^{-1}(x_0): \chi \in \Omega \} \cap B({\bf K},1,|\pi |^s)$, where
$x_0\ne 1$ and $x\ne 1$ for $g\ne 1$. Thus, if each $h_j\ne id$ in
$W_j$ has $x_0(h_j)\in V_s$, then the desired local one-parameter
subgroup $g^x$ will be found for $id \ne g=h_1...h_m\in W$. Take,
for example, $x_0=1+\theta ^s$.
\par We say that a family $A$ of functions $f: X\to \bf K$ separate
points of $X$ if for each $x\ne y\in X$ there exists $f\in A$ such
that $f(x)\ne f(y)$. If $G$ is a cyclic group of order $k$ and $a\in
G$ is such that $ \{ 1, a, a^2,..., a^{k-1} \} =G$ and $m$ is a
natural number mutually prime with $k$, $(m,k)=1$, then $\zeta :
G\to G$ such that $\zeta (a^l)=a^{lm}$ for each $l=0,...,k-1$ is the
automorphism of $G$, in particular, for $k=p^u-1$. Since the
multiplicative group $G={\bf F}_{p^u}^*$ is cyclic, then the family
of all $\phi $ separate points of $G$, hence the family $\Omega $ of
all continuous multiplicative characters $\chi $ of ${\bf K}^*$
separate points of ${\bf K}^*$.
\par In view of the Kaplansky theorem a subalgebra $A$ of $C^0(X,{\bf
K})$ containing all constant functions and separating points of a
locally compact totally disconnected Hausdorff space $X$ is dense in
$C^0(X,{\bf K})$ \cite{kaplpams,arsch}. From \cite{ami} it follows,
that $C^0(M,{\bf K}^n)$ has the polynomial basis $\{ Q_{\bar
m}(y)e_i: {\bar m}\in {\bf N_0}^n, i=1,...,n \} $ with $N_0 = \{
0,1,2,... \} $ and $e_i=(0,...,0,1,0,...)\in {\bf K}^n$ with $1$ on
the i-th place.
\par In view of the proof above up to the affine $C^{\infty }$
diffeomorphism it is sufficient to consider $W$ for $Diff^t(B)$,
where $B= B({\bf K}^n,0,1)$. We seek a solution in the form
\par $(4)$ $g_j^x(y)=y_j+\sum_{k=1}^{\infty }b_{j,k}(y,x)$\\
for each $j=1,...,n$ with the converging series of functions
$b_{j,k}(y,x): M\times B({\bf K},1,|\pi |^s)\to \bf K$ of class
$C^t$ by $y\in M$ and continuous by $x\in B({\bf K},1,|\pi |^s)$ for
each $j=1,...,n$ such that $\| b_{j,k}\|_{C^t_b(B,{\bf K})}<|\pi
|^{k_0}$ and $\lim_{k\to \infty } \max_{j=1,...,n} \| b_{j,k}
\|_{C^t(B,{\bf K})}=0$, where $g=(g_1,...,g_n)$. Each function
$b_{j,k}(y,x)$ has the decomposition \par $(5)$ $b_{j,k}(y,x)=\chi
_{j,k}(x)
\sum_{{\bar m}\in {\bf N_0}^n} c_{j,k,{\bar m}} Q_{\bar m}(y)$, \\
where $\chi _{j,k}(x)\in \Omega $, $y\in B$, $x\in B({\bf K},1,|\pi
|^s)$, $c_{j,k,{\bar m}}\in \bf K$. Without loss of generality we
can suppose that $|\chi _{j,k}(x)|\le 1$ on $B({\bf K},1,|\pi |^s)$
for each $j, k ,m$, since $B$ is the unit ball and $s\ge 1$. The
convergence of Series $(4)$ is equivalent to:
\par $(6)$ $\lim_{k+|{\bar m}|\to \infty }
\max_{j=1,...,n} \| c_{j,k,{\bar m}}Q_{\bar m}(y)\|_{C^{\alpha
}_b(M,{\bf K})}
\| \chi _{j,k}(x) \| _{C^0_b(B({\bf K},1,|\pi |^s),{\bf K})}=0$, \\
where $|{\bar m}|:=m_1+...+m_n$, $m=(m_1,...,m_n)$. Then Condition
$(1)$ is equivalent to $\sum_{k=1}^{\infty }c_{j,k,{\bar m}}=0$ for
each $\bar m$ and $j$, while Condition $(3)$ is equivalent to
$\sum_{\bar m}c_{j,{\bar m}}Q_{\bar m}(y) = g_j(y)-y_j$ with
$c_{j,{\bar m}} := \sum_{k=1}^{\infty }c_{j,k,{\bar m}}\chi
_{j,k}(x_0)$ , since $\chi _{j,k}(1)=1$, $\{ Q_{\bar m}(y): {\bar
m}\in {\bf N_0}^n \} $ are linearly independent over $\bf K$.

\par Consider two monotone increasing sequences $\{ s_v\in {\bf N}:
v\in {\bf N} \} $ and $\{ q_v\in {\bf N}: v\in {\bf N} \} $ such
that $q_v\ge s_v$ for each $v$. For each $q>l\in \bf N$ there exists
the quotient algebraic homomorphism of rings $\pi ^q_l: B({\bf
K},0,1)/B({\bf K},0,|\pi |^q)\to B({\bf K},0,1)/B({\bf K},0,|\pi
|^l)$ and $\pi _q: B({\bf K},0,1)\to B({\bf K},0,1)/B({\bf K},0,|\pi
|^q)$ such that $\pi ^q_l\circ \pi _q=\pi _l$. These algebraic
homomorphisms induce natural homomorphisms of functions in the same
notation such that $\pi _q^*(g) =: \mbox{ }_qg$, where $\mbox{
}_qg\in Hom (M_q)$, $M_q=\pi _q(M)$ is a finite set, $g\in W\subset
Diff^t(B)$, $\mbox{ }_qg: M_q\to M_q$ is a bijective surjective
mapping of $M_q$ onto $M_q$ (see Formulas 3.2$(4-7)$). Then also
$\pi _{s_v}(B({\bf K}^*,1,|\pi |^s))$ is the finite set for $s_v>s$.
For $s_v>1$ we have that $\pi _{s_v}({\bf K})$ is the discrete
commutative ring (see \S 3.1). If $x\in \bf K$ and $|x|=1$, then
$|x^{-1}|=1$. The multiplicative norm in $\bf K$ induces the norm in
$\pi _{s_v}({\bf K})$. If $|x-1|<1$, then
$x^{-1}=(1+(x-1))^{-1}=1+\sum_{l=1}^{\infty }(1-x)^l$, thus $\pi
_{s_v}(x)$ has the inverse in $\pi _{s_v}(B({\bf K},1,|\pi |^s))$,
when $|x-1|\le |\pi |^s$, $s\ge 1$, since $\pi _{s_v}((x-1)^l)=0$
for $ls\ge s_v$. Then $1\in \pi _{s_v}(B({\bf K},1,|\pi |^s)$ and
$\pi _{s_v}(x_1x_2)=\pi _{s_v}(x_1)\pi _{s_v}(x_2)\in \pi
_{s_v}(B({\bf K},1,|\pi |^s))$ for each $x_1, x_2\in B({\bf
K},1,|\pi |^s)$, since $|x_1x_2-1|\le \max (|x_1-1| |x_2|,
|x_2-1|)\le |\pi |^s$ and $|x_2|=1$. Thus, $B({\bf F}_{p^u},1,|\pi
|^s)$ for a natural number $s\ge 1$ is the multiplicative
commutative group and $B({\bf F}_{p^u},1,|\pi |^{s_1})/B({\bf
F}_{p^u},1,|\pi |^{s_2})$ for each natural numbers $1\le s_1<s_2$ is
the finite multiplicative commutative quotient group.
\par If $g\in W$ and
$k_0\ge 1$, then $|g(y)-y |\le \| g -id \|_{C^1}|y|\le |\pi
|^{k_0}|y|$ and $|g(y_1)-g(y_2)|\le \| g\|_{C^1}|y_1-y_2|$ and
$|g^{-1}(y_1)-g^{-1}(y_2)|\le \| g^{-1} \| _{C^1}|y_1-y_2|$, since
$t\ge 1$, consequently, $g$ is the isometry, that is, $|g(y)|=|y|$
for each $y\in M$. \par For given $s$ and $s_1$ choose $k_0$ and
$q_1\ge s_1$ such that for each $g\in W$ we have $\mbox{
}_{q_1}g=id$ on $M_{q_1}$. Then for a given $g\in W$ there exists
the algebraic homomorphism $\eta _1$ from $\pi _{s_1}(B({\bf
K},1,|\pi |^s))$ into $\pi _{q_1}(W)$ satisfying Conditions $(1-3)$
with $\pi _{s_1}(x_0)$ instead of $x_0$. Each $g\in W$ and each
$x\in \bf K$ is the projective limit of the inverse sequences $g =
pr-\lim \{ \mbox{ }_{q_v}g, \pi ^{q_v}_{q_l}, \Lambda _{q_1} \} $
and $x = pr-\lim \{ \mbox{ }_{s_v}x, \pi ^{s_v}_{s_l}, {\bf N} \} $,
since sequences $\{ q_v \} $ and $ \{ s_v \} $ are cofinal with $\bf
N$, where $v\ge l\in \bf N$ (see also \S 2.5 \cite{eng}). For
example, we can take $s_{v+1}=s_v+1$ for each $v\in \bf N$. Consider
$\mbox{ }_{q_{v+1}}g$, then $\pi ^{q_{v+1}}_{q_v}(\mbox{
}_{q_{v+1}}g)=\mbox{ }_{q_v}g$. Each $\mbox{ }_{q_v}g\in Hom
(M_{q_v})$ is the element of the symmetric group $S_{b_v}$, where
$b_v$ is the cardinality of the finite set $\pi _{q_v}(M)=M_{q_v}$
(see also Theorem 3.6 above). Thus, $\mbox{ }_{q_v}g$ has a
decomposition into a product of nonintersecting finite cycles. If
$y_v\in M_{q_v}$, then the cardinality of each $(\pi
^{q_{v+1}}_{q_v})^{-1}(y_v)$ is the same for each $y_v\in M_{q_v}$.
The homomorphism $\mbox{ }_{q_{v+1}}g$ on each subset $(\pi
^{q_{v+1}}_{q_v})^{-1}(y_v)$ acts as the isometry. Then each cycle
of $\mbox{ }_{q_v}g$ splits into the product of cycles or becomes a
cycle of greater length for $\mbox{ }_{q_{v+1}}g$. Moreover, $\mbox{
}_{q_{v+1}}g$ is the product of two homomorphisms $h$ and $f$, where
$f((\pi ^{q_{v+1}}_{q_v})^{-1}(y_v))=(\pi
^{q_{v+1}}_{q_v})^{-1}(y_v)$ for each $y_v\in M_{q_v}$ and $\pi
^{q_{v+1}}_{q_v}(f)=\pi _{q_v}^*(id)$ is the identity on $\pi
_{q_v}(M)$, while $h(z_{v+1})=h(y_{v+1,0})+ (z_{v+1}-y_{v+1,0})$
with a marked $y_{v+1,0}\in (\pi ^{q_{v+1}}_{q_v})^{-1}(y_v)$ and
for each $z_{v+1}\in (\pi ^{q_{v+1}}_{q_v})^{-1}(y_v)$. \par If
$\sigma \in \Sigma _l$ and $\sigma $ is a cycle of length $u$, then
its algebraic order is $u$, that is, the cyclic group $\{ \sigma ^a:
a\in {\bf N} \} $ is of order $u$. On the step $v$ there is the
commutative local subgroup $\mbox{ }_{q_v}g^{x_v}=\eta _v(x_v)$ for
each $x_v \in \pi _{s_v}(B({\bf K},1,|\pi |^s))$. Therefore, on the
$v+1$-th step choose $q_{v+1}$ sufficiently large such that a number
of nonintersecting cycles $\sigma _{v,j}$ or their length are
sufficiently large in the sense that a commutative subgroup of $\pi
_{q_{v+1}}(W)$ generated by finite products of $\sigma
_{v,j}^{b_j}$, $b_j\in \bf N$, is of sufficiently large order that
to provide the algebraic homomorphism $\eta _{v+1}$ from $\pi
_{s_{v+1}}(B({\bf K},1,|\pi |^s))$ into $\pi _{q_{v+1}}(W)$
satisfying Conditions $(1-3)$ with $\pi _{q_{v+1}}(x_0)$ instead of
$x_0$ and $\pi ^{s_{v+1}}_{s_v}(\eta _{v+1})=\eta _v$.
\par In accordance with Lemma 2.5.10 \cite{eng} if $ \{ \phi , f_{j'}
\} $ is a mapping of an inverse system $S= \{ X_j, \pi ^j_i, \Psi \}
$ into an inverse system $S' = \{ Y_{j'}, \pi ^{j'}_{i'}, \Psi ' \}
$ and all $f_{j'}$ are homeomorphisms, then the limit mapping $f =
pr-\lim \{ \phi , f_{j'} \} $ is also the homeomorphism of $X =
pr-\lim S$ onto $Y=pr-\lim S'$, where $\phi : \Psi '\to \Psi $ is a
nondecreasing function, $\Psi '$ and $\Psi $ are directed sets.
\par Using the same notation which can not cause a confusion we
mention the following. Each this induction step gives the
corresponding solution of $(4-6)$ with $\eta _v(\pi _{s_v}(x))\in
\pi _{q_v}(B({\bf K},1,|\pi |^s))$ and $\pi _{q_v}(y)$ and $\pi
_{q_v}(c_{j,k,{\bar m}})$, $\pi _{q_v}^*(Q_{\bar m}(y))$ and $\pi
_{s_v}(\chi _{j,k}(x))\in \pi _{q_v}(B({\bf K},1,|\pi |^s))$ instead
of $x, y, c_{j,k,{\bar m}}$, $Q_{\bar m}(y)$ and $\chi _{j,k}(x)$
respectively. If $f\in C^{\alpha }_b(M,{\bf K}^n)$ and $\| f-id \|
_{C^{\alpha }_b(M,{\bf K}^n)}< 1$, then $f\in Diff^{\alpha }(M)$. We
say that a subset $A$ in $M$ is an $\epsilon $-net, if for each
$y\in M$ there exists $z\in A$ such that $|z-y|<\epsilon $, where
$0<\epsilon <\infty $. We have $\lim_{v\to \infty }|\pi |^v=0$ and
representatives of $M_{q_v}$ in $M$ form a $\rho ^v$ net for
suitable subsequence $q_v$, which is denoted by the same notation,
where $0<\rho <1$ is a constant and $\rho $ is a parameter
characterizing a net $A_{q_v}$ corresponding to $\{ Q_{\bar m}:
{\bar m}\in {\bf N_0}^n \} $ for $m_j\le q_v$ for each $j=1,...,n$,
where each $z\in A_{q_v}$ is a zero of $Q_{\bar m}$ as soon as
$m_j\ge m_{j,0}$ for some $j=1,...,n$, where ${\bar
m}_0:=(m_{1,0},...,m_{n,0})$ corresponds to $z$. On the other hand,
$|\chi _{j,k,m}(x)|\le 1$ on $B({\bf K},1,|\pi |^s)$ for each $j, k
,m$, moreover, $x_0=1+\theta ^s$, $|x_0|=1$. Then the corresponding
to it series converges, since $\mbox{ }_{q_v}g^{\eta _v(x_v)}\in \pi
_{q_v}^*(W)$ for each $x_v\in \pi _{s_v}(B({\bf K},1,|\pi |^s))$ and
$W$ is complete relative to the $C^{\alpha }_b$ uniformity. Thus,
the inverse sequence of $\mbox{ }_{q_v}g^{\eta _v(\pi _{s_v}(x))}$
converges to a local multiplicative continuous one-parameter
subgroup $g^x(y)\in W$ relative to the $C^{\alpha }_b\times C^0_b$
uniformity which is $C^{\alpha }_b$ by $y\in M$ and $C^0_b$ by $x\in
B({\bf K},1,|\pi |^s)$.

\par {\bf 2. Theorem.} {\it If $M$ is a manifold on
the Banach space $X := c_0(\gamma _X,{\bf K})$ with a finite atlas
and charts with bounded $\phi _j(U_j)$ in $X$ for each $j\in \Lambda
_M$, where $char ({\bf K})=p>1$, then in $Diff^{\alpha }_b(M)$ in
each neighborhood of $id$ there are $g\ne id$ which does not belong
to any non-trivial one-parameter subgroup $g^y$ relative to the
additive group $({\bf K},+)$. Moreover, if $M$ is embedded into $X$
and \par $(i)$ $(h-h\circ g^{-1})$, $(h-h\circ g^{-1})\circ
g^2$,...,$(h-h\circ g^{-1})\circ g^{p-1}$ are $\bf K$-linearly
independent and $g\ne id$, where $h:=g-id$, then $g$ does not belong
to any one-parameter subgroup $\{ g^y: y\in ({\bf K},+) \} $.}
\par {\bf Proof.} From the conditions of this theorem it follows, that
$M$ has a $C^t$ embedding as an open bounded subset in $X$. In view
of Theorem 2.9 the group $Diff^{\alpha }_b(M)$ is metrizable. Since
$Diff^{\beta }_b(M)$ is everywhere dense in $Diff^{\alpha }_b(M)$
for each $\beta >\alpha $, then it is sufficient to consider the
case $t\ge 1$, $\alpha \in \{ t, [t] \} $. If $g\in Diff^{\alpha
}_b(M)$, then up to a diffeomorphism of manifolds $g|_U\in
Diff^{\alpha }_b(U)$ for any $U$ open in $M$, since $g: U\to
g(U)\subset M$. Let $U\subset M$ be a clopen bounded subset such
that $U\subset U_j$ for some $j\in \Lambda _M$, $g\in Diff^{\alpha
}_b(M)$ and $supp (g)\subset U$, then $g|_U\in Diff^{\alpha }_b(U)$.
Take in particular $U:=\phi _j^{-1}(B(X,0,|\pi |))$, where $\pi \in
\bf K$, $0<|\pi |<1$. Thus, if prove theorem for $U=B(X,0,|\pi |)$,
then it will be also true for $Diff^{\alpha }_b(M)$.
\par Consider in $Diff^{\alpha }_b(U)$ a left-invariant metric
$\rho ^{\alpha }$ (see Theorem 2.9).  Take $id\ne g\in W := \{ f\in
Diff^{\alpha }_b(U): \rho ^{\alpha }(f,e)\le |\pi | \} $, then
$g=id+h$, where $0< \| h \|_{C^{\alpha }_b(U,X)}\le |\pi |$. If we
consider $p^n$ as the element of $\bf K$, where $n\in {\bf N} := \{
1, 2,... \} $, then $p^n=0\in \bf K$, since $char ({\bf K})=p$. So
we need to have $g^{p^n}=g^0=id$ on $U$, where $g^1=g$,
$g^{k+1}=g^kg^1$,...,$g^{p^n}= g^{p^n-1}g^1$. Consider $(1+a)^k =
\sum_{s=0}^l {k\choose s}a^s$, where $a\in \bf K$, then
$(1+a)^{p^n}=1+a^{p^n}$ and inevitably $|(1+a)^{p^n}|=|1+a^{p^n}|$,
since each binomial coefficient ${{p^n}\choose l}$ is divisible on
$p$ for each $1\le l\le p^n-1$ and hence ${{p^n}\choose l}$ is equal
to zero in $\bf K$. Then $(x_i+ax_i^l)\circ (x_i+ax_i^l)=
x_i+ax_i^l+a\sum_{s=0}^l{l\choose s} x_i^{l+(l-1)s}a^s$ and so on.
In particular, $(x_i+ax_i^p)\circ (x_i+ax_i^p)=
x_i+2ax_i^p+x_i^{p^2}a^{p+1}$ for $l=p$ and so on. This shows, that
elements of the form $g(x)=id(x)+ax_i^le_i$ with $a\ne 0$ and
$1<l\in \bf N$ can not lie on any one-parameter subgroups, since
$g^p\ne id$.
\par Demonstrate in general for $M$ embedded into $X$, that each
$id\ne g\in W$ satisfying Condition $(i)$ does not belong to any
one-parameter subgroup. On the other hand, $g^2(x)=g\circ
g(x)=id\circ (id+h)(x) + h\circ (id+h)(x)=g(x) + h\circ g(x)$,
$g^n(x)=g^{n-1}(x)+h\circ g^{n-1}(x)$, consequently, $\rho ^{\alpha
}(g^n,g^{n-1})=\rho ^{\alpha }(g^2,g)\le |\pi |$ for each $n\in \bf
N$ and $\| h\circ g^{n-1} \| _{C^{\alpha }_b}\le |\pi |$. Then by
induction $g^n(x)=g(x)+ h\circ g(x) +... + h\circ g^{n-1}(x)$ for
each $n\ge 2$, hence $\rho ^{\alpha }(g^n,g)\le |\pi |$. Then $|
h(y)-h(x)|\le \| h \|_{C^{\alpha }_b(U,X)}|x-y|$ for each $x, y\in
U$, in particular, for $y=g(x)$, where $y-x=g(x)-x=h(x)$, hence $|
h\circ g(x)-h(x)|\le \| h \|_{C^{\alpha }_b(U,X)}|h(x)|\le |\pi
|^2$. Therefore, $|g^n(x)-g(x)-(n-1)h(x)|\le \| h \|_{C^{\alpha
}_b(U,X)}|h(x)| $ for each $x\in U$, where $|h(x)|\le \| h
\|_{C^{\alpha }_b(U,X)}|x|\le |\pi | |x|$, consequently,
$|g^n(x)-g(x)-(n-1)h(x)|\le (\| h \|_{C^{\alpha }_b(U,X)})^2|x|\le
|\pi |^2|x|$. Thus $|g^{p^k}(x)-id(x)| \le (\| h \|_{C^{\alpha
}_b(U,X)})^2|x|$, since $|p^kh(x)|=0$. Suppose that $g^{p^k}=id$ for
each $k\in \bf N$, then $0=g^{p}(x)-id(x)= g^p(x)-g(x)-(p-1)h(x)=
h(x)+h\circ g(x)+...+h\circ g^{p-1}(x)=(h\circ g(x)-h(x))+...+
(h\circ g^{p-1}(x)-h(x))=h(x)+ p(h\circ g^1(x)-h(x))+(p-1)(h\circ
g^2(x)-h\circ g(x))+...+ 3(h\circ g^{p-2}(x)-h\circ
g^{p-3}(x))+2(h\circ g^{p-1}(x)-h\circ g^{p-2}(x))- h\circ
g^{p-1}(x)$, since $ph(x)=0$ and $p(h\circ g-h)(x)=0$ identically
and $g^{p-1}(x)=g^{-1}(x)$. Therefore, it would be $0=|h(x)-h\circ
g^{-1}(x)|$ on $U$, since Condition $(i)$ is supposed to be
satisfied, that leads to the contradiction, since $h\ne 0$, $g\ne
id$, $h\circ g\ne h$. Thus $g^p(x)\ne id(x)$ for each $id\ne g\in W$
satisfying Condition $(i)$, consequently, $g$ does not belong to any
one-parameter subgroup $\{ g^y: y\in ({\bf K},+) \} $.
\par {\bf 2.1. Remark.} On the other hand, if $M\supset B({\bf
K}^n,x_0,r)$, then $g^y(x):=x+yz$ for $x\in B({\bf K}^n,x_0,r)$ and
$g^y(x):=x$ for $x\in M\setminus B({\bf K}^n,x_0,r)$ is the
non-trivial one-parameter subgroup for a marked $z\in B({\bf
K}^n,0,r)$, where $y\in B({\bf K},0,1)$, $0<r$, $x_0\in M$.
Therefore, for each neighborhood $W$ of $id$ there exists $0\ne z\in
B({\bf K}^n,0,r)$ such that $id\ne g\in W$.
\section{Topological perfectness of diffeomorphism groups.}
\par {\bf 1. Definition.} A group $G$ is called algebraically
perfect, if its commutator group $[G,G]$ coincides with $G$, where
$[G,G]$ is the minimal group generated by all commutators $[f,g] :=
f^{-1}g^{-1}fg$, $f, g\in G$. A group $G$ is called algebraically
simple, if it does not contain any normal subgroup other than $\{ e
\} $ and $G$. A topological group $G$ is called topologically
perfect if $cl_G[G,G]=G$, $G$ is called topologically simple if it
does not contain any closed normal subgroup other than $\{ e \} $ or
$G$, where $cl_GV$ denotes the closure of a subset $V$ in $G$.
\par {\bf 2. Remark.} It is well known that the symmetric group $S_n$
is perfect for $n\ne 2$ and $n\ne 6$, but it is not simple for $n\ge
3$, since the subgroup $A_n$ consisting of even permutations is its
normal subgroup different from $ \{ e \} $ and $S_n$ \cite{kargap}.
Below the topological perfectness and simplicity of the
diffeomorphism group $Diff^{\alpha }(M)$ is considered relative to
its $C^{\alpha }$ topology.
\par {\bf 3. Theorem.} {\it Let $M$ be a compact manifold, $t\ge 0$,
$\alpha \in \{ t, [t] \} $. Then $Diff^{\alpha }(M)$ supplied with
the $C^{\alpha }$ compact-open topology is topologically perfect.}
\par {\bf Proof.} In view of Theorem 3.6 above and Theorem 2.5
\cite{luseamb03} the diffeomorphism group $Diff^{\alpha }(M)$
algebraically is the projective limit of an inverse sequence $S= \{
G_q, \pi ^q_s, {\bf N} \} $ of finite groups $G_q$, where $G_q$ is
isomorphic with the symmetric group $S_b$ with $b=card (M_q)$,
$M_q=\pi _q(M)$. There exists $q_0\in \bf N$ such that for each
$q>q_0$ the cardinality of $M_q$ is greater than $6$, consequently,
$G_q$ is perfect for each $q>q_0$, since $S_b$ is perfect with
$b=card (M_q)$ by the H\"older Theorem 5.3.1 \cite{kargap}. Each
element $h\in Diff^t(M)$ is decomposable into the thread $ \{ h_q,
\pi ^q_s, {\bf N} \} $. Every $h_q$ is decomposable as a finite
product of commutators in $G_q$ for $q>q_0$. Therefore, the
commutator group $[G,G]$ is dense in $G$ (see also Theorem 3.7
above), since the projective limit of a thread of products of
commutators is the product of commutators, where $G=Diff^{\alpha
}(M)$.
\par {\bf 4. Corollary.} {\it Let $M$ be a locally compact manifold,
$t\ge 0$, $\alpha \in \{ t, [t] \} $. Then the diffeomorphism group
$Diff^{\alpha }(M)$ and the group $Diff^{\alpha }_c(M)$ of all
$C^{\alpha }$ diffeomorphisms with compact supports supplied with
the $C^{\alpha }$ compact-open topology are topologically simple.}
\par {\bf Proof.} The group $Diff^{\alpha }_c(M)$ is everywhere dense in
$Diff^{\alpha }(M)$. Therefore, it is sufficient to prove this
corollary for $Diff^{\alpha }_c(M)$.
\par The group $G=Diff^{\alpha }_c(M)$ satisfies
the Epstein system of axioms. Let $X$ be a paracompact Hausdorff
topological space, $G$ a group of homeomorphisms of $X$, and $\cal
U$ a basis of open sets for the topology of $X$. The Epstein axioms
are the following: \par $(E1)$ if $U\in \cal U$ and $g\in G$, then
$gU\in \cal U$;  \par $(E2)$ $G$ acts transitively on $\cal U$; \par
$(E3)$ let $g\in G$, $U\in \cal U$ and $\cal B$ be an open cover of
$X$, then there exist $m\in \bf N$ and $g_1,...,g_m\in G$ and
$V_1,...,V_m\in \cal B$ such that:
\par $(i)$ $g=g_mg_{m-1}...g_1$;
\par $(ii)$ $supp (g_i)\subset V_i$;
\par $(iii)$ $supp (g_i)\cup (g_{i-1}...g_1cl (U))\ne X$
for each $i: 1\le i\le n$.
\par The manifold $M$ has an $C^{\alpha }$ embedding $\theta $ as the
clopen subset into either ${\bf K}^n$ or a direct topological sum of
copies of ${\bf K}^n$ such that $\theta (M)$ is a disjoint union of
clopen balls, since $M$ is totally disconnected and has a base of
topology consisting of clopen compact subsets in $M$. As $\cal U$
take a system of all clopen proper subsets $U$ in $M$ for which
there exists $g\in G$ such that $\theta (gU)$ is a ball of finite
positive radius in ${\bf K}^n$. Since $M$ is totally disconnected
and modelled on ${\bf K}^n$ as the manifold, then $\cal U$ forms the
base of the topology of $M$. \par In particular, for the
diffeomorphism group $G=Diff^{\alpha }_c(M)$ and the family $\cal U$
of $M=X$ we have $gU\in \cal U$ for each $g\in G$, since
$g^{-1}(gU)=U$ and $g: M\to M$ is the continuous bijective
epimorphism together with $g^{-1}$. If $U=M$, then $g(M)=M$. A
disjoint clopen covering $\cal C$ of $\theta (M)$ by balls form a
$C^{\infty }$ atlas (moreover, it is the analytic atlas) such that
if $f|_U$ is of $C^{\alpha }$ class for each $U\in \cal C$, then $f$
is $C^{\alpha }$ on $\theta (M)$. If $U\in \cal U$ and $U\ne M$,
then $M\setminus U$ is a clopen nonvoid subset in $M$. Therefore, if
$U_1$ and $U_2$ are two clopen nonvoid subsets in $\cal U$ different
from $M$, then there exists $g\in G$ such that $g(U_1)=U_2$, since
each two clopen balls in ${\bf K}^n$ are analytically diffeomorphic.
Thus axioms $(E1,E2)$ are satisfied.
\par If $g\in G$, then $supp (g):= cl \{ x\in M: g(x)\ne x \} $
is compact and for an open covering $\cal B$ of $M$ there are
$V_1,...,V_m$ in $\cal B$ such that $supp (g)\subset V_1\cup ...
\cup V_m$. Suppose that $U\in \cal U$ and $cl (U)\ne M$, where $cl
(U)$ denotes the closure of $U$ in $M$. Since $U$ is clopen in $M$,
then $cl (U)=U$. The topological Tychonoff space $M$ is totally
disconnected and locally compact paracompact, since it is modelled
on ${\bf K}^n$. In view of Theorem 6.2.9 \cite{eng} it is strongly
zero dimensional. Therefore, $V_1,...,V_m$ has a refinement
$P_1,...,P_m$ consisting of clopen compact subsets in $M$ such that
$supp (g)\subset P_1\cup ... \cup P_m$. Take $P_0:=\emptyset $ and
$W_i := P_i\setminus _{j<i}P_j$, then $W_1,...,W_m$ is the disjoint
clopen covering of $supp (g)$, where $W_i\subset P_i\subset V_i$ for
each $i=1,...,m$. Therefore, $g(W_1)$,...,$g(W_m)$ is the clopen
disjoint covering of $g(supp (g))$. Each $\theta (W_i)$ is a finite
union of clopen balls $B_{i,j}$ in ${\bf K}^n$. Then put
$g_m|_{W_m}=g|_{W_m}$ and $g_m|_{M\setminus W_m}=id$,
$g_{m-1}|_{W_{m-1}}=g_m^{-1}g|_{W_{m-1}}$ and $g_{m-1}|_{M\setminus
W_{m-1}}=id$, $g_i|_{W_i}= g_{i+1}^{-1}...g_m^{-1}g|_{W_i}$ and
$g_i|_{M\setminus W_i}=id$ for each $i=1,...,m-2$, hence Conditions
$(i,ii)$ of $(E3)$ are satisfied. For each $1\le k<j<m$ there is the
identity
$g_j...g_k|_{W_k}=g_j...g_{k+1}g_{k+1}^{-1}...g_m^{-1}g|_{W_k}=
g_{j+1}^{-1}...g_m^{-1}g|_{W_k}$ and $g_m...g_k|_{W_k}=g|_{W_k}$ for
each $1\le k<m$. Then $g_{i-1}...g_1(U)\subset (U\setminus (W_1\cup
...\cup W_{i-1}))\cup (g_{i-1}...g_1(U\cap W_1)\cup
g_{i-1}...g_2(W_2\cap U)...\cup g_{i-1}(U\cap W_{i-1}))$ and
$W_i\cup g_{i-1}...g_1(U)\ne M$ for each $U\in \cal U$, since $U\ne
M$ and $g(U)\ne M$, consequently, Condition $(iii)$ of $(E3)$ is
also satisfied.
\par For each clopen compact $V$ in $M$ the group $G_V := \{
g\in G: supp (g)\subset V \} $ is topologically perfect by Theorem
3. Moreover, $G=\bigcup_V G_V$, where the union is by all clopen
compact subsets $V$ in $M$. The group $G$ is supplied with the
$C^{\alpha }$ compact-open topology. Therefore, $G$ is topologically
perfect. In view of the Epstein Theorem 2.2.1 and Corollary 2.2.2
\cite{banyaga} the commutator group $[G,G]$ is topologically simple,
since each $[G_V,G_V]$ is topologically simple. On the other hand,
$[G,G]$ is dense in $G$, hence $G$ is topologically simple. If $M$
is locally compact and noncompact, then $Diff^{\alpha }_c(M)$ is
everywhere dense in $Diff^{\alpha }(M)$, hence the latter group is
topologically simple relative to its $C^{\alpha }$ compact-open
topology.
\par {\bf 5. Remark.} If $M$ is a compact manifold, $t\ge 0$,
$\alpha \in \{ t, [t] \} $, then $Diff^{\alpha }(M)$ is not
algebraically simple being the projective limit of groups due to
Theorem 3.6.
\par {\bf 6. Theorem.} {\it Let $M$ be a $C^{\alpha }$
compact manifold, $t\ge 0$, $\alpha \in \{ t, [t] \} $. Then each
continuous automorphism of $Diff^{\alpha }(M)$ belongs to the group
of homeomorphisms $Hom (M,M) =: Hom (M)$.}
\par {\bf Proof.} Let $\psi $ be a continuous automorphism of
$Diff^{\alpha }(M)$. In view of Theorem 3.6 above or Theorem 2.5
\cite{luseamb03} for each $v\in \bf N$, $v\ge s$, there exists the
quotient mapping $\pi _v^* : Diff^{\alpha }(M)\to \Sigma _{b_v}$,
where $\Sigma _m$ is the symmetric group of the set $ \{ 1,...,m \}
$, $b_v$ is the cardinality of the finite set $M_v=\pi _v(M)$, $\pi
_v^*$ is induced by the quotient mapping $\pi _v: {\bf K}\to {\bf
K}/B({\bf K},0,|\pi |^v)$ with the help of polynomial expansions by
Formulas 3.2$(4,5)$. Consider $M$ embedded into ${\bf K}^n$. The set
of diffeomorphisms $f$ such that $f-id$ is a piecewise affine on
balls of the covering of $M$ is contained in $Diff^{\alpha }(M)$,
consequently, $\pi _v^*(Diff^{\alpha }(M))=S_{b_v}$ is the
epimorphism.  This also follows from Formulas 3.2$(4-7)$. The
automorphism $\psi $ induces the automorphism $\psi _v$ of $\pi
_v^*(Diff^{\alpha }(M))$ such that $\pi _v^*(\psi (g))=:\psi _v(\pi
_v^*(g))$, since then $\pi _v^*(\psi (gh))=\pi _v^*(\psi (g)\psi
(h))=\pi _v^*(\psi (g))\pi _v^*(\psi (h))=\psi _v(\pi _v^*(g))\psi
_v(\pi _v^*(h))$ for each $g, h\in Diff^{\alpha }(M)$. Consider $v$
sufficiently large such that $b_v>6$. In view of the H\"older
Theorem 5.3.1 \cite{kargap} the automorphism $\psi _v$ is internal:
$\psi _v(a)=h_vah_v^{-1}$ for each $a\in \pi _v^*(Diff^{\alpha
}(M))$, where $h_v\in \pi _v^*(Diff^{\alpha }(M))$ is a marked
element. This defines the inverse sequence $\{ h_v, \pi ^v_l,
\Lambda _s \} $ such that $\pi ^v_l(h_v)=h_l$ for each $v\ge l\in
\Lambda _s$, where $\pi ^v_l\circ \pi _v^*=\pi _l^*$ and $\pi ^v_l:
G_v\to G_l$ are algebraic epimorphisms of discrete groups $G_v:=\pi
_v^*(Diff^{\alpha }(M))$ for each $v\ge l\in \Lambda _s$, $\pi ^v_l$
are written without star for simplicity of notation. Each $h_v:
M_v\to M_v$ is the homeomorphism and each $\psi _v: G_v\to G_v$ is
the homeomorphism. A limit of an inverse mapping system of
homeomorphic mappings is a homeomorphism by Proposition 2.5.10
\cite{eng}. Therefore, $h=pr-\lim \{ h_v, \pi ^v_l, \Lambda _s \} $
is the element of $Diff_w(M)$ (see Theorems 3.6, 3.7 above and
Theorem 2.5 \cite{luseamb03}), moreover, $\psi (g)=\lim \{ h_v\pi
_v^*(g)h_v^{-1}, \pi ^v_l, \Lambda _s \} =hgh^{-1}$ for each $g\in
Diff^{\alpha }(M)$. On the other hand, $Diff_w(M)$ is algebraically
isomorphic with $Hom(M)$, hence $h\in Hom(M)$.
\par {\bf 7. Remark.} Theorem 6 is not true for a locally compact
noncompact manifold $M$ for the group $Diff_c^{\alpha }(M)$ of
compactly supported diffeomorphisms of $Diff^{\alpha }(M)$, since
then $Diff_c^{\alpha }(M)$ has the external automorphisms $\phi
_f(g):=fgf^{-1}$ for $f\in Diff^{\alpha }(M)\setminus Diff_c^{\alpha
}(M)$, where $g\in Diff_c^{\alpha }(M)$. Moreover, $Diff^{\alpha
}_c(M)$ is the proper normal subgroup in $Diff^{\alpha }(M)$, but
$Diff^{\alpha }_c(M)$ for a locally compact noncompact manifold $M$
is the proper non-closed subgroup everywhere dense in $Diff^{\alpha
}(M)$.

\section{Projective decomposition of loop groups.}
\par {\bf 1.} Let as in \S 2.1 $\bar M$ and $N$ be two compact manifolds
over a locally compact non-archimedean infinite field $\bf K$ with a
multiplicative non trivial norm relative to which $\bf K$ is
complete as the uniform space and $Diff^{\alpha }_0({\bar M})$ be a
subgroup in $Diff^{\alpha }({\bar M})$ of all elements $\psi \in
Diff^{\alpha }({\bar M})$ such that $\psi (s_0)=s_0$, where $s_0$ is
a marked point in $\bar M$, $\alpha \in \{ t, [t] \} $. Denote by
$C^{\alpha }_0(M,N)$ a subspace in $C^{\alpha }({\bar M},N)$ of all
elements $f \in C^{\alpha }({\bar M},N)$ such that $\lim_{|\zeta
_1|+...+|\zeta _n|\to 0}{\bar \Phi }^v(f-w_0) (s_0;h_1,...,h_n;\zeta
_1,...,\zeta _n)=0$ for $\alpha =t$ or $\lim_{|\zeta _1|+...+|\zeta
_n|\to 0}{\Upsilon }^v(f-w_0) (s_0^{[n]})=0$ for $\alpha =[t]$ for
each $v \in \{ 0,1,..., t \} $, where $M={\bar M}\setminus s_0$ and
$w_0({\bar M})= \{ y_0 \} $, $x^{[k+1]}=(x^{[k]},v^{[k]},\zeta _k)$
(see \S 2.6 \cite{luanmbp} and Section 2 above). Geometric loop
monoids $\Omega _{\alpha }(M,N)$ and loop groups $L_{\alpha }(M,N)$
for $C^{\alpha }$ classes of mappings were constructed in
\cite{luanmbp}. The same construction is for $char ({\bf K})=p>0$.
\par {\bf Theorem.} {\it Let $\Omega _{\alpha }(M,N)$ be a commutative
loop monoid, then the quotient mappings $\pi _k$ induce the
corresponding inverse sequence $\{ \Omega (M_k,N_k): k\in {\bf N} \}
$ such that $\Omega ^w(M,N):=pr-\lim_k \Omega (M_k,N_k)$ is the
commutative compact topological monoid, where $\pi _k^*: \Omega
_{\alpha }(M,N)\to \Omega (M_k,N_k)$, $\pi ^l_k: \Omega (M_l,N_l)\to
\Omega (M_k,N_k)$ are surjective mappings for each $l\ge k$, $\Omega
(M_k,N_k)=\{ f_k: f_k\in N_k^{M_k}, f_k(s_{0,k})=y_{0,k} \}
/K_{\alpha ,k}$, $K_{\alpha ,k}$ is an equivalence relation induced
by an equivalence relation $K_{\alpha }$. Moreover, $\Omega ^w(M,N)$
is a compactification of $\Omega _{\alpha }(M,N)$ relative to the
projective weak topology $\tau _w$.}
\par {\bf Proof.} In view of Corollary 3.3
$\pi _k(C^{\alpha }_0(M,N))$ is isomorphic with $\{ f_k: f_k\in
N_k^{M_k}, f_k(s_{0,k})=y_{0,k} \} $, where the quotient mapping is
denoted by $\pi _k$ both for $M$ and $N$, since it is induced by the
same ring homomorphism $\pi _k: {\bf K} \to {\bf K}/B({\bf K},0,1)$,
$s_{0,k}:=\pi _k(s_0)$ and $y_{0,k}:=\pi _k(y_0)$, where $k\ge
s=\max (s(M),s(N))$. Then $\pi _k^*(Diff^{\alpha }_0(M))$ is
isomorphic with $Hom_0(M_k):=\{ \psi _k: \psi _k\in Hom(M_k), \psi
_k(s_{0,k})=s_{0,k} \} $. All of this is also applicable with the
corresponding changes to classes of smoothness $C^{\alpha }$ (or
$C({\alpha })$ in the notation of \cite{luanmbp}), where $\alpha =t$
or $\alpha =[t]$ with substitution of ${\bar {\Phi }}^k$ on
$\Upsilon ^k$ in the latter case. If $f$ and $g$ are two $K_{\alpha
}$-equivalent elements in $C^{\alpha }_0(M,N)$, that is, there are
sequences $f_n$ and $g_n$ in $C^{\alpha }_0(M,N)$ converging to $f$
and $g$ respectively and also a sequence $\psi _n \in Diff^{\alpha
}_0(M)$ such that $f_n(x)=g_n(\psi _n(x))$ for each $x\in M$, then
$\pi _k^*(f_n)=:f_{n,k}$ and $g_{n,k}:=\pi _k^*(g_n)$ converge to
$\pi _k^*(f)$ and $\pi _k^*(g)$ respectively and also $\psi
_{n,k}:=\pi _k^*(\psi _n)\in Hom_0(M_k)$. From the equality
$f_{n,k}(x(k))=g_{n,k}(\psi _{n,k}(x(k)))$ for each $n\in \bf N$ and
$x(k)\in M_k$ it follows, that the equivalence relation $K_{\alpha
}$ induces the corresponding equivalence relation $K_{\alpha ,k}$ in
$\pi _k^*(C^{\alpha }_0(M,N))$ such that classes $<\pi
_k^*(f)>_{K,\alpha ,k}$ of $K_{\alpha ,k}$-equivalent elements are
closed. Each element $f_k\in \pi _k^*(C^{\alpha }_0(M,N))$ is
characterized by the equality $f_k(s_{0,k})=y_{0,k}$. This induces
the quotient mapping $\pi _k^*: \Omega _{\alpha }(M,N)\to \Omega
(M_k,N_k)$ and surjective mappings $\pi ^l_k: \Omega (M_l,N_l) \to
\Omega (M_k,N_k)$ for each $l\ge k$. Each $\Omega (M_k,N_k)$ is the
finite discrete set, since each $N_k^{M_k}$ is the finite discrete
set. This produces the inverse sequence of finite discrete spaces,
hence the limit $\Omega ^w(M,N) := pr-\lim \{ \Omega (M_k,N_k), \pi
^k_l, \Lambda _s \} $ of the inverse sequence is compact and totally
disconnected. It remains to verify that $\Omega ^w(M,N)$ is the
commutative topological monoid with the unit element and the
cancelation property.
\par From the equality $M={\bar M}\setminus \{ s_0 \} $, it follows
that $M_k={\bar M}_k,$ since for each $k\in \bf N$ there exists
$x\in M$ such that $x+B({\bf K^{\sf m}},0,p^{-k})\ni s_0$. Moreover,
$M_k$ and $N_k$ are finite discrete spaces. Then $\pi _k(M\vee
M)=M_k\vee M_k$, where $A\vee B:= A\times \{ b_0 \} \cup \{ a_0 \}
\times B\subset A\times B$ is the wedge product of pointed spaces
$(A,a_0)$ and $(B,b_0)$, $A$ and $B$ are sets with marked points
$a_0\in A$ and $b_0\in B$. The composition operation is defined on
threads $ \{ <f_k>_{K,\alpha ,k}: k \in {\bf N} \} $ of the inverse
sequence in the following way. There was fixed a $C^{[\infty
]}$-diffeomorphism $\chi : M\vee M\to M$ \cite{luanmbp}. Let $x\in
M$, then $\pi _k(x)\in M_k$ and $\chi ^{-1}(U)\in M\vee M$, where
$U:=\pi _k^{-1}(x+B({\bf K},0, p^{-k}))\cap M.$ On the other hand
$\chi ^{-1}(U)$ is a disjoint union of balls of radius $p^{-2k}$ in
$B({\bf K^{2m}},0,1)$, hence there is defined a surjective mapping
$\chi _k: M_{2k}\vee M_{2k} \to M_k$ induced by $\chi $, $\pi _k$
and $\pi _{2k}$ such that $\chi _k(\chi ^{-1}(U))=\pi _k(x)$. If $f$
and $g\in C^{\alpha }(M,N)$, then $f\vee g\in C^{\alpha }((M\vee
M),N)$ and $\chi (f\vee g)\in C^{\alpha }(M,N)$ as in \S 2.6
\cite{luanmbp}. Hence $\chi _k(f_{2k}\vee g_{2k})\in C^{\alpha
}(M_k,N_k)$ and inevitably $\chi _k(<f_{2k}\vee g_{2k}>_{K,\alpha
,2k}) =\chi _k(<f_{2k}>_{K,\alpha ,2k}\vee <g_{2k}>_{K,\alpha ,2k})
\in \Omega (M_k,N_k)$.
\par There exists a one to one correspondence between
elements $f\in C_w({\bar M},N)$ and $ \{ f_k: k \} \in \{ N_k^{M_k}:
k\in \Lambda _s \} $. Therefore, $pr-\lim_k \Omega (M_k,N_k)$
algebraically this is the commutative monoid with the cancelation
property. Let $U$ be a neighborhood of $e$ in $\Omega ^w(M,N)$, then
there exists $U_k=\pi _k^{-1}(V_k)$ such that $V_k$ is open in
$\Omega (M_k,N_k)$, $e\in U_k$ and $U_k\subset U$. On the other hand
there exists $U_{2k}=\pi _{2k}^{-1}(V_{2k})$ such that $V_{2k}$ is
open in $\Omega (M_{2k},N_{2k})$, $e\in U_{2k}$ and
$U_{2k}+U_{2k}\subset U_k$. Therefore, $(f+U_{2k})+(g+U_{2k})\subset
f+g+U_k\subset f+g+U$ for each $f, g\in \Omega ^w(M,N)$,
consequently, the composition in $\Omega ^w(M,N)$ is continuous.
Since $C^{\alpha }_0(M,N)$ is dense in $C_{0,w}({\bar M},N)$, then
$\Omega _{\alpha }(M,N)$ is dense in $\Omega ^w(M,N)$ relative to
the projective weak topology $\tau _w$.
\par {\bf 2. Corollary.} {\it The loop group $L_{\alpha }(M,N)$
has a non-archimedean compactification $L^w(M,N)$ relative to the
projective weak topology $\tau _w$.}
\par {\bf Proof.}
Using the Grothendieck construction we get a compactification
$L^w(M,N)={\bar F}/{\bar B}$ of a loop group $L_{\alpha }(M,N)$,
where $\bar F$ is a closure in $(\Omega ^w(M,N))^{\bf Z}$ of a free
commutative group $F$ generated by $\Omega ^w(M,N)$ and $\bar B$ is
a closure of a subgroup $B$ generated by all elements
$[a+b]-[a]-[b]$, since the product of compact spaces is compact by
the Tychonoff theorem.
\par {\bf 3.} Let now $s_0=0$ and $y_0=0$ be two marked points in
the compact manifolds $\bar M$ and $N$ embedded into $\bf K^{\sf m}$
and $\bf K^{\sf n}$ respectively. There is defined the following
$C^{[\infty ]}$-diffeomorphism $inv: ({\bf K^{\sf m}})' \to ({\bf
K^{\sf m}})'$ for $({\bf K^{\sf m}})':={\bf K^{\sf m}}\setminus \{
x:$ $\mbox{ there exists }$ $j \mbox{ with }$ $x_j=0 \} $ such that
$inv (x_1,...,x_{\sf m})= (x_1^{-1},...,x_{\sf m}^{-1})$. Let
$M'=M\cap ({\bf K^{\sf m}})'$, then $inv (M')$ is locally compact
and unbounded in $\bf K^{\sf m}$, consequently, $\pi _k (inv
(M'))=(inv (M'))_k$ is a discrete infinite subset in $\bf K_k^{\sf
m}$ for each $k\in \bf N$. Analogously $\pi _k (inv (M'\vee
M'))=(inv (M'\vee M'))_k \subset \bf K_k^{\sf 2m}$. There exists a
$C^{[\infty ]}$-diffeomorphism $\chi : M\vee M\to M$ such that
$inv\circ \chi \circ inv$ is the $C^{[\infty ]}$-diffeomorphism of
$inv (M'\vee M')$ with $inv (M')$ and it induces bijective mappings
$\chi _k$ of $inv ((inv (M'\vee M'))_k)$ with $inv ((inv (M'))_k)$
for each $k\in \bf N$ such that ${\hat \pi }^l_k\circ \chi _l=\chi
_k$ for each $l\ge k$, where ${\hat \pi }^l_k:=inv\circ \pi
^l_k\circ inv $. This produces inverse sequences of discrete spaces
$inv ((inv (M'))_k)=:{\hat M}_k$, $inv ((inv (M'\vee M'))_k)= {\hat
M}_k\vee {\hat M}_k$ and their bijections $\chi _k$ such that
$pr-\lim_k{\hat M}_k$ is homeomorphic with $M'$ and $pr-\lim_k \chi
_k$ is equal to $\chi $ up to the homeomorphism, since
$pr-\lim_k{\bf K_k^{\sf m}}=\bf K^{\sf m}$ (see also about
admissible modifications and polyhedral expansions in
\cite{luumpe,lufpmpe}). If $\psi \in Diff^{\alpha }_0({\bar M})$,
then ${\hat \psi }\in Diff^{\alpha }({\hat M})$. Let $J_{f,k}:= \{
h_k: h_k=f_k\circ \psi _k, \psi _k\in Hom({\hat M}_k), \psi
_k(s_{0,k})=s_{0,k} \} $ for $f_k\in N_k^{{\hat M}_k}$ with
$\lim_{x\to 0}f_k(x)=0$, then $J_{f,k}$ is closed and ${\hat \pi
}_k^*(<f>_{K,\alpha })\subset J_{f,k}$. Therefore, $g_k$ and $f_k$
are ${\hat K}_{\alpha ,k}$-equivalent if and only if there exists
$\psi _k\in Hom({\hat M}_k)$ such that $\psi _k(s_{0,k})=s_{0,k}$
and $g_k(x)=f_k(\psi _k(x))$ for each $x\in {\hat M}_k$. Let $\Omega
({\hat M}_k,N_k):={\hat \pi }_k^*(\Omega _{\alpha }(M,N))$.
\par {\bf Theorem.} {\it The set of $\Omega ({\hat M}_k,N_k)$
forms an inverse sequence \\ $S =\{ \Omega ({\hat M}_k,N_k); {\hat
\pi }^l_k; k\in \Lambda _s \} $ such that $pr-\lim S=: \Omega
^{i,w}(M,N)$ is an associative topological loop monoid with the
cancelation property and the unit element $e$. There exists an
embedding of $\Omega _{\alpha }(M,N)$ into $\Omega ^{i,w}(M,N)$ such
that $\Omega _{\alpha }(M,N)$ is dense in $\Omega ^{i,w}(M,N)$
relative to the projective weak topology $\tau _{i,w}$.}
\par {\bf Proof.} Let ${U'}_i$ be an analytic disjoint atlas
of $inv(M')$, $f\in C^{\alpha }(inv(M'),{\bf K})$, $\psi \in
Diff^{\alpha }(inv(M'))$, then each restriction $f|_{{U'}_i}$ has
the form $f|_{{U'}_i}(x)=\sum_mf_{i,m}{\bar Q}_{i,m}(x)$ for each
$x\in {U'}_i$, where ${\bar Q}_{i,m}$ are basic Amice polynomials
for ${U'}_i$, $f_{i,m}\in \bf K$. Therefore $f$ is a combination
$f=\nabla _if|_{{U'}_i}$, hence ${\hat \pi }_k^*(f\circ \psi (x))=
\sum_m[{\hat \pi }_k^*(f_{i,m})\nabla _{(i,\psi _k(x(k)) \in {\hat
\pi }_k({U'}_k)}{\bar Q}_{i,m,k}(\psi _k(x(k)))]$ and inevitably
${\hat \pi }_k^*((f\circ \psi )(x))=f_k\circ \psi _k(x(k))$, where
${\bar Q}_{i,m,k}:={\hat \pi }_k^*({\bar Q}_{i,m})$, $x\in inv(M')$
and $x(k)={\hat \pi }_k(x)$.
\par As in \S 2.6.2 \cite{luanmbp} we choose an infinite atlas
$At'(M):= \{ ({U'}_j,{\phi '}_j): j\in {\bf N} \} $ such that ${\phi
'}_j: {U'}_j\to B(X,{y'}_j,{r'}_j)$ are homeomorphisms, $\lim_{k\to
\infty }{r'}_{j(k)}=0$, $\lim_{k\to \infty } {y'}_{j(k)}=0$ for an
infinite sequence $\{ j(k)\in {\bf N}: k\in {\bf N} \} $ such that
$cl_{\bar M}[\bigcup_{k=1}^{\infty }{U'}_{j(k)}]$ is a clopen
neighborhood of zero in $\bar M$, where $cl_{\bar M}A$ denotes the
closure of a subset $A$ in $\bar M$. We take
$|{y'}_{j(k)}|>{r'}_{j(k)}$ for each $k$, hence $inv
(B(X,{y'}_j,{r'}_j)\cap X')=B(X,{y'}_j^{-1},{r'}_j^{-1})\cap X'$ and
$\bigcup_k inv({U'}_{j(k)}\cap X')$ is open in $X'$, where $X={\bf
K}^{\sf m}$. For an atlas $At'(M\vee M):=\{ (W_l,\alpha _l): l\in
{\bf N} \} $ with homeomorphisms $\alpha _l: W_l\to B(X,z_l,a_l)$,
$\lim_{k\to \infty } a_{l(k)}=0$, $\lim_{k\to \infty }z_{l(k)}=0$
for an infinite sequence $\{ l(k)\in {\bf N}: k\in {\bf N} \} $ such
that $cl_{\bar M\vee \bar M}[\bigcup_{k=1}^{\infty }W_{l(k)}]$ is a
clopen neighborhood of $0\times 0$ in $\bar M\vee \bar M$ we also
choose $|z_l|>a_l$ for each $l$, where $card ( {\bf N}\setminus \{
l(k): k\in {\bf N} \} ) =card ( {\bf N}\setminus \{ j(k): k\in {\bf
N} \} )$. \par Then we take $\chi (W_{l(k)})={U'}_{j(k)}$ for each
$k\in \bf N$ and $\chi (W_l)={U'}_{\kappa (l)}$ for each $l\in ({\bf
N}\setminus \{ l(k): k\in {\bf N} \} )$, where $\kappa : ({\bf
N}\setminus \{ l(k): k\in {\bf N} \} ) \to ({\bf N}\setminus \{
j(k): k\in {\bf N} \} )$ is a bijective mapping such that $p^{-1}\le
{r'}_{j(k)}/a_{l(k)}\le p$ for each $k$ and $p^{-1}\le {r'}_{\kappa
(l)}/a_l\le p$ for each $l\in ({\bf N}\setminus \{ l(k): k\in {\bf
N} \} )$. We can choose the locally affine mapping $\chi $ on
$M={\bar M}\setminus \{ s_0 \} $ such that ${\bar {\Phi }}^n\chi =0$
or $\Upsilon ^n\chi =0$ for each $n\ge 2$ and
$B(X',{y'}_l^{-1},{r'}_l^{-1})$ are diffeomorphic with $inv
({U'}_l\cap X')$ and $B(X'\vee X',z_l^{-1},a_l^{-1})$ are
diffeomorphic with $inv (W_l\cap (X'\vee X'))$.
\par This induces the diffeomorphisms ${\hat \chi }:=
inv\circ \chi \circ inv: {\hat M}\vee {\hat M}\to {\hat M}$ and
${\hat \chi }^*: C^{\alpha }_0(({\hat M} \vee {\hat M}, \infty
\times \infty ), (N,y_0))\to C^{\alpha }_0(({\hat M},\infty
),(N,y_0))$, since each ${\bar {\Phi }}^n(f\vee g)({\hat \chi
}^{-1})$ or $\Upsilon ^n(f\vee g)({\hat \chi }^{-1})$ has an
expression through ${\bar {\Phi }}^l(f\vee g)$ and ${\bar {\Phi
}}^j({\hat \chi }^{-1})$ or $\Upsilon ^l(f\vee g)$ and $\Upsilon
^j({\hat \chi }^{-1})$ respectively with $l, j\le q$ and $q$
subordinated to $\alpha $, where ${\hat M}:= inv (M')$ and
conditions defining the subspace $C^{\alpha }_0 (({\hat M},\infty
),(N,y_0))$ differ from that of $C^{\alpha }_0 ((M,s_0),(N,y_0))$ by
substitution of $\lim_{x\to s_0}$ on $\lim_{|x|\to \infty }$. Then
$\lim_{|x|\to \infty }|{\hat \chi }(x)|=\infty $, consequently,
there exists $k_0\in \bf N$ such that ${\hat \chi }_k: {\hat
M}_k\vee {\hat M}_k \to {\hat M}_k$ are bijections for each $k\ge
k_0$, where ${\hat \chi }_k:={\hat \pi }_k\circ {\hat \chi }$. If
$\psi \in Diff^{\alpha }(\bar M)$ and $\psi (0)=0$, then
$\lim_{|x|\to \infty }{\hat \psi }(x)=\infty $ and $\lim_{|x|\to
\infty }{\hat \psi }^{-1}(x)=\infty .$ Then considering ${\hat \psi
}_k$ we get an equivalence relation $K_{\alpha ,k}$ in $\{ f_k:
f_k\in N_k^{{\hat M}_k}, \lim_{|x|\to \infty } f_k(x)=0 \} $ induced
by $K_{\alpha }$, where ${\hat M}_k$ is supplied with the quotient
norm induced from the space $X$, since $X'\subset X$, $x\in {\hat
M}_k$. \par Let $J_k$ denotes the quotient mapping corresponding to
$K_{\alpha ,k}.$ Therefore analogously to \S 2.6 \cite{luanmbp} we
get, that $\Omega ({\hat M}_k,N_k)$ are commutative monoids with the
cancelation property and the unit elements $e_k$, since $\Omega
({\hat M}_k,N_k)= \{ f_k: f_k\in C^0({\hat M}_k,N_k), \lim_{|x|\to
\infty }f_k(x)=0 \} /{\hat K}_{\alpha ,k}$ and mappings ${\hat \pi
}^l_k: ({{\bf K^{\sf m}})'}_l\to ({{\bf K^{\sf m}})'}_k$ and
mappings $\pi ^l_k: {\bf K^{\sf n}}_l\to {\bf K^{\sf n}}_k$ induce
mappings ${\hat \pi }^l_k: \Omega ({\hat M}_l,N_l) \to \Omega ({\hat
M}_k,N_k)$ for each $l\ge k$. Let the topology in $\{ f_k: f_k\in
C^0({\hat M}_k,N_k), \lim_{|x|\to \infty }f_k(x)=0 \} $ be induced
from the Tychonoff product topology in $N_k^{{\hat M}_k}$ and
$\Omega ({\hat M}_k,N_k)$ be in the quotient topology. \par The
space $N_k^{{\hat M}_k}$ is metrizable by the Baire metric $\rho
(x,y):=|\pi |^{-j}$, where $j=\min \{ i: x_i\ne y_i,
x_1=y_1,...,x_{i-1}=y_{i-1} \} $, $x=(x_l: x_l\in N_k, l\in {\bf N}
)$, ${\hat M}_k$ is enumerated as $\bf N$, $\pi \in \bf K$, $0<|\pi
|<1$ is the generator of the valuation group $\Gamma _{\bf K}$.
Therefore, $\Omega ({\hat M}_k,N_k)$ is metrizable and the mapping
$(f_k,g_k)\to f_k\vee g_k$ is continuous, hence the mapping
$(J_k(f_k),J_k(g_k))\to J_k(f_k)\circ J_k(g_k)$ is also continuous.
Then $J_k(w_{0,k})$ is the unit element, where $w_{0,k}({\hat
M}_k)=0$. Hence $\Omega ^{i,w}(M,N)$ is the commutative monoid with
the cancelation property and the unit element. Certainly
$\prod_k\Omega ({\hat M}_k, N_k)$ is the topological monoid and
$pr-\lim S$ is a closed in it topological totally disconnected
monoid. For each $f\in C^{\alpha }_0(M,N)$ there exists an inverse
sequence $\{ f_k: f_k={\hat \pi }_k^*(f), k\in \Lambda _s \} $ such
that $f(x)=pr-\lim_kf_k(x(k))$ for each $x\in M'$, but $M'$ is dense
in $M$. Therefore there exists an embedding $\Omega ^{\alpha
}(M,N)\hookrightarrow \Omega ^{i,w}(M,N)$, hence $\Omega ^{\alpha
}(M,N)$ is dense in $\Omega ^{i,w}(M,N)$ relative to the projective
weak topology $\tau _{i,w}$, since $C^{\alpha }(M,N)$ is dense in
$C_w(M,N)$ relative to the $\tau _w$ topology.
\par {\bf 4. Corollary.} {\it The inverse sequence of loop monoids induces
the inverse sequence of loop groups $S_L:=\{ L({\hat M}_k,N_k);
{\hat \pi }^l_k; \Lambda _s \} $. Its projective limit
$L^{i,w}(M,N):= pr-\lim S_L$ is a commutative topological totally
disconnected group and $L_{\alpha }(M,N)$ has an embedding in it as
a dense subgroup.}
\par {\bf Proof.} Due to the Grothendieck construction the
inversion operation $f_k\mapsto f_k^{-1}$ is continuous in $L({\hat
M}_k,N_k)$ and homomorphisms ${\hat \pi }^l_k$ and ${\hat \pi }_k$
have continuous extensions from loop submonoids onto loop groups
$L({\hat M}_k,N_k)$. Each monoid $\Omega ({\hat M}_k,N_k)$ is
totally disconnected, since $N_k^{{\hat M}_k}$ is totally
disconnected and $\Omega ({\hat M}_k,N_k)$ is supplied with the
quotient ultrametric, hence the free Abelian group $F_k$ generated
by $\Omega ({\hat M}_k,N_k)$ is also totally disconnected and
ultramertizable, consequently, $L({\hat M}_k,N_k)$ is
ultrametrizable. Evidently their inverse limit is also
ultrametrizable and the equivalent ultrametric can be chosen with
values in ${\tilde \Gamma }_{\bf K}:=\{ |z|: z\in {\bf K} \} $,
where ${\tilde \Gamma }_{\bf K}\cap (0,\infty )$ is discrete in
$(0,\infty ):=\{ x: 0<x<\infty , x\in {\bf R} \} $. Then the
projective limit (that is, weak) topology of $L^{i,w}(M,N)$ is
induced by the projective weak topology of $C_w(M,{\bf K})$.
\par {\bf 5. Theorem.} {\it For each prime number $p$
the loop group $L_{\alpha }(M,N)$ in its weak topology inherited
from $L^{i,w}(M,N)$ has the non-archimedean compactification
isomorphic with ${\bf Z_p}^{\aleph _0},$ moreover, $L^{i,w}(M,N)$
has the compactification $(\nu {\bf Z})^{\aleph _0}$, where $\nu
{\bf Z}$ is the one-point Alexandroff compactification of $\bf Z$.}
\par {\bf Proof.} The projective ring homomorphism $\pi _k:
{\bf K}\to {\bf K_k}$ induces \par ${\hat \pi }_k^* ( {\bar \Phi
}^m(f(x;h_1,...,h_m; \zeta _1,...,\zeta _m))={\bar \Phi
}^mf_k(x(k);h_1(k),...,h_m(k); \zeta _1(k),...,\zeta _m(k))$\\ and
${\hat \pi }_k^* ( {\Upsilon }^m(f(x^{[m]}))={\Upsilon
}^{[m]}f_k(x(k)^{[m]})$, \\ where $m\in \bf N$, ${\bar \Phi }^mf_k$
and $\Upsilon ^mf_k$ are defined for the field of fractions
generated by ${\bf K_k}$, since $\bf K$ is the commutative field
(see also \cite{bacht} and \S \S 2.1-2.6 \cite{luanmbp}). Then the
condition
$$\lim_{|x|\to \infty }{\bar \Phi }^mf(x;h_1,...,h_m;
\zeta _1,...,\zeta _m)=0\mbox{ or }$$

$$\lim_{|x|\to \infty }{\Upsilon }^mf(x^{[m]})=0$$

implies the condition
$$\lim_{|x(k)|\to \infty }{\bar \Phi }^mf_k(x(k);h_1(k),...,h_m(k);
\zeta _1(k),...,\zeta _m(k))=0\mbox{ or }$$
$$\lim_{|x(k)|\to \infty }{\Upsilon }^mf_k(x(k)^{[m]})=0$$
respectively, where $x^{[1]}=(x,v^{[0]},\zeta _1)$, $x^{[m+1]}:=
(x^{[m]},v^{[m]},\zeta _{m+1})$. Therefore, $supp (f_k):={\hat
M}_k^f:= \{ x(k): f_k(x(k))\ne 0 \} $ is a finite subset of the
discrete space ${\hat M}_k$ for each $k\in \bf N$. Then evidently,
${\hat \pi }_k^*(<g>_{K,\alpha })$ is a closed subset in $N_k^{{\hat
M}_k}$ for each $g\in C^{\alpha }_0(({\hat M},\infty ),(N,0))$,
since for each limit point $f_k$ of ${\hat \pi }_k^*(<g>_{K,\alpha
})$ its support is the finite subset in ${\hat M}_k$. Let $k_0$ be
such that $N_{k_0}\ne \{ 0 \} $, then this is also true for each
$k\ge k_0$. If $f_k\notin {\hat \pi }_k^*(<w_0>_{K,\alpha })$ and
$k\ge k_0$, then $f_k^{\vee n}\notin {\hat \pi }_k^*(<w_0>_{K,\alpha
})$ for each $n\in \bf N$, where $f_k^{\vee n}:=f_k\vee ... \vee
f_k$ denotes the $n$-times wedge product, since $\| f^{\vee
n}\|_{C^{\alpha }} \ge \| f \|_{C^{\alpha }} >0$ and $\| f_k^{\vee
n}\|_{C({\bf K^{\sf m}_k},{\bf K^{\sf n}_k})} \ge \| f \|_{C({\bf
K^{\sf m}_k},{\bf K^{\sf n}_k})} >0,$ where $C({\bf K^{\sf
m}_k},{\bf K^{\sf n}_k})=\pi _k^* (C^{\alpha }_b({\bf K^{\sf
m}},{\bf K^{\sf n}}))$ is the quotient module over the ring $\bf
K_k$. Each ${\hat \pi }_k^*(<f>_{K,\alpha })$ can be presented as
the following composition $z_1b_1+...+z_lb_l$ in the additive group
$L({\hat M}_k,N_k)$, where each $b_i$ corresponds to ${\hat \pi
}_k^*(<g_i>_{K,\alpha })$ and the embedding of $\Omega ({\hat
M}_k,N_k)$ into $L({\hat M}_k,N_k)$, $z_i\in \{ -1, 0, 1 \} $,
$l=card ({\hat M}_k^f)$, ${\hat M}_k^{g_i}$ are singletons for each
$i=1,...,l$. \par Each ${\hat M}_k$ is the finite discrete set as
well as $N_k$. For each $x\ne y\in {\hat M}_k$ there exists $\psi
\in Hom_0(M_k)$ such that $\psi (x)=y$, where $0$ corresponds to
$s_0$ for convenience of the notation. Using the group $Hom_0(N_k)$
we get that $L({\hat M}_k,N_k)$ is isomorphic with ${\bf Z}^{n_k}$,
where $n_k=card(N_k)>1$. For each prime number $p>1$ there exists
the $p$-adic completion of $\bf Z$ which is ${\bf Z_p}$. In view of
Corollary 4 $L_{\alpha }(M,N)$ has the non-archimedean completion
isomorphic with ${\bf Z}_p^{\aleph _0}$, since $\bf Z$ is dense in
$\bf Z_p$ and $pr-\lim_k{\bf Z}^{n_k}= {\bf Z}^{\aleph _0}.$ \par On
the other hand, we can take the multiplicative subgroup $ \{ \theta
^l: l\in {\bf Z} \} $ of the locally compact field ${\bf
F_{p^u}}(\theta )$ which gives the embedding $\phi $ of $\bf Z$ into
${\bf F_{p^u}}(\theta )$, where $\theta ^0=1$. The completion of
$\phi ({\bf Z})$ in ${\bf F_{p^u}}(\theta )$ is $\phi ({\bf Z})\cup
\{ 0 \} $ which is the one-point Alexandroff compactification $\nu
\bf Z$ of ${\bf Z}$. This gives the non-archimedean completion $(\nu
{\bf Z})^{\aleph _0}$ of $L_{\alpha }(M,N)$. Moreover, ${\bf
Z}_p^{\aleph _0}$ and $(\nu {\bf Z})^{\aleph _0}$ are compact as
products of compact spaces and $L^{i,w}(M,N)$ has the aforementioned
embeddings into them.
\par {\bf 6. Note.}
Using quotient mappings $\eta _{p,s}: {\bf Z}\to {\bf Z}/p^s\bf Z$
we get that $L_{\alpha }(M,N)^{\aleph _0}$ has the compactification
equal to $\{ \prod_{p\in {\sf P}}{\bf Z_p}^{\aleph _0} \} \times
(\nu {\bf Z})^{\aleph _0}$ relative to the product Tychonoff
topology, where $\sf P$ denotes the set of all prime numbers $p>1$,
$s\in \bf N$. These compactifications produce characters of
$L_{i,w}(M,N)$, since each compact Abelian group has only
one-dimensional irreducible unitary representations \cite{hew}. On
the other hand, there are irreducible continuous representations of
compact groups in non-archimedean Banach spaces \cite{roosch}. Among
them there are infinite-dimensional \cite{diar2,robert}. Moreover,
in their initial $C^{\alpha }$ topologies diffeomorphism and loop
groups also have infinite-dimensional irreducible unitary
representations \cite{lutmf99,luanmbp}. At the same time topologies
of $L_{\alpha }(M,N)$ and $L^w(M,N)$ or $L^{i,w}(M,N)$ are
incomparable, since the topologies of $C^{\alpha }(M,N)$ and
$C_w(M,N)$ are incomparable (see Theorem 3.7 above). \par Projective
limits of groups obtained above have the non-archimedean origin
related with non-archimedean families of semi-norms on spaces of
continuous or more narrow classes of functions between manifolds
over ultra-normed fields. Generally, if a topological space $X$ has
a projective limit decomposition $X= pr-\lim \{ X_{\alpha }, \pi
^{\alpha }_{\beta }, \Lambda \} $, then if $f_{\beta }: X_{\beta
}\to Y$ is a continuous function into a topological space $Y$, then
$f := f_{\beta }\circ \pi _{\beta }: X\to Y$ is a continuous
function, where the mapping $\pi ^{\alpha }_{\beta }: X_{\alpha }\to
X_{\beta }$ is continuous for each $\alpha \ge \beta \in \Lambda $,
$\Lambda $ is a directed set, $\pi ^{\alpha }_{\beta }\circ \pi
_{\alpha } = \pi _{\beta }$, $\pi _{\alpha } : X\to X_{\alpha }$ is
continuous and epimorphic. Therefore, $f_{\alpha } = {\tilde \pi
}^{\beta }_{\alpha }(f_{\beta }) := f_{\beta }\circ \pi ^{\alpha
}_{\beta }$ for each $\alpha \ge \beta $ generate the inductive
limit $ind-\lim \{ C(X_{\beta },Y); {\tilde \pi }^{\beta }_{\alpha
}; \Lambda \} $, where $C(X,Y)$ denotes the family of all continuous
mappings from $X$ into $Y$. On the other hand, if $Y=pr-\lim \{
Y_{\gamma }, p^{\gamma }_{\delta }, \Psi \} $, then one gets
$pr-\lim \{ C(X,Y_{\gamma }); p^{\gamma }_{\delta }, \Psi \} $. Then
these two constructions can be combined with repeated application of
projective and inductive limits, which may be dependent on the order
of taking limits. If $card (\Psi )\ge \aleph _0$, then a suitable
box topology in $\prod_{\gamma \in \Psi } C(X,Y_{\gamma })$ is
strictly stronger than a weak topology in it and in its projective
limit subspace (see also \cite{nari}).
\par {\bf 7. Theorem.} {\it If $\bar M$ and $N$ are compact
manifolds, $\alpha \in \{ \infty , [\infty ] \} $, then $L_{\alpha }
(M,N)$ is the $C^{\alpha }$ Lie group.}
\par {\bf Proof.} The uniform space $C^{\alpha }((M,s_0),(N,y_0))$
has the structure of the $C^{\alpha }$ manifold, since $M$ and $N$
are $C^{\alpha }$ manifolds, where $\alpha \in \{ \infty , [\infty ]
\} $. Therefore, $\Omega _{\alpha }(M,N)$ and $L_{\alpha }(M,N)$ are
$C^{\alpha }$ manifolds. The wedge product $(f,g)\mapsto f\vee g$ in
$C^{\alpha }((M,s_0),(N,y_0))$ is the $C^{\alpha }$ mapping, since
$M={\bar M}\setminus \{ s_0 \} $. Using the quotient mapping by
closures of equivalence relation caused by the action of
$Diff^{\alpha }_0(M)$ becomes the $C^{\alpha }$ manifold and
$C^{\alpha }$ monoid with the $C^{\alpha }$ smooth composition.
Using the construction of $L_{\alpha }(M,N)$ we get, that $L_{\alpha
}(M,N)$ is the $C^{\alpha }$ Lie group (see also for more details
\cite{luanmbp}).
\par {\bf 8. Remark.} Theorem 7 can be generalized in the $C^{\alpha }_b$
class for noncompact $C^{\alpha }_b$ manifolds $M$ and $N$.
\section{Appendix.}
\par {\bf 1. Lemma.} {\it  Let either
$f, g\in C^{[n]}(U,Y)$, where $U$ is an open subset in $X$, $Y$ is
an algebra over $\bf K$, or $f\in C^{[n]}(U,{\bf K})$ and $g\in
C^{[n]}(U,Y)$, where $Y$ is a topological vector space over $\bf K$,
then
\par $(1)$ $(fg)^{[n]}(x^{[n]}) = ({\Upsilon }\otimes {\hat P} +
{\hat {\pi }}\otimes {\Upsilon })^n
.(f\otimes g)(x^{[n]})$ \\
and $(fg)^{[n]}\in C^0(U^{[n]},Y)$, where $({\hat {\pi
}}^kg)(x^{[k]}) := g\circ \pi ^0_1\circ \pi ^1_2\circ ...\circ \pi
^{k-1}_k(x^{[k]})$, ${\hat P}^ng := P_nP_{n-1}...P_1g$, $\pi
^{k-1}_k(x^{[k]}) := x^{[k-1]}$, $(A\otimes B).(f\otimes g):=
(Af)(Bg)$ for $A, B\in L(C^n(U,Y),C^m(U,Y))$, $m\le n$, $(A_1\otimes
B_1)...(A_k\otimes B_k).(f\otimes g):=(A_1...A_k\otimes
B_1...B_k).(f\otimes g):= (A_1...A_kf)(B_1...B_kg)$ for
corresponding operators, ${\Upsilon
}^nf := f^{[n]}$, $(P_kg)(x^{[k]}) := g(x^{[k-1]}+v^{[k-1]}t_k)$, \\
${\hat P}^k{\hat {\pi }}^{a_1}\Upsilon ^{b_1}...{\hat {\pi
}}^{a_l}\Upsilon ^{b_l}g = P_{k+s} ...P_{s+1}{\hat {\pi
}}^{a_1}\Upsilon ^{b_1}...{\hat {\pi }}^{a_l}\Upsilon ^{b_l}g$ with
$s=b_1+...+b_l-a_1-...-a_l\ge 0$, $a_1,...,a_l, b_1,...,b_l \in \{
0, 1, 2, 3,... \} $.}
\par {\bf Proof.} Let at first $n=1$, then
\par $(2)$ $(fg)^{[1]}(x^{[1]})=[(fg)(x+vt)-(fg)(x)]/t =
[(f(x+vt)-f(x))g(x+vt) + f(x)(g(x+vt)-g(x))]/t= ({\Upsilon
}^1f)(x^{[1]})(P_1g)(x^{[1]}) + ({\hat {\pi
}}^0_1f)(x^{[1]}){\Upsilon }^1g(x^{[1]})$, \\
since ${\hat {\pi }}^0_1(x^{[1]})=x$ and $P_1$ is the composition of
the projection ${\hat {\pi }}^0_1$ and the shift operator on $vt$.
Let now $n=2$, then applying Formula $(2)$ we get: \par $(3)$
$(fg)^{[2]}(x^{[2]})= ((fg)^{[1]}(x^{[1]}))^{[1]}(x^{[2]}) =
(\Upsilon ^1(f^{[1]}(x^{[1]})(x^{[2]}))g(x+(v^{[0]}+v^{[1]}_2t_2)
(t_1+v^{[1]}_3t_2) + v^{[1]}_1t_2) + f^{[1]}(x^{[1]})
g^{[1]}(x+v^{[0]}t_1,v^{[1]}_1+v^{[1]}_2 (t_1+v^{[1]}_3t_2),t_2) +
f^{[1]}(x,v^{[1]}_1,t_2) g^{[1]}(x^{[1]}+v^{[1]}_1t_2) + f(x)
g^{[2]}(x^{[2]})$, \\
where $v^{[k]} = (v^{[k]}_1, v^{[k]}_2, v^{[k]}_3)$ for each $k\ge
1$ and $v^{[0]}=v^{[0]}_1$ such that $x^{[k]}+v^{[k]}t_{k+1} =
(x^{[k]}+v^{[k]}_1t_{k+1}, v^{[k-1]}+ v^{[k]}_2t_{k+1},
t_k+v^{[k]}_3t_{k+1})$ for each $1\le k\in \bf Z$. For $n=3$ we get
\par $(4)$ $(fg)^{[3]}(x^{[3]})=[(\Upsilon ^3f)({\hat P}^3g) +
({\hat {\pi }}^1\Upsilon ^2f)(\Upsilon ^1{\hat P}^2g)+ (\Upsilon ^1
({\hat {\pi }}^1\Upsilon ^1f))({\hat P}^1{\Upsilon ^1}{\hat
P}^1g)$\\
$+({\hat {\pi }}^2\Upsilon ^1f)(\Upsilon ^2{\hat P}^1g)+ (\Upsilon
^2{\hat {\pi }}^1f)({\hat P}^2\Upsilon ^1g)+ ({\hat {\pi
}}^1{\Upsilon }^1{\hat {\pi }}^1f) (\Upsilon ^1{\hat P}^1\Upsilon
^1g)$ \\  $+ (\Upsilon ^1({\hat {\pi }}^2f))({\hat
P}^1\Upsilon ^2g)+ ({\hat {\pi }}^3f) (\Upsilon ^3g)](x^{[3]})$, \\
since by our definition ${\hat P}^k{\hat {\pi }}^{a_1}\Upsilon
^{b_1}...{\hat {\pi }}^{a_l}\Upsilon ^{b_l}g = P_{k+s}
...P_{s+1}{\hat {\pi }}^{a_1}\Upsilon ^{b_1}...{\hat {\pi
}}^{a_l}\Upsilon ^{b_l}g$ with $s=b_1+...+b_l-a_1-...-a_l\ge 0$,
$a_1,...,a_l, b_1,...,b_l \in \{ 0, 1, 2, 3,... \} $. \par
Therefore, Formula $(1)$ for $n=1$ and $n=2$ and $n=3$ is
demonstrated by Formulas $(2-4)$. If $f, g \in C^0(U^{[k]},Y)$, $a,
b\in \bf K$, then $(P_k(af+bg))(x^{[k]}) :=
(af+bg)(x^{[k-1]}+v^{[k-1]}t_k)=$ $af(x^{[k-1]}+v^{[k-1]}t_k)+
bg(x^{[k-1]}+v^{[k-1]}t_k)$, moreover, ${\hat {\pi
}}^k(af+bg)(x^{[k]}) = (af+bg)\circ \pi ^0_1\circ \pi ^1_2\circ
...\circ \pi ^{k-1}_k(x^{[k]})=(af+bg)(x)= a f(x)+ b g(x) = a{\hat
{\pi }}^kf(x^{[k]}) + b {\hat {\pi }}^kg(x^{[k]})$ for each
$x^{[k]}\in U^{[k]}$, hence ${\hat {\pi }}^k$ and $P_k$ and ${\hat
P}^k$ are $\bf K$-linear operators for each $k\in \bf N$. Suppose
that Formula $(1)$ is proved for $n=1,...,m$, then for $n=m+1$ it
follows by application of Formula $(2)$ to both sides of Formula
$(1)$ for $n=m$:
\par $(fg)^{m+1}(x^{[m+1]})=((fg)^{[m]}(x^{[m]}))^{[1]}(x^{[m+1]})=
(({\Upsilon }\otimes {\hat P} + {\hat {\pi }}\otimes {\Upsilon })^m
.(f\otimes g)(x^{[m]}))^{[1]}(x^{[m+1]})= ({\Upsilon }\otimes {\hat
P} + {\hat {\pi }}\otimes {\Upsilon })^{m+1} .(f\otimes
g)(x^{[m+1]})$, \\
since $x^{[m+1]}=(x^{[m]})^{[1]}$ and more generally
$x^{[m+k]}=(x^{[m]})^{[k]}$ for each nonnegative integers $m$ and
$k$ such that $\pi ^{k-1}_k(x^{[m+k]})=x^{[m+k-1]}$ for $k\ge 1$;
$\Upsilon ^k$, ${\hat P}^k$ and ${\hat {\pi }}$ are $\bf K$-linear
operators on corresponding spaces of functions (see above and Lemma
2.3) and \par $({\Upsilon }\otimes {\hat P} + {\hat {\pi }}\otimes
{\Upsilon })^{m+1} .(f\otimes g)(x^{[m+1]})=$
\\  $\sum_{a_1+...+a_{m+1}+b_1+...+b_{m+1}=m+1} ({\Upsilon }^{a_1}\otimes
{\hat P}^{a_1})$\\ $({\hat {\pi }}^{b_1}\otimes {\Upsilon
}^{b_1})... ({\Upsilon }^{a_{m+1}}\otimes {\hat P}^{a_{m+1}}) ({\hat
{\pi }}^{b_{m+1}}\otimes {\Upsilon }^{b_{m+1}}).(f\otimes
g)(x^{[m+1]}) $, \\
where $a_j$ and $b_j$ are nonnegative integers for each
$j=1,...,m+1$, $(A_1\otimes B_1)...(A_k\otimes B_k).(f\otimes g):=
(A_1...A_k\otimes B_1...B_k).(f\otimes
g):=(A_1...A_kf)(B_1...B_kg).$
\par {\bf 2. Note.} Consider the projection
\par $(1)$ $\psi _n: X^{m(n)}\times {\bf K}^{s(n)}\to
X^{l(n)}\times {\bf K}^n$, \\
where $m(n)=2m(n-1)$, $s(n)=2s(n-1)+1$, $l(n)=n+1$ for each $n\in
\bf N$ such that $m(0)=1$, $s(0)=0$, $m(n)=2^n$,
$s(n)=1+2+2^2+...+2^{n-1}=2^n-1$. Then $m(n)$, $s(n)$, $l(n)$ and
$n$ correspond to number of variables in $X$, $\bf K$ for $\Upsilon
^n$, in $X$ and $\bf K$ for ${\bar {\Phi }}^n$ respectively.
Therefore, $\psi (x^{[n]})=x^{(n)}$ and $\psi _n(U^{[n]})=U^{(n)}$
for each $n\in \bf N$ for suitable ordering of variables. Thus
${\bar {\Phi }}^nf(x^{(n)})= {\hat {\psi }}_n{\Upsilon
}^nf(x^{[n]})=f^{[n]}(x^{[n]})|_{W^{(n)}}$, where ${\hat {\psi
}}_ng(y) := g(\psi _n(y))$ for a function $g$ on a subset $V$ in
$X^{l(n)}\times {\bf K}^n$ for each $y\in \psi _n^{-1}(V)\subset
X^{m(n)}\times {\bf K}^{s(n)}$, $W^{(n)}=U^{(n)} \times 0$, $0\in
X^{m(n)-l(n)}\times {\bf K}^{s(n)-n}$ for the corresponding ordering
of variables.
\par {\bf 3. Corollary.} {\it Let either $f, g\in C^n(U,Y)$, where $U$
is an open subset in $X$, $Y$ is an algebra over $\bf K$, or $f\in
C^n(U,{\bf K})$ and $g\in C^n(U,Y)$, where $Y$ is a topological
vector space over $\bf K$, then
\par $(1)$ ${\bar {\Phi }}^n(fg)(x^{(n)}) =
({\bar {\Phi }}\otimes {\hat P} + {\hat {\pi }}\otimes {\bar {\Phi }
})^n .(f\otimes g)(x^{(n)})$ \\
and ${\bar {\Phi }}^n(fg)\in C^0(U^{(n)},Y)$. In more details:
\par $(2)$ ${\bar {\Phi }}^n(fg)(x^{(n)}) =
\sum_{0\le a, 0\le b, a+b=n}\sum_{j_1<...<j_a; s_1<...<s_b; \{
j_1,...,j_a \} \cup \{ s_1,...,s_b \} = \{ 1,...,n \} }$ \\
${\bar {\Phi }}^af(x;v_{j_1},...,v_{j_a};t_{j_1},...,t_{j_a}) {\bar
{\Phi }}^bg(x+v_{j_1}t_{j_1}+...+v_{j_a}t_{j_a};v_{s_1},...,v_{s_b};
t_{s_1},...,t_{s_b})$.}
\par {\bf Proof.} The operator ${\hat {\psi }}_n$ is $\bf K$-linear,
since ${\hat {\psi }}_n(af+bg)(y)=(af+bg)(\psi _n(y))=af(\psi
_n(y))+ bg(\psi _n(y))$ for each $a, b\in \bf K$ and functions $f,
g$ on a subset $V$ in $X^{l(n)}\times {\bf K}^n$ and each $y\in \psi
_n^{-1}(V)\subset X^{m(n)}\times {\bf K}^{s(n)}$. Mention that the
restrictions of ${\hat {\pi }}^{k-1}_k$ and $P_k$ on $W^{(k)}$ gives
$\pi ^{k-1}_k(x^{(k)}) := x^{(k-1)}$ and $(P_kg)(x^{(k)}) :=
g(x^{(k-1)}+v_kt_k)$ in the notation of \S 1.1. The application of
the operator ${\hat {\psi }}_n$ to both sides of Equation 1$(1)$
gives Equation $(1)$ of this corollary, since ${\hat {\psi
}}_n\Upsilon ^n = {\bar {\Phi }}^n$ for each nonnegative integer
$n$, where $\Upsilon ^0=I$ and ${\bar {\Phi }}^0=I$ and ${\hat {\psi
}}_0=I$ are the unit operators.
\par {\bf 4. Lemma.} {\it Let $f_1,...,f_k\in C^{[n]}(U,Y)$, where $U$
is an open subset in $X$, either $Y$ is an algebra over $\bf K$, or
$f_1,...,f_{k-1}\in C^{[n]}(U,{\bf K})$ and $f_k\in C^{[n]}(U,Y)$,
where $Y$ is a topological vector space over $\bf K$, then
\par $(1)$ $(f_1...f_k)^{[n]}(x^{[n]}) = [\sum_{\alpha =0}^{k-1}
{\hat {\pi }}^{\otimes \alpha }\otimes \Upsilon \otimes {\hat
P}^{\otimes (k-\alpha -1)}]^n.(f_1\otimes ...\otimes f_k)(x^{[n]})$ \\
and $(f_1...f_k)^{[n]}\in C^0(U^{[n]},Y)$, where \par ${\hat {\pi
}}^{\otimes \alpha }\otimes \Upsilon \otimes {\hat P}^{\otimes
(k-\alpha -1)}.(f_1\otimes ... \otimes f_k) := ({\hat {\pi
}}(f_1...f_{\alpha }))(\Upsilon f_{\alpha +1})({\hat P}(f_{\alpha
+2}...f_k))$, where ${\hat {\pi }}^0:=I$, ${\hat P}^0=I$ is the unit
operator, ${\hat {\pi }}f_0:=1$, ${\hat P}f_{k+1}:=1$ (see Lemma
1).}
\par {\bf Proof.} Consider at first $n=1$ and apply Formula 1$(1)$
by induction to appearing products of functions, then
\par $(2)$ $\Upsilon ^1(f_1...f_k)(x^{[1]})=[(\Upsilon
^1(f_1...f_{k-1}))(P_1f_k) + ({\hat {\pi
}}^1(f_1...f_{k-1}))(\Upsilon ^1f_k)](x^{[1]})=$\\  $[(\Upsilon
^1(f_1...f_{k-2}))(P_1f_{k-1})(P_1f_k) + ({\hat {\pi
}}^1(f_1...f_{k-2}))(\Upsilon ^1f_{k-1})(P_1f_k)$\\ $+ ({\hat {\pi
}}^1(f_1...f_{k-1}))(\Upsilon ^1f_k)](x^{[1]}) =...$\\
$=(\sum_{\alpha =0}^{k-1}({\hat {\pi }}^1)^{\otimes \alpha }\otimes
\Upsilon ^1\otimes P_1^{\otimes
(k-\alpha -1)}).(f_1\otimes ...\otimes f_k)$, \\
where $A^{\otimes \alpha }\otimes B\otimes C^{\otimes (k-\alpha
-1)}.(f_1\otimes ... \otimes f_k) := (A(f_1...f_{\alpha
}))(Bf_{\alpha +1})(C(f_{\alpha +2}...f_k))$ for operators $A, B$
and $C$ and each nonnegative integer $\alpha $, where $A^0:=I$,
$C^0=I$ is the unit operator, $Af_0:=1$, $Cf_{k+1}:=1$, in
particular, $A={\hat {\pi }}^1$, $B=\Upsilon ^1$, $C=P_1$. Thus,
acting by induction on both sides by $\Upsilon ^1$ from Formula
$(2)$ we get Formula $(1)$ of this lemma, since the product of $n$
terms $\Upsilon ^1...\Upsilon ^1$  is equal to $\Upsilon ^n$.
\par {\bf 5. Corollary.} {\it Let $f_1,...,f_k\in C^n(U,Y)$, where $U$
is an open subset in $X$, either $Y$ is an algebra over $\bf K$, or
$f_1,...,f_{k-1}\in C^n(U,{\bf K})$ and $f_k\in C^n(U,Y)$, where $Y$
is a topological vector space over $\bf K$, then
\par $(1)$ ${\bar {\Phi }}^n(f_1...f_k)(x^{(n)}) =
[\sum_{\alpha =0}^{k-1} {\hat {\pi }}^{\otimes \alpha }\otimes {\bar
{\Phi }} \otimes {\hat P}^{\otimes (k-\alpha -1)}]^n.
(f_1\otimes ...\otimes f_k)(x^{(n)})$ \\
and ${\bar {\Phi }}^n(f_1...f_k)\in C^0(U^{(n)},Y)$, where
\par ${\hat {\pi }}^{\otimes \alpha }\otimes {\bar {\Phi }}\otimes
{\hat P}^{\otimes (k-\alpha -1)}.(f_1\otimes ... \otimes f_k) :=
({\hat {\pi }}(f_1...f_{\alpha }))({\bar {\Phi }}f_{\alpha
+1})({\hat P}(f_{\alpha +2}...f_k))$ (see Lemma 3).}
\par {\bf Proof.} Applying operator ${\hat {\psi }}_n$ from Note 2
to both sides of Equation 4$(1)$ we get Formula $(1)$ of this
Corollary.
\par {\bf 6. Lemma.} {\it Let $u\in C^{[n]}({\bf K}^s,{\bf K}^m)$,
$u({\bf K}^s)\subset U$ and $f\in C^{[n]}(U,Y)$, where $U$ is an
open subset in ${\bf K}^m$, $s, m\in \bf N$, $Y$ is a $\bf K$-linear
space, then
\par $(1)$ $(f\circ
u)^{[n]}(x^{[n]})=
[\sum_{j_1=1}^m...\sum_{j_n=1}^{m(n)}(A_{j_n,v^{[n-1]},t_n}...
A_{j_1,v^{[0]},t_1} f\circ u) (\Upsilon ^1\circ
p_{j_n}{\hat S}_{j_{n-1}+1,v^{[n-2]}t_{n-1}}$ \\
$...{\hat S}_{j_1+1,v^{[0]}t_1}u^{n-1}) (P_n\Upsilon ^1\circ
p_{j_{n-1}}{\hat S}_{j_{n-2}+1,v^{[n-3]}t_{n-2}}... {\hat
S}_{j_1+1,v^{[0]}t_1} u^{n-2})...(P_n...P_2\Upsilon ^1\circ
p_{j_1}u)$ \\  $+  \sum_{j_1=1}^m...\sum_{j_{n-1}=1}^{m(n-1)} ({\hat
{\pi }}^1(A_{j_{n-1},v^{[n-2]},t_{n-1}}...A_{j_1,v^{[0]},t_1}f\circ
u) [\sum_{\alpha =0}^{n-2} {\hat {\pi }}^{\otimes \alpha }\otimes
{\Upsilon } \otimes {\hat P}^{\otimes (n-\alpha -2)}]$ \\ $
((\Upsilon ^1\circ p_{j_{n-1}}{\hat
S}_{j_{n-2}+1,v^{[n-3]}t_{n-2}}...{\hat
S}_{j_1+1,v^{[0]}t_1}u^{n-2})
\otimes ... \otimes (P_{n-1}...P_2\Upsilon ^1\circ p_{j_1}u))$ \\
$+ [\sum_{\alpha =0}^{n-2} {\hat {\pi }}^{\otimes \alpha }\otimes
{\Upsilon } \otimes {\hat P}^{\otimes (n-\alpha -2)}]
(\sum_{j_1=1}^m...\sum_{j_{n-2}=1}^{m(n-2)}({\hat {\pi
}}^1(A_{j_{n-2}, v^{[n-3]},t_{n-2}}...A_{j_1,v^{[0]},t_1}f\circ
u))\otimes [\sum_{\alpha =0}^{n-3} {\hat {\pi }}^{\otimes \alpha
}\otimes {\Upsilon } \otimes {\hat P}^{\otimes (n-\alpha -3)}]
((\Upsilon ^1\circ p_{j_{n-2}}{\hat
S}_{j_{n-3}+1,v^{[n-4]}t_{n-3}}... {\hat
S}_{j_1+1,v^{[0]}t_1}u^{n-3})\otimes ...\otimes
(P_{n-2}...P_2\Upsilon ^1\circ p_{j_1}u))+...$
\\ $+ [\sum_{\alpha =0}^2 {\hat {\pi }}^{\otimes \alpha }\otimes
{\Upsilon } \otimes {\hat P}^{\otimes (2-\alpha )}]^{n-3} \{
\sum_{j_1=1}^m\sum_{j_2=1}^{m(2)} ({\hat {\pi }}^1
A_{j_2,v^{[1]},t_2} A_{j_1,v^{[0]},t_1}f\circ u) ({\Upsilon
^1}\otimes {\hat P}^1 + {\hat {\pi }}^1\otimes {\Upsilon ^1})
((\Upsilon ^1\circ p_{j_2}{\hat S}_{j_1+1,v^{[0]}t_1}u)\otimes
(P_2\Upsilon ^1\circ
p_{j_1}u)) \} $ \\
$+ ({\Upsilon }\otimes {\hat P} + {\hat {\pi }}\otimes {\Upsilon
})^{n-2} \{ \sum_{j_1=1}^m({\hat {\pi }}^1A_{j_1,v^{[0]},t_1} f\circ
u)\otimes (\Upsilon ^2\circ p_{j_1}u) \} ](x^{[n]})$ \\
and $f\circ u\in C^0(({\bf K}^s)^{[n]},Y)$, where \par $S_{j,\tau
}u(y):=(u_1(y),...,u_{j-1}(y),u_j(y+\tau _{(s)}),u_{j+1}(y+ \tau
_{(s)}),...,u_m(y+\tau _{(s)}))$, $u=(u_1,...,u_m)$, $u_j\in \bf K$
for each $j=1,...,m$, $y\in {\bf K}^s$, $\tau =(\tau _1,...,\tau
_k)\in {\bf K}^k$, $k\ge s$, $\tau _{(s)}:=(\tau _1,...,\tau _s)$,
$p_j(x):=x_j$, $x=(x_1,...,x_m)$, $x_j\in \bf K$ for each
$j=1,...,m$, ${\hat S}_{j+1,\tau }g(u(y),\beta ):=g(S_{j+1,\tau
}u(y),\beta )$, $y\in {\bf K}^s$, $\beta $ is some parameter,
$A_{j,v,t}:=({\hat S}_{j+1,vt} \otimes t\Upsilon ^1\circ
p_j)^*\Upsilon ^1_j$, where $\Upsilon ^1$ is taken for variables
$(x,v,t)$ or corresponding to them after actions of preceding
operations as $\Upsilon ^k$, $\Upsilon ^1_jf(x,v_j,t) :=
[f(x+e_jv_jt)-f(x)]/t$, $(B\otimes A)^*\Upsilon ^1f_i\circ u^i
(x,v,t) := \Upsilon ^1_jf_i(Bu^i,v,Au^i)$, $B: {\bf K}^{m(i)}\to
{\bf K}^{m(i)}$, $A: {\bf K}^{m(i)}\to \bf K$, $e_j=
(0,...,0,1,0,...,0)\in {\bf K}^{m(i)}$ with $1$ on $j$-th place;
$m(i)=m+i-1$, $j_i=1,...,m(i)$; $u^1:=u$, $u^2:=(u^1,t_1\Upsilon
^1\circ p_{j_1}u^1)$,...,$u^n=(u^{n-1},t_{n-1}\Upsilon ^1\circ
p_{j_{n-1}}u^{n-1})$, $A_{j_1,v^{[0]},t_1}f\circ u =: f_1\circ u^1$,
$A_{j_n,v^{[n-1]},t_n}f_{n-1}\circ u^{n-1} =: f_n\circ u^n$, ${\hat
S}_*\Upsilon ^1f(z) := \Upsilon ^1f({\hat S}_*z)$.}
\par {\bf Proof.} At first consider $n=1$, then $(f\circ
u)^{[1]}(t_0,v,t) = [f(u(t_0+vt))-f(u(t_0))]/t$, where $t_0\in {\bf
K}^s$, $t\in \bf K$, $v\in {\bf K}^s$. Though we consider here the
general case mention, that in the particular case $s=1$ one has
$t_0\in \bf K$, $v\in \bf K$. Then \par $(f\circ u)^{[1]}(t_0,v,t)=
[f(u(t_0+vt)) - f(u_1(t_0), u_2(t_0+vt),...,u_m(t_0+vt))]/t +
[f(u_1(t_0),u_2(t_0+vt),u_3(t_0+vt),...,u_m(t_0+vt))-
f(u_1(t_0),u_2(t_0),u_3(t_0+vt),...,u_m(t_0+vt))]/t+$ \\  $...+
[f(u_1(t_0),...,u_{m-1}(t_0),u_m(t_0+vt)) - f(u(t_0))]/t$,\\ where
$u=(u_1,...,u_m)$, $u_j\in \bf K$ for each $j=1,...,m$. Since
$u_j(t_0+vt)-u_j(t_0) = tu_j^{[1]}(t_0,v,t)$, hence \par $(f\circ
u)^{[1]}(t_0,v,t) = \Upsilon
^1f((u_1(t_0),u_2(t_0+vt),...,u_m(t_0+vt)),e_1,t\Upsilon
^1u_1(t_0,v,t))$ \\  $\Upsilon ^1u_1(t_0,v,t) + \Upsilon
^1f((u_1(t_0),u_2(t_0),u_3(t_0+vt),...,u_m(t_0+vt)),e_2,t\Upsilon
^1u_2(t_0,v,t))$\\  $\Upsilon ^1u_2(t_0,v,t) +...+ \Upsilon
^1f(u(t_0),e_m,t\Upsilon ^1u_m(t_0,v,t))\Upsilon ^1u_m(t_0,v,t)$, \\
since $u_j\in \bf K$ for each $j=1,...,m$ and $\bf K$ is the field,
where $e_j=(0,...,0,1,0,...,0)\in {\bf K}^m$ with $1$ on $j$-th
place for each $j=1,...,m$. With the help of shift operators it is
possible to write the latter formula shorter:
\par $(2)$  $\Upsilon ^1(f\circ u)(y,v,t) = \sum_{j=1}^m
{\hat S}_{j+1,vt} \Upsilon ^1f(u(y),e_j,t\Upsilon ^1\circ
p_ju(y,v,t)) (\Upsilon ^1\circ p_ju(y,v,t))$, \\
where $p_j(x):=x_j$, $x=(x_1,...,x_m)$, $x_j\in \bf K$ for each
$j=1,...,m$, ${\hat S}_{j+1,\tau }g(u(y),\beta ):=g(S_{j+1,\tau
}u(y),\beta )$, $y\in {\bf K}^s$, $\tau \in {\bf K}^k$, $k\ge s$,
$\beta $ is some parameter. Introduce operators $A_{j,v,t}:=({\hat
S}_{j+1,vt} \otimes t\Upsilon ^1\circ p_j)^*\Upsilon ^1_j$, where
$\Upsilon ^1$ is taken for variables $(y,v,t)$ or corresponding to
them after actions of preceding operators as $\Upsilon ^k$
remembering that $y^{[k]}, v^{[k]}\in ({\bf K}^s)^{[k]}$, $t\in \bf
K$, $v^{[k]}=(v^{[k]}_1,v^{[k]}_2,v^{[k]}_3)$ with $v^{[k]}_1,
v^{[k]}_2\in ({\bf K}^s)^{[k-1]}$, $v^{[k]}_3\in {\bf K}^k$ for each
$k\ge 1$, in particular, $v^{[0]}=v^{[0]}_1$ for $k=0$, $\Upsilon
^1_jf(x,v,t) := [f(x+e_j v_jt)-f(x)]/t$, $(B\otimes A)^*\Upsilon
^1f_i\circ u^i (y,v,t) := \Upsilon ^1_jf_i(Bu^i,v,Au^i)$, $B: {\bf
K}^{m(i)}\to {\bf K}^{m(i)}$, $A: {\bf K}^{m(i)}\to \bf K$. For
example, in the particular case of $s=1$ we have $v^{[k]}\in ({\bf
K})^{[k]}$. Therefore, in the general case Formula $(2)$ takes the
form:
\par $(3)$ $\Upsilon ^1f\circ u(y,v,t)=\sum_{j=1}^m(A_{j,v,t}f\circ
u)(\Upsilon ^1\circ p_ju)(y,v,t)$.
\par Take now $n=2$, then
\par $\Upsilon ^2f\circ u(y^{[2]})= \Upsilon ^1 \sum_{j=1}^m
[(A_{j,v,t}f\circ u)(\Upsilon ^1\circ p_ju)(y,v,t)](y^{[2]})$. \\
In the square brackets there is the product, hence from Formula
1$(1)$ and Lemma 2.3 we get:
\par $(4)$ $\Upsilon ^2f\circ u(y^{[2]})=
\sum_{j=1}^m[(\Upsilon ^1A_{j,v^{[0]},t}f\circ u)(P_2\Upsilon
^1\circ p_ju) + ({\hat {\pi }}^1A_{j,v^{[0]},t}f\circ u)(\Upsilon
^2\circ p_ju)](y^{[2]})$. \\
Then from Formula $(3)$ applied to terms $A_{j,v,t}f\circ u$ it
follows, that $\Upsilon ^1A_{j_1,v^{[0]},t_1}f\circ u (y^{[2]})=
\sum_{j_2=1}^{m(2)}(A_{j_2,v^{[1]},t_2}A_{j_1,v^{[0]},t_1}f\circ u)
(\Upsilon ^1\circ p_{j_2} S_{j_1+1,v^{[0]}t_1}u)(y^{[2]})$, where
$v^{[0]}=v$, $t_1=t$ (see also Lemma 1). Therefore,
\par $(5)$ $\Upsilon ^2f\circ u(y^{[2]})=
[\sum_{j_1=1}^m \sum_{j_2=1}^{m(2)} (A_{j_2,v^{[1]},t_2}
A_{j_1,v^{[0]},t_1} f\circ u) (\Upsilon ^1 \circ p_{j_2}{\hat
S}_{j_1+1,v^{[0]}t_1}u)(P_2\Upsilon ^1\circ p_{j_1}u)
+\sum_{j_1=1}^m({\hat {\pi }}^1A_{j_1,v^{[0]},t_1}f\circ u)
(\Upsilon ^2 \circ p_{j_1}u)](y^{[2]})$. \\
Then for $n=3$ applying Formulas $(3)$ and 4$(1)$ to $(5)$ we get:
\par $(6)$ $\Upsilon ^3f\circ u(y^{[3]})=
[\sum_{j_1=1}^m\sum_{j_2=1}^{m(2)}\sum_{j_3=1}^{m(3)}
(A_{j_3,v^{[2]},t_3} A_{j_2,v^{[1]},t_2} A_{j_1,v^{[0]},t_1}f\circ
u)$\\ $(\Upsilon ^1\circ p_{j_3}{\hat S}_{j_2+1,v^{[1]}t_2} {\hat
S}_{j_1+1,v^{[0]}t_1}u^2) (P_2\Upsilon ^1\circ p_{j_2}{\hat
S}_{j_1+1,v^{[0]}t_1}u) (P_3P_2\Upsilon ^1\circ p_{j_1}u) +$ \\
$\sum_{j_1=1}^m\sum_{j_2=1}^{m(2)} [({\hat {\pi
}}^1(A_{j_2,v^{[1]},t_2} A_{j_1,v^{[0]},t_1}f\circ u)) (\Upsilon
^2\circ p_{j_2}{\hat S}_{j_1+1,v^{[0]}t_1}u) (P_3P_2\Upsilon ^1\circ
p_{j_1}u) +$\\  $({\hat {\pi }}^1\{ (A_{j_2,v^{[1]},t_2}
A_{j_1,v^{[0]},t_1} f\circ u) (\Upsilon ^1\circ p_{j_2}{\hat
S}_{j_1+1,v^{[0]}t_1}u) \} (\Upsilon ^1P_2\Upsilon ^1\circ
p_{j_1}u)] +$\\  $\sum_{j_1=1}^m\sum_{j_3=1}^{m(3)}
(A_{j_3,v^{[2]},t_3}{\hat {\pi }}^1A_{j_1,v^{[0]},t_1}f\circ u)
(\Upsilon ^1\circ p_{j_3}{\hat S}_{j_1+1,v^{[0]}t_1}u) (P_3\Upsilon
^2\circ p_{j_1}u)$ \\  $+ \sum_{j_1=1}^m ({\hat {\pi
}}^2A_{j_1,v^{[0]},t_1}f\circ u)
(\Upsilon ^3\circ p_{j_1}u)] (y^{[3]})$. \\
Thus Formula $(1)$ is proved for $n=1, 2, 3$. Suppose that it is
true for $k=1,...,n$ and prove it for $k=n+1$. Applying Formula
4$(1)$ to both sides of $(1)$ we get:
\par $(7)$ $\Upsilon ^{n+1}f\circ u (y^{[n+1]}) =
[\sum_{j_1=1}^m...\sum_{j_{n+1}=1}^{m(n+1)}
(A_{j_{n+1},v^{[n]},t_{n+1}}... A_{j_1,v^{[0]}, t_1} f\circ u)$ \\
$(\Upsilon ^1\circ p_{j_{n+1}}{\hat S}_{j_n+1,v^{[n-1]}t_n}...{\hat
S}_{j_1+1,v^{[0]}t_1}u^n) (P_{n+1}\Upsilon ^1\circ p_{j_n}{\hat
S}_{j_{n-1}+1,v^{[n-2]}t_{n-1}}... {\hat S}_{j_1+1,v^{[0]}t_1}
u^{n-1})...$\\  $(P_{n+1}...P_2\Upsilon ^1\circ p_{j_1}u) +
\sum_{j_1=1}^m...\sum_{j_n=1}^{m(n)} ({\hat {\pi
}}^1(A_{j_n,v^{[n-1]},t_n}...A_{j_1,v^{[0]},t_1}f\circ u)$\\
$\Upsilon ^1((\Upsilon ^1\circ p_{j_n}{\hat
S}_{j_{n-1}+1,v^{[n-2]}t_{n-1}} ...{\hat
S}_{j_1+1,v^{[0]}t_1}u^{n-1}) ... (P_n...P_2\Upsilon ^1\circ
p_{j_1}u))+$\\ $\Upsilon
^1(\sum_{j_1=1}^m...\sum_{j_{n-1}=1}^{m(n-1)}({\hat {\pi
}}^1(A_{j_{n-1}, v^{[n-2]},t_{n-1}}...A_{j_1,v^{[0]},t_1}f\circ
u))\Upsilon ^1 ((\Upsilon ^1\circ p_{j_{n-1}}{\hat
S}_{j_{n-2}+1,v^{[n-3]}t_{n-2}}$\\  $... {\hat
S}_{j_1+1,v^{[0]}t_1}u^{n-2})...(P_{n-1}...P_2\Upsilon ^1\circ
p_{j_1}u))+... + \Upsilon ^{n-2}\{ \sum_{j_1=1}^m
\sum_{j_2=1}^{m(2)}$\\  $({\hat {\pi }}^1 A_{j_2,v^{[1]},t_2}
A_{j_1,v^{[0]},t_1}f\circ u)\Upsilon ^1((\Upsilon ^1\circ
p_{j_2}{\hat S}_{j_1+1,v^{[0]}t_1}u)(P_2\Upsilon ^1\circ p_{j_1}u))
\} +$ \\  $\Upsilon ^{n-1} \{ \sum_{j_1=1}^m{\hat {\pi
}}^1A_{j_1,v^{[0]},t_1} f\circ u)(\Upsilon ^2\circ p_{j_1}u)
\} ](y^{[n+1]})= $ \\
$[\sum_{j_1=1}^m...\sum_{j_{n+1}=1}^{m(n+1)}
(A_{j_{n+1},v^{[n]},t_{n+1}}... A_{j_1,v^{[0]},t_1} f\circ u)
(\Upsilon ^1\circ
p_{j_{n+1}}{\hat S}_{j_n+1,v^{[n-1]}t_n}$ \\
$...{\hat S}_{j_1+1,v^{[0]}t_1}u^n) (P_{n+1}\Upsilon ^1\circ
p_{j_n}{\hat S}_{j_{n-1}+1,v^{[n-2]}t_{n-1}}... {\hat
S}_{j_1+1,v^{[0]}t_1} u^{n-1})...(P_{n+1}...P_2\Upsilon ^1\circ
p_{j_1}u)$ \\  $+ \sum_{j_1=1}^m...\sum_{j_n=1}^{m(n)} ({\hat {\pi
}}^1(A_{j_n,v^{[n-1]},t_n}...A_{j_1,v^{[0]},t_1}f\circ u)
[\sum_{\alpha =0}^{n-1} {\hat {\pi }}^{\otimes \alpha }\otimes
{\Upsilon } \otimes {\hat P}^{\otimes (n-\alpha -1)}]$ \\ $
((\Upsilon ^1\circ p_{j_n}{\hat
S}_{j_{n-1}+1,v^{[n-2]}t_{n-1}}...{\hat
S}_{j_1+1,v^{[0]}t_1}u^{n-1})\otimes ... \otimes(P_n...P_2\Upsilon
^1\circ p_{j_1}u))$ \\  $+ [\sum_{\alpha =0}^{n-1} {\hat {\pi
}}^{\otimes \alpha }\otimes {\Upsilon } \otimes {\hat P}^{\otimes
(n-\alpha -1)}] (\sum_{j_1=1}^m...\sum_{j_{n-1}=1}^{m(n-1)} ({\hat
{\pi }}^1(A_{j_{n-1}, v^{[n-2]},t_{n-1}}...A_{j_1,v^{[0]},t_1}f\circ
u))\otimes [\sum_{\alpha =0}^{n-2} {\hat {\pi }}^{\otimes \alpha
}\otimes {\Upsilon } \otimes {\hat P}^{\otimes (n-\alpha -2)}]$ \\
$((\Upsilon ^1\circ p_{j_{n-1}}{\hat
S}_{j_{n-2}+1,v^{[n-3]}t_{n-2}}... {\hat
S}_{j_1+1,v^{[0]}t_1}u^{n-2})\otimes ...\otimes
(P_{n-1}...P_2\Upsilon ^1\circ p_{j_1}u))$
\\ $+ [\sum_{\alpha =0}^2 {\hat {\pi }}^{\otimes \alpha }\otimes
{\Upsilon } \otimes {\hat P}^{\otimes (2-\alpha )}]^{n-2} \{
\sum_{j_1=1}^m\sum_{j_2=1}^{m(2)} ({\hat {\pi }}^1
A_{j_2,v^{[1]},t_2} A_{j_1,v^{[0]},t_1}f\circ u) ({\Upsilon
^1}\otimes {\hat P}^1 + {\hat {\pi }}^1\otimes {\Upsilon ^1})
((\Upsilon ^1\circ p_{j_2}{\hat S}_{j_1+1,v^{[0]}t_1}u)\otimes
(P_2\Upsilon ^1\circ p_{j_1}u)) \} $ \\
$+ ({\Upsilon }\otimes {\hat P} + {\hat {\pi }}\otimes {\Upsilon
})^{n-1} \{ \sum_{j_1=1}^m({\hat {\pi }}^1A_{j_1,v^{[0]},t_1} f\circ
u)\otimes (\Upsilon ^2\circ p_{j_1}u) \} ](y^{[n+1]}).$ \\
Mention that in general $(\Upsilon ^{n+1}f\circ u)(y^{[n+1]})$ may
depend nontrivially on all components of the vector $y^{[n+1]}$
through several terms in Formula $(7)$. Thus Formula $(1)$ of this
Lemma is proved by induction.
\par {\bf 7. Corollary.} {\it Let $u\in C^n({\bf K}^s,{\bf K}^m)$,
$u({\bf K}^s)\subset U$ and $f\in C^n(U,Y)$, where $U$ is an open
subset in ${\bf K}^m$, $s, m\in \bf N$, $Y$ is a $\bf K$-linear
space, then
\par $(1)$ ${\bar {\Phi }}^n(f\circ
u)(x^{(n)})= [\sum_{j_1=1}^m...\sum_{j_n=1}^{m(n)}
(B_{j_n,v^{(n-1)},t_n}... B_{j_1,v^{(0)},t_1} f\circ u)$\\ $({\bar
{\Phi }}^1\circ p_{j_n}{\hat S}_{j_{n-1}+1,v^{(n-2)}t_{n-1}}...{\hat
S}_{j_1+1,v^{(0)}t_1}u^{n-1}) (P_n{\bar {\Phi }}^1\circ
p_{j_{n-1}}{\hat S}_{j_{n-2}+1,v^{(n-3)}_0,t_{n-2}}... {\hat
S}_{j_1+1,v^{(0)}t_1}u^{n-2})$ \\  $...(P_n...P_2{\bar {\Phi
}}^1\circ p_{j_1}u) + \sum_{j_1=1}^m...\sum_{j_{n-1}=1}^{m(n-1)}
({\hat {\pi }}^1(B_{j_{n-1},v^{(n-2)},t_{n-1}}...B_{j_1,v^{(0)},t_1}
f\circ u) [\sum_{\alpha =0}^{n-2} {\hat {\pi }}^{\otimes \alpha
}\otimes {\bar {\Phi }} \otimes {\hat P}^{\otimes (n-\alpha -2)}]$
\\ $(({\bar {\Phi }} ^1\circ p_{j_{n-1}}{\hat
S}_{j_{n-2}+1,v^{(n-3)}t_{n-2}}...{\hat
S}_{j_1+1,v^{(0)}t_1}u^{n-2})
\otimes ... \otimes(P_{n-1}...P_2{\bar {\Phi }}^1\circ p_{j_1}u))$ \\
$+ [\sum_{\alpha =0}^{n-2} {\hat {\pi }}^{\otimes \alpha }\otimes
{\bar {\Phi }} \otimes {\hat P}^{\otimes (n-\alpha -2)}]
(\sum_{j_1=1}^m...\sum_{j_{n-2}=1}^{m(n-2)}({\hat {\pi
}}^1(B_{j_{n-2},v^{(n-3)}, t_{n-2}}...B_{j_1,v^{(0)},t_1}f\circ
u))\otimes [\sum_{\alpha =0}^{n-3} {\hat {\pi }}^{\otimes \alpha
}\otimes {\bar {\Phi }} \otimes {\hat P}^{\otimes (n-\alpha -3)}]
(({\bar {\Phi }} ^1\circ p_{j_{n-2}}{\hat
S}_{j_{n-3}+1,v^{(n-4)}t_{n-3}}... {\hat
S}_{j_1+1,v^{(0)}t_1}u^{n-3})\otimes ...\otimes (P_{n-2}...P_2{\bar
{\Phi }}^1\circ p_{j_1}u))+...$
\\ $+ [\sum_{\alpha =0}^2 {\hat {\pi }}^{\otimes \alpha }\otimes
{\bar {\Phi }} \otimes {\hat P}^{\otimes (2-\alpha )}]^{n-3} \{
\sum_{j_1=1}^m\sum_{j_2=1}^{m(2)} ({\hat {\pi }}^1
B_{j_2,v^{(1)},t_2} B_{j_1,v^{(0)},t_1}f\circ u) ({\bar {\Phi
}}^1\otimes {\hat P}^1 + {\hat {\pi }}^1\otimes {\bar {\Phi }}^1)
(({\bar {\Phi }}^1\circ p_{j_2}{\hat S}_{j_1+1,v^{(0)}t_1}u)\otimes
(P_2{\bar {\Phi }}^1\circ p_{j_1}u)) \} $ \\
$+ ({\bar {\Phi }}\otimes {\hat P} + {\hat {\pi }}\otimes {\bar
{\Phi }})^{n-2} \{ \sum_{j_1=1}^m({\hat {\pi }}^1B_{j_1,v^{(0)},t_1}
f\circ u)\otimes ({\bar {\Phi }}^2\circ p_{j_1}u) \} ](x^{(n)})$ \\
and $f\circ u\in C^0(({\bf K}^s)^{(n)},Y)$ (see notation of Lemma
9), where $B_{j,v,t}:=({\hat S}_{j+1,vt} \otimes t{\bar {\Phi
}}^1\circ p_j)^*{\bar {\Phi }}^1_j$, where ${\bar {\Phi }}^1$ is
taken for variables $(x,v,t)$ or corresponding to them after actions
of preceding operations as ${\bar {\Phi }}^k$, ${\bar {\Phi
}}^1_jf(x,v,t) := [f(x+e_jv_jt)-f(x)]/t$, $(B\otimes A)^*{\bar {\Phi
}}^1f_i\circ u^i (x,v,t) := {\bar {\Phi }}^1_jf_i(Bu^i,v,Au^i)$, $B:
{\bf K}^{m(i)}\to {\bf K}^{m(i)}$, $A: {\bf K}^{m(i)}\to \bf K$,
$m(i)=m+i-1$, $j_i=1,...,m(i)$, $u^1=u$, $u^2:=(u^1,t_1{\bar {\Phi
}}^1\circ p_{j_1}u^1)$, $u^n:=(u^{n-1},t_{n-1}{\bar {\Phi }}^1\circ
p_{j_{n-1}}u^{n-1})$, ${\hat S}_* {\bar {\Phi }}^1f(x):= {\bar {\Phi
}}^1f({\hat S}_*x)$.}
\par {\bf Proof.} The restriction of operators of Lemma 6
on $W^{(n)}$ from Note 2 gives Formula $(1)$ of this corollary,
where $v^{(k)}\in ({\bf K}^s)^k\times {\bf K}^k$.

\end{document}